\documentclass[a4paper,reqno,12pt]{amsart}
\usepackage[all]{xy}           
\usepackage{amssymb}           
\usepackage{hyperref}
\usepackage{eucal}
\usepackage{epsfig}
\usepackage{pst-grad} 
\usepackage{pst-plot} 
\numberwithin{equation}{section}

\newtheorem{definition}{Definition}[section]
\newtheorem{lemma}[definition]{Lemma}
\newtheorem{theorem}[definition]{Theorem}
\newtheorem{proposition}[definition]{Proposition}
\newtheorem{corollary}[definition]{Corollary}
\newtheorem{remarkth}[definition]{Remark}
\newtheorem{example}[definition]{Example}

\renewcommand{\emph}[1]{{\bfseries\itshape{#1}}}

\makeatletter
\newcommand\prol{\@ifstar{\@proldf}{\@prolpf}}  
\def\@prolpf{\@ifnextchar[{\@prolpf@wrt}{\@prolpf@}}
\def\@prolpf@wrt[#1]#2{\@ifnextchar[{\@prolpf@wrt@at{#1}{#2}}{\@prolpf@wrt@{#1}{#2}}}
\def\@prolpf@wrt@at#1#2[#3]{\prolsymbol^{#1}_{#3}#2}
\def\@prolpf@wrt@#1#2{\prolsymbol^{#1}#2}
\def\@prolpf@#1{\@ifnextchar[{\@prolpf@at{#1}}{\@prolpf@@{#1}}}
\def\@prolpf@at#1[#2]{\prolsymbol_{#2}#1}
\def\@prolpf@@#1{\prolsymbol#1}
\def\@proldf{\@ifnextchar[{\@proldf@wrt}{\@proldf@}}
\def\@proldf@wrt[#1]#2{\@ifnextchar[{\@proldf@wrt@at{#1}{#2}}{\@proldf@wrt@{#1}{#2}}}
\def\@proldf@wrt@at#1#2[#3]{\prolsymbol^{*#1}_{#3}#2}
\def\@proldf@wrt@#1#2{\prolsymbol^{*#1}#2}
\def\@proldf@#1{\@ifnextchar[{\@proldf@at{#1}}{\@proldf@@{#1}}}
\def\@proldf@at#1[#2]{\prolsymbol^*_{#2}#1}
\def\@proldf@@#1{\prolsymbol^*#1}
\def\prolsymbol{\mathcal{T}}
\makeatother

\begin{document}

\title[Cosserat media]{Lie groupoids and algebroids applied to the study of uniformity and homogeneity of Cosserat media}

\author[V. M. Jim\'enez]{V\'ictor Manuel Jim\'enez}
\address{V\'ictor Manuel Jim\'enez:
Instituto de Ciencias Matem\'aticas (CSIC-UAM-UC3M-UCM),
c$\backslash$ Nicol\'as Cabrera, n. 13-15, Campus Cantoblanco, UAM
28049 Madrid, Spain} \email{victor.jimenez@icmat.es}

\author[M. de Le\'on]{Manuel de Le\'on}
\address{Manuel de Le\'on: Instituto de Ciencias Matem\'aticas (CSIC-UAM-UC3M-UCM),
c$\backslash$ Nicol\'as Cabrera,n. 13-15, Campus Cantoblanco, UAM
28049 Madrid, Spain} \email{mdeleon@icmat.es}

\author[M. Epstein]{Marcelo Epstein}
\address{Marcelo Epstein:
Department of Mechanical Engineering. University of Calgary. 2500 University Drive NW, Calgary, Alberta, Canada, T2N IN4} \email{epstein@enme.ucalgary.ca}

\thanks{This work has been partially supported by MINECO Grants  MTM2013-42-870-P and the ICMAT Severo Ochoa projects SEV-2011-0087 and SEV-2015-0554.
V.M.~Jiménez wishes to thank MINECO for a FPI-PhD Position.}
 \subjclass[2000]{}

\begin{abstract}
A Lie groupoid, called \textit{second-order non-holonomic material Lie groupoid}, is associated in a natural way to any Cosserat media. This groupoid is used to give a new definition of homogeneity which does not depend on a reference crystal. The corresponding Lie algebroid, called \textit{second-order non-holonomic material Lie algebroid}, is used to characterize the homogeneity property of the material. We also relate these results with the previously ones in terms of non-holonomic second-order $\overline{G}$-structures.
\end{abstract}

\maketitle

\tableofcontents

\section{Introduction}

\indent{In the theory of Continuum Mechanics a \textit{body} is presented by a three-dimensional manifold $\mathcal{B}$ which can be covered with just one chart (see \cite{JEMARS}). A \textit{configuration} $\psi$ is an embbeding of $\mathcal{B}$ in $\mathbb{R}^{3}$ and it is very usual to identify the body with one of the configurations $\psi_{0}$ which is called \textit{reference configuration}. A change of configuration $\psi \circ \psi_{0}^{-1}$ is said to be a \textit{deformation}.}\\
\indent{A problem of importance from the point of view of theoretical and practical physics is the following: given a mechanical response as a function of the positions on the body, how to decide if all the points of the body are made of the same material. W. Noll developed a geometric theory to deal with the properties of the material in \cite{WNOLL} (see also \cite{CTRUE,CCWAN,CCWANSEG}). As the Noll's and Wang's approach shows us, the use of $G-$structures has redefined the formulation and facilited the derivation of specific results (see for example \cite{MAREMDL2,MAREMDL3,MAREMDL4,MELZA}). In fact, the lack of integrability of the associated $G-$structure manifests the presence of inhomogeneities (such as dislocations).}\\
\indent{Thus, we could say that the theory of inhomogeneities of simple materials is well established in terms of differential geometric structures. However, there are many non-simple materials. In fact, materials like granular solids, rocks or bones cannot be modelled without extra kinematic variables \cite{GCAPR,ACERI}. The theory of \textit{generalized media} was introduced by Eug\`ene and Fran\c{c}ois Cosserat between 1905 and 1910 in \cite{COSSE}. The Cosserat's associated to each point of the body a family of vector (\textit{directors}). In a more mathematical way, a \textit{Cosserat continuum} can be described as a manifold of dimension $m$ and a family of vector $n$ fields on the manifold. Some of the developments of the theory can be found in Maugin \cite{GAMAU1,GAMAU3}, Kr$\ddot{\text{o}}$ner \cite{EKRON} or Eringen \cite{ACERI}.}\\
\indent{The geometrical structures which are necessary to develop a rigorous theory has been available for some time. Actually, the notion of director given by Cosserats is closely related with \textit{frame bundles}. In 1950 C. Ehresmann (see \cite{CELC11,CELC12,CELC13,CELC14,CELC15}) formalized the notion of \textit{principal bundles} and studied many frame bundles associated in a natural way to an arbitrary manifold: non-holonomic and holonomic frame bundles. Thus, we can intepret a Cosserat medium as a linear frame bundle $F \mathcal{B}$ of a manifold $\mathcal{B}$ which can be covered with just one chart (see \cite{MAREMDL}). Then, a configuration of $F\mathcal{B}$ is an embedding $\Psi : F \mathcal{B} \rightarrow F \mathbb{R}^{3}$ of principal bundles such that the induced Lie group morphism is the identity map. Again, we fix a configuration $\Psi_{0}$, as the reference configuration, and a deformation is a change of configurations.}\\
\noindent{The constitutive elastic law is now written as}
$$W = W \left( X, F \right),$$
\noindent{where $X$ is a point of the Cosserat medium and $F$ is the gradient
of a deformation at a point $X$. This constitutive equation permits us to associate to each point two points $X,Y$ of $\mathcal{B}$ the family of \textit{material isomorphisms from $X$ to $Y$}, i.e., the local principal bundle isomorphisms from $X$ to $Y$ that are invariants by $W$.}\\
\indent{This paper is devoted to study the geometrical description of Cosserat media using the notions of \textit{Lie groupoids} and \textit{Lie algebroids}. The notion of \textit{non-holonomic material groupoid of second order} associated to a material body $\mathcal{B}$ arises in a very natural way. Actually, the collection of all the $1-$jets of material isomorphisms constitutes a groupoid over the body $\mathcal{B}$. This approach considering $1-$jet prolongations has been developed in \cite{MAREMDL}. In that paper \cite{MAREMDL}, the authors have characterized the homogeneity of a Cosserat medium in terms of the integrability of the natural associated non-holonomic $\overline{G}-$structures of second order (see also \cite{MAREMDL2,MAREMDL3,MAREMDL4,MAREMDL5,MAREMDL6,MAREMDL7}).}\\
\indent{Thus, as in the case of simple materials (see \cite{VMJIMM}), a Cosserat medium has associated in a natural way a groupoid, called second-order non-holonomic material groupoid. The goal of our research is to introduce a new definition of homogeneity, based on the second-order non-holonomic material groupoid, and to characterize some materials properties like uniformity and homogeneity in terms of the Lie algebroid associated to the second-order non-holonomic material groupoid. An interesting feature is that our Lie algebroid can be described as the space of derivations on the tangent bundle of a manifold $M$ \cite{KMG}.}\\
\indent{Groupoids are an old topic coming back to the German mathematician Heinrich Brandt \cite{HBU}; Brandt was led to this concept by his studies in number theory. Of course, it is possible to find groupoids present implicitly in many early works, for instance the notion of continuous group of local transformations used by Sophus Lie leads to Lie groupoids no less to Lie groups \cite{AFUJI}. }\\
\indent{One should remark that a groupoid can be defined as an small \textit{category} in which all morphisms are invertible. However, the notion of category were introduced some years later (1945) by Eilenberg and MacLane \cite{SESM}.}\\
\indent{On the other hand, The concept of Lie groupoid was defined by Ehresmann in a collection of articles \cite{CELC,CELP,CES,CEC} where the author added topological and differential structures to the groupoid. Pradines redefined the notion of Lie groupoid in \cite{JPRA}. In this paper, as an infinitesimal version of Lie groupoids, Pradines defined the notion of Lie algebroid. Thus, Pradines extended the construction of the Lie algebra associated to a Lie group.  This construction allows to associate a Lie algebroid to any Lie groupoid, generating a functor between these categories. This functor was given by Pradines \cite{JPRATT} and is detailed by Mackenzie in \cite{KML} for the case of Lie algebroids with the same base and by Higgins and Mackenzie in \cite{PJHKM} for the case of Lie algebroids with different bases.}\\\\
\indent{In this way, arose the natural question of extending the three Lie's fundamental theorems (see \cite{JJDUI}):}\\\\
\textbf{Lie's first fundamental theorem}\\
Any integrable Lie algebra can be integrated to a simply connected Lie group.\\\\
\textbf{Lie's second fundamental theorem}\\
Any morphism between integrable Lie algebras can be integrated to a morphism of Lie groups.\\\\
\textbf{Lie's third fundamental theorem}\\
Any Lie algebra can be integrated to a Lie group.\\\\
\indent{It is possible to extend Lie's first fundamental theorem and Lie's second fundamental theorem. Finally, in order to generalize Lie's third fundamental theorem, in \cite{JPRALP} the next question appears: is any Lie algebroid integrable? i.e., is any Lie algebroid associated to a Lie groupoid? Pradines thought in \cite{JPRALP} that the answer was affirmative and, for a long time, it was believed that he was right. However, Almeida and Molino found a counterexample in \cite{RAPM}, i.e., there exist not integrable Lie algebroids. This result has a great importance because it resists with the one for Lie groups and Lie algebras. In \cite{MCRLF} the authors give necessary and sufficient conditions for the integrability of any Lie algebroid. We will need to use some results of integrabiliy to characterize the homogeneity in terms of the Lie algebroid associated to the second-order non-holonomic material groupoid.}\\
\indent{This paper is structured as follows: Section 2 is devoted to recall the main notions on principal bundles. As an important case, we present the notion of linear frame bundle $FM$ of a manifold $M$. Associated to this concept, we introduce the notion of \textit{the second-order non-holonomic frame bundle $\overline{F}^{2}M$ of a manifold} $M$. We also deal with \textit{the second-order non-holonomic $\overline{G}-$structures} and its associated notion of integrability.}\\
\indent{In Section 3 we describe a method to prolongate a pair of ordinary \textit{parallelisms} to obtain a \textit{non-holonomic parallelism of second order}. This method can be used to define a more general notion of integrability: \textit{integrable prolongation}.}\\
\indent{Groupoids and Lie groupoids are the matter of Section 4. We have developed a general introduction to this topic. We consider with more detail the Lie groupoid $\Pi^{1} \left( M , M \right)$ of $1-$jets of local diffeomorphisms of a manifold $M$. As a subgroupoid of $\Pi^{1} \left( FM , FM \right)$ we consider the Lie groupoid $J^{1} \left( FM \right)$ of local automorphisms of $FM$ over the identity map of $Gl \left( n, \mathbb{R}\right)$. Finally, the \textit{non-holonomic groupoid of second order $\tilde{J}^{1} \left( FM \right)$ over $M$} is presented as a quotient space. In Section 5 we introduce Lie algebroids and construct the Lie algebroid associated to an arbitrary Lie groupoid. The construction is similar to the corresponding one for Lie groups and Lie algebras. Many examples, as the \textit{$1-$jets algebroid $A \Pi^{1} \left( M , M \right)$ over $M$} and the \textit{non-holonomic algebroid of second order $A\tilde{J}^{1} \left( FM \right)$ over $M$}, are also exhibited. An alternative description of the Lie algebroid $A \Pi^{1} \left( M , M \right)$ associated to the Lie groupoid $\Pi^{1} \left( M , M \right)$ is given in Section 6 in terms of derivations of the tangent bundle $TM$. To do this rigurously we need the corresponding notion of exponential map that generalizes that for Lie groups.}\\
\indent{In Section 7 we construct the non-holonomic material groupoid of second order $\overline{\Omega} \left( \mathcal{B} \right)$ of a Cosserat media $\mathcal{B}$ as a reduced subgroupoid of $\tilde{J}^{1} \left(F \mathcal{B} \right)$. We use this object to characterize the uniformity and to define a new notion of homogeneity of the Cosserat medium $\mathcal{B}$.}\\
\indent{Section 8 deals with the notion of integrability of reduced subgroupoid of second-order non-holonomic groupoid by introducing the concept of \textit{standard flat subgroupoid of $\tilde{J}^{1} \left( \mathbb{R}^{n} \right)$}. We also give a precise construction of prolongation of sections to generalize the notion of integrability of reduced subgroupoids of $\tilde{J}^{1} \left(F M \right)$. A similar development is made for Lie algebroids. The results of Section 8 are applied to characterize the homogeneity of Cosserat media.}\\
\indent{Finally, in Section 10 we recall the characterization of homogeneity using the notion of $\overline{G}-$structures \cite{MAREMDL} and relate them with the current one in terms of his groupoid. In particular, the homogeneity defined in this paper is equivalent to the homogeneity in terms of non-holonomic $\overline{G}-$structures of second order for some reference crystal.}

\section{Principal bundles}
Firstly, we will present the notion of \textit{principal bundle} (see \cite{KONOM}) and, as particular cases, we will introduce the concepts of \textit{second-order non-holonomic frame bundle} and \textit{non-holonomic $\overline{G}-$structure of second order}. This will permit us to introduce the notion of \textit{integrability}.\\
\begin{definition}
\rm
Let $P$ be a manifold and $G$ be a Lie group which acts over $P$ by the right satisfying:

\begin{itemize}
\item[(i)] The action of $G$ is free, i.e., 
$$xg = x \Leftrightarrow g = e,$$
where $e \in G$ is the identity of $G$.
\item[(ii)] The canonical projection $\pi: P \rightarrow M = P/G$, where $P/G$ is the space of orbits, is a surjective submersion.
\item[(iii)] $P$ is locally trivial, i.e., $P$ is locally a product $U \times G$, where $U$ is an open set of $M$. More precisely, there exists a diffeomorphism $\Phi : \pi^{-1}\left(U\right) \rightarrow U \times G$, such that $\Phi \left(u\right) = \left(\pi \left(u\right) , \phi\left(u\right)\right)$, where the map $\phi : \pi^{-1}\left(U \right) \rightarrow G$ satisfies that
$$\phi \left(ua \right) = \phi \left( u \right)a, \ \forall u \in U, \ \forall a \in G.$$
$\Phi$ is called a \textit{trivialization on $U$}.

\end{itemize}
\end{definition}
A principal bundle will be denoted by $P\left(M,G\right)$, or simply $\pi: P \rightarrow M$ if there is no ambiguity about to the structure group $G$. $P$ is called the \textit{total space}, $M$ is the \textit{base space}, $G$ is the \textit{structure group} and $\pi$ is the \textit{projection}. The closed submanifold $\pi^{-1}\left(x\right)$, $x \in M$ will be called the \textit{fibre over $x$}. For each point $u \in P$, we have $\pi^{-1}\left(x\right) \triangleq  uG$, where $\pi \left(u\right) = x$, and $u G$ will be called the \textit{fibre through $u$}. Every fibre is diffeomorphic to $G$, but this diffeomorphism depends on the choice of the trivialization.\\

Now, we will define the morphisms of this category.
\begin{definition}
\rm
Given $P\left(M,G\right)$ and $P'\left(M',G'\right)$ principal bundles, a principal bundle morphism from $P\left(M,G\right)$ to $P'\left(M',G'\right)$ consists of a differentiable map $\Phi: P \rightarrow P'$ and a Lie group homomorphism $\varphi : G \rightarrow G'$ such that
$$\Phi \left(ua\right) = \Phi\left(u\right) \varphi\left(a\right), \ \forall u \in P, \ \forall a \in G.$$
\end{definition}
\noindent Notice that, in this case, $\Phi$ maps fibres into fibres and it induces a differentiable map $\phi :M \rightarrow M'$ by the equality $ \phi\left(x\right) = \pi\left(\Phi\left(u\right)\right)$, where $u \in \pi^{-1}\left(x\right)$.\\
If these maps are embeddings, the principal bundle morphism will be called \textit{embedding}. In such case, we can identify $P$ with $\Phi\left(P\right)$, $G$ with $\varphi \left(G\right)$ and $M$ with $\phi\left(M\right)$ and $P\left(M,G\right)$ is said to be a \textit{subbundle} of $P'\left(M',G'\right)$. Furthermore, if $M=M'$ and $\varphi = Id_{M}$, then $P\left(M,G\right)$ will be called a \textit{reduced subbundle} and we also say that $G'$ \textit{reduces} to the subgroup $G$.\\
Finally, a principal bundle morphism is called \textit{isomorphism} if it can be inverted by another principal  bundle morphism.\\

\begin{example}
\rm
Given a manifold $M$ and $G$ a Lie group, we can consider $M \times G$ as a principal bundle over $M$ with projection $pr_{1} : M \times G \rightarrow M$ and structure group $G$. The action is given by,
$$\left(x,a\right)b = \left(x,ab\right),  \ \forall x \in M , \ \forall a,b \in G.$$
This principal bundle is called a \textit{trivial principal bundle}.
\end{example}
Now, we will introduce an important example of principal bundle, the \textit{frame bundle}. In order to do that, we will start with the following definition.
\begin{definition}
\rm
Let $M$ be a manifold. A \textit{linear frame} at the point $x \in M$ is an ordered basis of $T_{x}M$.
\end{definition}

\begin{remarkth}
\rm
Alternatively, a linear frame at $x$ can be viewed as a linear isomorphism $z : \mathbb{R}^{n} \rightarrow T_{x}M$ identifying a basis on $T_{x}M$ as the image of the canonical basis of $\mathbb{R}^{n}$ by $z$.\\
There is a third way to interpret a linear frame by using the theory of jets. Indeed, a linear frame $z$ at $x \in M$ may be considered as the 1-jet $j^{1}_{0,x}\phi$ of a local diffeomorphism $\phi$ from an open neighbourhood of $0$ in $\mathbb{R}^{n}$ onto an open neighbourhood of $x$ in $M$ such that $\phi\left(0\right) = x$. So, $z= T_{0}\phi$.\\

\end{remarkth}
Thus, we denote by $FM$ the set of all linear frames at all the points of $M$. We can view $FM$ as a principal bundle over $M$ with the structure group $Gl\left(n , \mathbb{R}\right)$ and projection $\pi_{M} : FM \rightarrow M$ given by
$$ \pi_{M} \left(j^{1}_{0,x}\phi\right) = x, \ \forall j^{1}_{0,x}\phi \in FM.$$
This principal bundle is called the \textit{frame bundle on} $M$. Let $\left(x^{i}\right)$ be a local coordinate system on an open set $U \subseteq M$. Then we can introduce local coordinates $\left(x^{i}, x^{i}_{j}\right)$ over $FU \subseteq FM$ such that
\begin{equation}\label{27}
x^{i}_{j} \left( j_{0,x}^{1}\psi \right) =  \dfrac{\partial  \left(x^{i}\circ \psi\right)}{\partial x^{j}_{|0} }.
\end{equation}
If $\psi : N \rightarrow M$ is a local diffeomorphism, we denote by $F \psi : FN \rightarrow FM$ the local isomorphism induced from $\psi$, and defined by
$$F \psi \left( j_{0, \phi \left(0\right)}^{1} \phi \right) = j_{0, \psi \left( \phi \left(0\right)\right)}^{1}\left( \psi \circ \phi\right).$$
We denote by $e_{1}$ the frame $j_{0,0}^{1}Id_{\mathbb{R}^{n}} \in F \mathbb{R}^{n}$, where $Id_{\mathbb{R}^{n}}$ is the identity map on $\mathbb{R}^{n}$. Let $ \Psi : F \mathbb{R}^{n}  \rightarrow FM$ be a local isomorphism of principal bundles such that its domain contains 
$e_{1}$ and the induced isomorphism on Lie groups is the identity, i.e.,
$$\Psi \left( z \cdot g\right) = \Psi \left(z\right) \cdot g, \ \forall z \in Dom\left(\Psi\right) \subseteq F \mathbb{R}^{n}, \ \forall g \in Gl\left(n, \mathbb{R}^{n}\right).$$
We denote by $\psi : \mathbb{R}^{n} \rightarrow M$ the local diffeomorphism induced by $\Psi$. We recall that
$$ \psi \circ \pi_{\mathbb{R}^{n}} = \pi_{M} \circ \Psi.$$
The collection of all $1-$jets $j_{e_{1}, \Psi \left( e_{1}\right)}^{1} \Psi$ is a manifold which will be denoted by $\overline{F}^{2}M$. Of course, $j_{e_{1}, \Psi \left( e_{1}\right)}^{1} \Psi$ can be identified with a linear frame at the point $\Psi \left(e_{1}\right)$ since $T_{e_{1}} \Psi  : T_{e_{1}} \left( F \mathbb{R}^{n} \right) \cong \mathbb{R}^{n+n^{2}} \rightarrow T_{\Psi \left( e_{1}\right)}FM$ is a linear isomorphism, and we have $\overline{F}^{2}M \subset F\left(FM\right)$.\\
There are two canonical projections $\overline{\pi}^{2}_{1} : \overline{F}^{2}M \rightarrow FM$ and $\overline{\pi}^{2} : \overline{F}^{2}M \rightarrow M$ given by:
\begin{itemize}
\item $\overline{\pi}^{2}_{1} \left(j_{e_{1} , \Psi \left( e_{1} \right)}^{1} \Psi\right) = \Psi \left(e_{1}\right)$
\item $\overline{\pi}^{2} \left(j_{e_{1} , \Psi \left( e_{1} \right)}^{1} \Psi\right) = \psi \left(0\right)$
\end{itemize}
Of course, we have $ \overline{\pi}^{2} = \pi_{M} \circ \overline{\pi}^{2}_{1}$. We can show that $\overline{F}^{2}M$ is a principal bundle over $FM$ with canonical projection $\overline{\pi}^{2}_{1}$ and structure group,
$$ \overline{G}^{2}_{1}\left(n\right) := {\{ j_{e_{1} , e_{1}}^{1} \Psi \in \overline{F}^{2}\mathbb{R}^{n} / \ \Psi \left(e_{1}\right)= e_{1} \} = \overline{\pi}^{2}_{1}}^{-1} \left(e_{1}\right).$$
Notice that $\overline{G}^{2}_{1}\left(n\right)$ is a Lie subgroup of $Gl\left(n+n^{2} , \mathbb{R}\right)$ acting on $\overline{F}^{2}M$ by composition of jets.\\
We also have that $\overline{F}^{2}M$ is a principal bundle over $M$ with canonical projection $\overline{\pi}^{2}$ and structure group
$$ \overline{G}^{2}\left(n\right) := \{ j_{e_{1} , \Psi \left(e_{1}\right)}^{1} \Psi \in \overline{F}^{2}\mathbb{R}^{n} / \ \psi \left(0\right)= 0 \} = {\overline{\pi}^{2}}^{-1} \left(0\right),$$
which, again, acts on $\overline{F}^{2}M$ by composition of jets.\\
The principal bundle $\overline{F}^{2} M$ will be called the \textit{non-holonomic frame bundle of second order} and its elements will be called \textit{non-holonomic frames of second order}. There are more principal bundles defined over the $1-$jets of local isomorphisms $j^{1}_{e_{1} , Z} \Psi$ on $FM$ (\textit{holonomic} and \textit{semi-holonomic}). To know about the relations between them see \cite{MDELAM}\\
\begin{remarkth}
\rm
Notice that there exists a canonical projection $\tilde{\pi}^{2}_{1} : \overline{F}^{2}M \rightarrow FM$ defined by
$$ \tilde{\pi}^{2}_{1} \left(j^{1}_{e_{1}, Z} \Psi \right) = j_{0,z}^{1} \psi.$$
Observe that $\tilde{\pi}^{2}_{1}$ is a principal bundle morphism from $ \overline{F}^{2}M$ to $FM$ according to the diagram

\vspace{0.5cm}
\begin{picture}(375,50)(50,40)
\put(180,20){\makebox(0,0){$FM$}}
\put(250,25){$\pi_{M}$}               \put(200,20){\vector(1,0){90}}
\put(310,20){\makebox(0,0){$M$}}
\put(145,50){$\overline{\pi}^{2}_{1}$}                  \put(180,70){\vector(0,-1){40}}
\put(320,50){$\pi_{M}$}                  \put(310,70){\vector(0,-1){40}}
\put(180,80){\makebox(0,0){$\overline{F}^{2}M$}}
\put(250,85){$\tilde{\pi}^{2}_{1}$}               \put(200,80){\vector(1,0){90}}
\put(310,80){\makebox(0,0){$FM$}}
\end{picture}

\vspace{35pt}

\end{remarkth}

As we did with $FM$, having a local coordinate system $\left(x^{i}\right)$ on an open set $U \subseteq M$, we can introduce local coordinates $\left(x^{i}, x^{i}_{j}\right)$ over $FU \subseteq FM$ and, therefore, we can also introduce local coordinates \linebreak$\left(\left(x^{i}, x^{i}_{j}\right),x^{i}_{,j},x^{i}_{,jk},x^{i}_{j,k}, x^{i}_{j,kl}\right)$ over $F\left(FU\right)$ such that
\begin{itemize}
\item $x^{i}_{,j} \left( j_{e_{1},Z}^{1}\Psi \right) =  \dfrac{\partial \left(x^{i}\circ \Psi\right)}{\partial x^{j}_{|e_{1}} }$
\item $x^{i}_{,jk} \left( j_{e_{1},Z}^{1}\Psi \right) =  \dfrac{\partial \left(x^{i}\circ \Psi\right)}{\partial {x^{j}_{k}}_{|e_{1}} }$
\item $x^{i}_{j,k} \left( j_{e_{1},Z}^{1}\Psi \right) =  \dfrac{\partial \left(x^{i}_{j}\circ \Psi\right)}{\partial {x^{k}}_{|e_{1}} }$
\item $x^{i}_{j,kl} \left( j_{e_{1},Z}^{1}\Psi \right) =  \dfrac{\partial \left(x^{i}_{j}\circ \Psi\right)}{\partial {x^{k}_{l}}_{|e_{1}} }$
\end{itemize}
Hence, if we restrict to $\overline{F}^{2}U$ we have that
\begin{itemize}
\item $x^{i}_{,jk} =  0$
\item $x^{i}_{j,kl} =x^{i}_{k}\delta_{l}^{j}$
\end{itemize}
So, the induced coordinates on $\overline{F}^{2}U$ are 
\begin{equation}\label{138}
\left(\left(x^{i}, x^{i}_{j}\right),x^{i}_{,j},x^{i}_{j,k}\right).
\end{equation}
Then, locally
\begin{itemize}
\item $\pi_{M}\left(x^{i}, x^{i}_{j}\right) = x^{i}$
\item $\pi_{FM}\left(\left(x^{i}, x^{i}_{j}\right),x^{i}_{,j},x^{i}_{,jk},x^{i}_{j,k}, x^{i}_{j,kl}\right) = \left(x^{i}, x^{i}_{j}\right)$
\item $\overline{\pi}^{2}_{1} \left(\left(x^{i}, x^{i}_{j}\right),x^{i}_{,j},x^{i}_{j,k}\right) = \left(x^{i} , x^{i}_{j}\right)$
\item $\overline{\pi}^{2} \left(\left(x^{i}, x^{i}_{j}\right),x^{i}_{,j},x^{i}_{j,k}\right) = x^{i}$
\item $\tilde{\pi}^{2}_{1} \left(\left(x^{i}, x^{i}_{j}\right),x^{i}_{,j},x^{i}_{j,k}\right) = \left(x^{i}, x^{i}_{,j}\right)$
\end{itemize}

Notice that there exists a canonical isomorphism $\overline{F}^{2} \mathbb{R}^{n} \cong \mathbb{R}^{n} \times \overline{G}^{2}\left(n\right)$. In fact, let us define a global section $\overline{s} : \mathbb{R}^{n} \rightarrow \overline{F}^{2} \mathbb{R}^{n}$ as follows,
$$\overline{s} \left(x\right) = j_{e_{1} , {e_{1}}_{x}}^{1} F \tau_{x},$$
where $\tau_{x}$ denote the translation on $\mathbb{R}^{n}$ by the vector $x$. Then, locally
$$ \overline{s} \left(x^{i} \right) = \left(\left(x^{i} , \delta^{i}_{j}\right), \delta^{i}_{j}, 0\right).$$
So, a non-holonomic frame of second order $\overline{Z}$ at a point $x\in \mathbb{R}^{n}$ may be written in a unique way as
$$ \overline{Z} = \overline{s} \left(x\right) \cdot \overline{g},$$
where $\overline{g} \in \overline{G}^{2} \left(n\right)$. We have thus obtained a principal bundle isomorphism $\overline{F}^{2} \mathbb{R}^{n} \cong \mathbb{R}^{n} \times \overline{G}^{2} \left(n\right)$. Now, if $\overline{G}$ is a Lie subgroup of $\overline{G}^{2} \left(n\right)$, we can transport $\mathbb{R}^{n} \times \overline{G}$ by this isomorphism to obtain a $\overline{G}-$reduction of $\overline{F}^{2}\mathbb{R}^{n}$.\\

\begin{definition}
\rm
Let $\overline{G}$ be a Lie subgroup of $\overline{G}^{2}\left(n\right)$. A \textit{second-order non-holonomic $\overline{G}-$structure} $\overline{\omega}_{\overline{G}}\left(M\right)$ is a reduced subbundle of $\overline{F}^{2}M$ with structure group $\overline{G}$.
\end{definition}
Hence, the $\overline{G}-$reduction of $\overline{F}^{2}\mathbb{R}^{n}$ obtained above is a second-order non-holonomic $\overline{G}-$structure on $\mathbb{R}^{n}$ which will be called the \textit{standard flat second-order non-holonomic $\overline{G}-$structure}.\\

Next, we will introduce the notion of \textit{integrability of a second-order non-holonomic $\overline{G}-$structure}.\\

\begin{definition}\label{41}
\rm
Let $\overline{\omega}_{\overline{G}}\left(M\right)$ be a second-order non-holonomic \linebreak$\overline{G}-$structure on $M$. $\overline{\omega}_{\overline{G}}\left(M\right)$ is said to be \textit{integrable} if it is locally isomorphic to the trivial principal bundle $\mathbb{R}^{n} \times \overline{G}$, or equivalently, it is locally isomorphic to the standard flat $\overline{G}-$structure on $\mathbb{R}^{n}$.
\end{definition}

What we mean by "locally isomorphic" is that for each point $x \in M$, there exists a local chart through $x$, $\psi_{U}: U \rightarrow \overline{U}$ such that induces an isomorphism of principal bundles given by
$$ \Psi_{U} : \overline{\omega}_{\overline{G}} \left(U\right)  \rightarrow \overline{U} \times \overline{G},$$
where
$$ \Psi_{U} \left( j_{e_{1} , Z}^{1} \Psi \right) = \left( \psi_{U} \left(z\right), j_{e_{1} ,Z_{0}}^{1} \left(F\left(\tau_{-\psi_{U} \left(z\right)} \circ \psi_{U}\right) \circ \Psi\right) \right),$$
with $ \pi_{M}\left(Z\right)=z$.\\\\

\begin{remarkth}
\rm
There exists an alternative definition of second-order non-holonomic frames (see \cite{WNO}). Consider a differentiable map $\phi : U \rightarrow FM$ defined on some open neighbourhood of $0$ in $\mathbb{R}^{n}$ such that $\pi_{M} \circ \phi$ is an embedding. Then the $1-$jet $j_{0 , \phi \left(0\right)}^{1} \phi$ is a non-holonomic frame of second order at $x = \pi_{M}\left(\phi \left(0\right)\right)$. In fact, given $\phi$ we define a local principal bundle isomorphism $\Phi : F \mathbb{R}^{n} \rightarrow FM$ over $U$ given by
$$ \Phi \left(r , R\right) =\phi \left(r\right) R,$$
where $r \in \mathbb{R}^{n}$ and $ R \in Gl \left(n , \mathbb{R} \right)$. Thus, $j_{e_{1},Z}^{1}\Phi$ defines a non-holonomic frame of second order at $x$. Conversely, having a local principal bundle isomorphism $\Phi : F \mathbb{R}^{n} \rightarrow FM$ over an open set $U$, we define $\phi$ as follows:
$$ \phi \left(r\right) = \Phi \left(r , e\right) ,$$
where $r \in \mathbb{R}^{n}$ and $ e \in Gl \left(n , \mathbb{R} \right)$ is the identity.

\end{remarkth}

Any second-order non-holonomic $\{\overline{e}\}-$structure on $M$, with $\overline{e}$ the identity of $\overline{G}^{2}\left(n\right)$, will be called \textit{non-holonomic parallelism of second-order}. It is easy to show that any non-holonomic parallelism of second-order is, indeed, a global section of the second-order non-holonomic frame bundle $ \overline{F}^{2}M$. So, we can consider \textit{integrable sections}. Let $\left(x^{i}\right)$ be a local coordinate system on an open set $U \subseteq M$ and \linebreak$\left(\left(x^{i} , x^{i}_{j}\right), x^{i}_{,j} , x^{i}_{j,k}\right)$ be the induced coordinates on $\overline{F}^{2}U$ we have that (locally) any integrable sections can be written as follows

$$\overline{P}\left(x^{i}\right) = \left(\left(x^{i}, \delta^{i}_{j}\right) , \delta^{i}_{j}, 0\right).$$

Indeed, we can show that

\begin{proposition}\label{31}
Let $\overline{\omega}_{\overline{G}}\left(M\right)$ be a second-order non-holonomic \linebreak$\overline{G}-$structure on $M$. $\overline{\omega}_{\overline{G}}\left(M\right)$ is \textit{integrable} if, and only if, for each point $x \in M$ there exists a local coordinate system $\left(x^{i}\right)$ on $M$ such that the local section, 
\begin{equation}\label{32}
\overline{P}\left(x^{i}\right) = \left(\left(x^{i}, \delta^{i}_{j}\right), \delta^{i}_{j}, 0\right),
\end{equation}
takes values into $\overline{\omega}_{\overline{G}}\left(M\right)$.
\end{proposition}
Notice that, in a similar way to the case of the integrable \linebreak$G-$structures in the frame bundle (see \cite{VMJIMM}), Eq. (\ref{32}) is equivalent to the following: for each $z \in M$, there exists a local chart $\left(\psi_{U} , U \right)$ over $z $ such that for all $x \in U$ 
\begin{equation}\label{33}
\overline{P}\left(x\right) = j^{1}_{e_{1},X} \left(F\psi_{U}^{-1} \circ F\tau_{\psi_{U} \left(x\right)}\right) = j^{1}_{e_{1},X} \left(F\left(\psi_{U}^{-1} \circ\tau_{\psi_{U} \left(x\right)}\right)\right),
\end{equation}
where $\tau_{\psi_{U} \left(x\right)}$  denote the translation on $\mathbb{R}^{n}$ by the vector $\psi_{U} \left(x\right)$.\\\\
Next, we shall describe a particular subbundle of $\overline{F}^{2} M$. Consider the non-holonomic frames of second order given by $j^{1}_{e_{1},X} \left( F \psi \right)$, where $\psi : \mathbb{R}^{n} \rightarrow M$ is a local diffeomorphism at $0$. These kind of frames are called \textit{holonomic frames of second order} or \textit{second order frame bundle} (by short). The set of all holonomic frames of second order is denoted by $F^{2} M$ and it is called \textit{second-order holonomic frame bundle}. The restrictions of $\overline{\pi}^{2}_{1}$ and $\overline{\pi}^{2}$ to $F^{2}M$ are denoted by $\pi^{2}_{1}: F^{2}M \rightarrow FM$ and $\pi^{2}: F^{2}M \rightarrow M$. $\pi^{2}_{1}$ endows to $F^{2}M$ with a principal bundle structure with structure group $G^{2}_{1} \left( n \right)$, which is the set of all $1-$jets of local isomorphism of the form $F \psi$, where $\psi : \mathbb{R}^{n} \rightarrow \mathbb{R}^{n}$ is a local diffeomorphism with $F\psi \left( e_{1} \right) = e_{1}$ (equivalently, $j_{0,0}^{1} \phi = e_{1}$).\\
$\pi^{2} : F^{2}M \rightarrow M$ is also a principal bundle which structure group $G^{2} \triangleq \left( \pi^{2} \right)^{-1} \left( 0 \right)$.\\
We deduce that $\pi^{2}_{1}$ (resp. $\pi^{2}$) is a principal subbundle of $\overline{\pi}^{2}_{1}$ (resp. $\overline{\pi}^{2}$). So, restricting the isomorphism $\overline{F}^{2} \mathbb{R}^{n} \cong \mathbb{R}^{n} \times \overline{G}^{2} \left( n \right)$ we have that
$$ F^{2} \mathbb{R}^{n} \cong \mathbb{R}^{n} \times G^{2} \left( n \right).$$
Then, for each Lie subgroup $G$ of $G^{2} \left( n \right)$ we obtain a $G-$reduction of $F^{2} \mathbb{R}^{n}$ which is isomorphic to $ \mathbb{R}^{n} \times G$.
\begin{definition}\label{153}
\rm
Let $G$ be a Lie subgroup of $G^{2}\left(n\right)$. A \textit{second-order holonomic $G-$structure} $\omega_{G}\left(M\right)$ is a reduced subbundle of $F^{2}M$ with structure group $G$.
\end{definition}
Hence, the $G-$reduction of $F^{2}\mathbb{R}^{n}$ obtained above is a second-order holonomic $G-$structure on $\mathbb{R}^{n}$ which will be called the \textit{standard flat second-order holonomic $G-$structure}.\\

Note that each second-order holonomic $G-$structure $\omega_{G}\left(M\right)$ can be seen as a second-order non-holonomic $G-$structure. So, the notion of integrability will be the same.\\
A \textit{holonomic parallelism of second order} is a second-order holonomic trivial structure or, equivalently, a global section of $\pi^{2} : F^{2}M \rightarrow M$. So, we will also speak about \textit{integrable sections of $F^{2}M$}.\\

Observe that, by definition, any integrable non-holonomic parallelism of second order is in fact holonomic.\\
Summarizing we have the following sequence of Lie subgroups:
\begin{itemize}
\item[] $G^{2} \left( n \right) \subset \overline{G}^{2} \left( n \right) \subset Gl\left( n , \mathbb{R} \right) \times Gl\left( n +n^{2}, \mathbb{R} \right),$

\item[] $G^{2}_{1} \left( n \right) \subset \overline{G}^{2}_{1} \left( n \right) \subset  Gl\left( n +n^{2}, \mathbb{R} \right),$
\end{itemize}
and the following sequence of principal bundles:
$$ F^{2} M \subset \overline{F}^{2} M \subset F \left( F M \right),$$
over $FM$ and
$$ F^{2} M \subset \overline{F}^{2} M,$$
over $M$.\\
Let $\left( x \right)^{i}$ be a local coordinate system on an open $U \subseteq M$. Then, we can induce local coordinates over $\overline{F}^{2} M$ $\left(\left(x^{i}, x^{i}_{j}\right),x^{i}_{,j},x^{i}_{j,k}\right)$ (see Eq. (\ref{138})). Then, if we restrict to $F^{2}M$ we have that
$$ x^{i}_{j} = x^{i}_{,j} \ \ \ \ \ ; \ \ \ \ \ x^{i}_{j,k} = x^{i}_{k,j}.$$
Thus, we may obtain local coordinates on $F^{2}M$ denoted as follows:
\begin{equation}\label{139}
\left( x^{i} , x^{i}_{j} , x^{i}_{jk} \right), \ \ x^{i}_{jk} \triangleq x^{i}_{j,k} = x^{i}_{k,j}.
\end{equation}

\section{Non-holonomic prolongations of parallelisms of second order}

Let $M$ be a manifold and $\overline{P}$ be a section of the second-order non-holonomic frame bundle $\overline{F}^{2}M$. Then, $\overline{P}\left(x^{i} \right)= \left(\left(x^{i}, P^{i}_{j}\right) , P^{i}_{,j} , R^{i}_{j,k}\right)$ induces two sections $P$ and $Q$ of $FM$ (i.e. induces two ordinary parallelisms on $M$) by projecting $\overline{P}$ via the two canonical projections $\overline{\pi}^{2}_{1}$ and $\tilde{\pi}^{2}_{1}$, i.e.,
$$ P = \overline{\pi}^{2}_{1} \circ \overline{P}, \ \ \ \ \ Q = \tilde{\pi}^{2}_{1} \circ \overline{P}.$$
So, we obtain
$$ P \left(x^{i} \right) = \left(x^{i} , P^{i}_{j}\right) , \ \ \ \ \ Q\left(x^{i} \right) = \left(x^{i} , P^{i}_{,j} \right).$$
Conversely, let $P,Q: M \rightarrow FM$ be two sections of $FM$. Hence, $P$ (resp. $Q$) defines a family of $n$ (where $n$ is the dimension of $M$) linearly independent vector fields $\{P_{1} , \cdots , P_{n}\}$ (resp. $\{ Q_{1}, \cdots , Q_{n} \}$).\\
We define a horizontal subspace $H_{P\left(x\right)}$ at the point $P\left(x\right)$ by translating the basis $\{ Q_{a}\left(x\right)\}$ at $x$ into a set of linearly independent tangent vectors at $P\left(x\right)$,
$$ \{ T_{x} P \left( Q_{a} \left(x\right)\right) \}.$$
Locally,
$$ T_{x} P \left( \sum_{i} Q_{a}^{i} \left(x\right)\dfrac{\partial}{\partial x^{i}_{|x} }\right) = \sum_{i,r,s}Q_{a}^{i}\left(x\right)\dfrac{\partial}{\partial x^{i}_{|P\left(x\right)} } + Q^{i}_{a}\left(x\right)\dfrac{\partial P^{r}_{s}}{\partial x^{i}_{|x} }\dfrac{\partial}{\partial {x^{r}_{s}}_{|P\left(x\right)} },$$
where,
$$P_{a}\left(x\right) = \sum_{i}P_{a}^{i}\left(x\right)\dfrac{\partial}{\partial x^{i}_{|x} }, \ \ \ \ \ Q_{a}\left(x\right) = \sum_{i}Q_{a}^{i}\left(x\right)\dfrac{\partial}{\partial x^{i}_{|x} }.$$
By completing this set of linearly independent tangent vectors to a basis of $T_{P\left(x\right)} FM$ we obtain a second-order non-holonomic frame at $x$. We have so obtained a section of $\overline{F}^{2}M$ (i.e. a non-holonomic parallelism of second order on $M$), which is denoted by $P^{1}\left(Q\right)$.\\
\begin{definition}
\rm
A non-holonomic parallelism of second order $\overline{P}$ is said to be a \textit{prolongation} if $\overline{P} = P^{1}\left(Q\right)$ where $P$ and $Q$ are the induced ordinary parallelisms.

\end{definition}
The local expression of $P^{1}\left(Q\right)$ becomes
\begin{equation}\label{34}
P^{1}\left(Q\right) \left(x^{i}\right) = \left(\left(x^{i} , P^{i}_{j}\right), Q^{i}_{j}, \sum_{l} Q^{l}_{k}\dfrac{\partial P^{i}_{j}}{\partial x^{l} }\right).
\end{equation}
Hence, a section of $\overline{F}^{2}M$, $\overline{P}\left(x^{i} \right)= \left(\left(x^{i}, P^{i}_{j}\right) , Q^{i}_{j} , R^{i}_{j,k}\right)$, is a second-order non-holonomic prolongation if, and only if,
\begin{equation}\label{35}
R^{i}_{j,k} = \sum_{l} Q^{l}_{k}\dfrac{\partial P^{i}_{j}}{\partial x^{l} }.
\end{equation}

\begin{remarkth}\label{36}
\rm
Now, we will describe another way to construct $P^{1}\left(Q\right)$ which will be useful in what follows. Let $P,Q:M \rightarrow FM$ be two sections and we denote
$$ Q \left(x\right) = j_{0,x}^{1}\psi_{x}.$$
Then, for each $a =1, \cdots, n$
$$Q_{a}\left(x\right) = T_{0} \psi_{x} \left( \dfrac{\partial }{\partial x^{a}_{|0} }\right),$$
which implies that
$$T_{x}P\left(Q_{a}\left(x\right)\right) = T_{0}\left(P \circ \psi_{x}\right) \left( \dfrac{\partial }{\partial x^{a}_{|0} }\right).$$
Taking into account this equality we construct the following map

$$
\begin{array}{rccl}
\overline{P \circ \psi_{x}} : & FU  & \rightarrow & FV \\
& j_{0,v}^{1}f   &\mapsto &  P\left(\psi_{x}\left(v\right)\right) \cdot j_{0,0}^{1}\left( \tau_{-v} \circ f\right).
\end{array}
$$
where $\psi_{x} : U \rightarrow V$. It easy to show that $\overline{P \circ \psi_{x}}$ is an isomorphism of principal bundle over $\psi_{x}$ with inverse given by
$$ j_{0,w}^{1}g \in FV \mapsto j_{0, \psi_{x}^{-1} \left(w\right)}^{1} \tau_{\psi_{x}^{-1} \left(w\right)} \cdot [P \left(w\right)]^{-1} \cdot j_{0,w}^{1} g.$$
Thus, we can define,
$$
\begin{array}{rccl}
\overline{P } : & M  & \rightarrow & \overline{F}^{2}M \\
& x   &\mapsto &  j_{e_{1},Z}^{1}\left(\overline{P  \circ \psi_{x}}\right)
\end{array}
$$
which satisfies
\begin{itemize}
\item[(i)] $\overline{\pi}_{1}^{2} \circ \overline{P} \left(x\right) = \overline{P \circ \psi_{x}} \left( e_{1}\right) = P\left(x\right)$
\item[(ii)] $\tilde{\pi}_{1}^{2} \circ \overline{P} \left(x\right) = j_{0,z}^{1}\psi_{x} = Q\left(x\right)$
\item[(iii)] $\dfrac{\partial \left(P^{i}_{j} \circ  \Psi_{x}\right)}{\partial x^{k}_{|0} }= d{ P^{i}_{j}}_{|x} \circ \dfrac{\partial  \Psi_{x}}{\partial x^{k}_{|0} } = d{ P^{i}_{j}}_{|x} \circ \left( Q_{k}^{1}\left(x\right), \cdots ,Q_{k}^{n}\left(x\right)\right)$\\
Then, by definition of induced coordinates, $R^{i}_{j,k}$ is given by
$$R^{i}_{j,k}\left(x\right) =  \sum_{l} Q_{k}^{l}\left(x\right) \dfrac{\partial P^{i}_{j}}{\partial x^{l}_{|x} }.$$

\end{itemize}
Therefore $\overline{P} = P^{1}\left(Q\right)$.
\end{remarkth}
Notice that, if $Q$ is integrable, then there exists local coordinates $\left(x^{i}\right)$ on $M$ such that
$$ Q^{i}_{j} = \delta^{i}_{j},$$
and hence, 
$$R^{i}_{j,k} =\dfrac{\partial P^{i}_{j}}{\partial x^{k} },$$
where $P^{1}\left(Q\right)= \left(x^{i},P^{i}_{j},Q^{i}_{j}, R^{i}_{j,k}\right)$. In such a case, $P^{1}\left(Q\right)$ is said to be \textit{integrable}. However, in general $P^{1} \left( Q \right)$ is not integrable as a paralellism (see Proposition \ref{116}).\\
In this case for each $z \in M$, there exists a local chart $\left(\psi_{U} , U\right)$ over $z$ such that for each $x \in U$
$$P^{1}\left(Q\right) \left(x\right) = j_{e_{1},Z}^{1} \left( \overline{P  \circ \left(\psi_{U}^{-1} \circ \tau_{\psi_{U} \left( x \right)}\right)} \right) = 
j_{e_{1},Z}^{1} \left( \overline{P  \circ \psi_{U}^{-1}} \circ F\tau_{\psi_{U} \left( x \right) } \right).$$
In fact, for each $ j_{0,v}^{1} f \in F\psi_{U} \left(U\right)$ we have
$$\overline{P  \circ \left(\psi_{U}^{-1} \circ \tau_{\psi_{U} \left(x\right)}\right)} \left( j_{0,v}^{1} f\right) = P\left(\psi_{U}^{-1}\left(v+ \psi_{U} \left(x\right)\right)\right) \cdot j_{0,0}^{1} \left(\tau_{-v} \circ f\right)$$
$$ =\left(\overline{P  \circ \psi_{U}^{-1}} \circ F\tau_{\psi_{U} \left(x\right)}\right) \left( j_{0,v}^{1}f\right).$$
Thus, let $P^{1}\left(Q\right)$ be a non-holonomic integrable prolongation of second-order; then, for each $z \in M$, there exists a local principal bundle isomorphism $\Psi$ from an open set $F U \subseteq F M$ with $z \in U$ to an open subset of $F \mathbb{R}^{n}$ such that for all $x \in U$
\begin{equation}\label{37}
P^{1}\left(Q\right) \left(x\right) = j^{1}_{e_{1} , Z} \left(\Psi^{-1} \circ F\tau_{\psi \left(x\right)}\right),
\end{equation}
where $\psi$ is the induced map of $\Psi$ onto the base manifolds. On the other hand, using $\psi$ as a local chart we can prove that Eq. (\ref{37}) implies that $P^{1}\left(Q\right)$ is a non-holonomic integrable prolongation of second order.\\
This equality reminds us Eq. (\ref{33}), for second-order non-holonomic integrable sections. Indeed, a second-order non-holonomic integrable prolongation $P^{1}\left(Q\right)$ satisfying Eq. (\ref{37}) is integrable if, and only if, we have
$$j_{Z,{e_{1}}_{x}}^{1} \Psi = j_{Z,{e_{1}}_{x}}^{1} F \psi,$$
for all $x$, where ${e_{1}}_{x} =j_{0,x}^{1} \tau_{\psi \left(x\right)}$. Thus, a second-order non-holonomic prolongation is integrable if, and only if, takes values in the holonomic frame, i.e., the only integrable prolongations in $F^{2}M$ are the integrable sections.\\

In general, we can think about the second-order non-holonomic integrable prolongations as an intermediate step beetwen sections and integrable sections. The following result is obvious.
\begin{proposition}\label{116}
Let $\overline{P}$ be a section of $\overline{F}^{2}M$. $\overline{P}$ is integrable if, and only if, $\overline{P} = P^{1}\left(Q\right)$ is a second-order non-holonomic integrable prolongation and $P = Q$. In particular, a second-order non-holonomic integrable prolongation $P^{1}\left(Q\right)$ is integrable if, and only if, $P = Q$.
\end{proposition}

This result provides us with examples of second-order no-holonomic integrable prolongations which are not integrable. Indeed, any prolongation of two different ordinary parallelisms is not a second-order no-holonomic integrable prolongation.\\

Now, to end this section, we will define the concept prolongation for non-holonomic $\overline{G}-$structures of second order.

\begin{definition}\label{38}
\rm
Let $\overline{\omega}_{\overline{G}} \left(M\right)$ be a second-order non-holonomic \linebreak$\overline{G}-$structure. $\overline{\omega}_{\overline{G}} \left(M\right)$ is a non-holonomic integrable prolongation of second-order if we can cover $M$ by local non-holonomic integrable prolongations of second order which take values in $\overline{\omega}_{\overline{G}} \left(M\right)$.
\end{definition}

Notice that Definition \ref{38} can be expressed as follows: For any point $x \in M$ there exists a local coordinate system $\left(x^{i}\right)$ over an open set $U \subseteq M$ which contains $x$ such that the local section on $U$
$$P^{1}\left(Q\right) \left(x^{i}\right) = \left(x^{i},P^{i}_{j},\delta^{i}_{j}, \dfrac{\partial P^{i}_{j}}{\partial x^{k} }\right),$$
is contained in $\overline{\omega}_{\overline{G}} \left(M\right)$.\\

As we have noticed, a second-order non-holonomic $\overline{G}-$structure which is contained in $F^{2}M$ is integrable if, and only if, it is an integrable prolongation.\\

\begin{remarkth}\label{39}
\rm

Let $\overline{\omega}_{\overline{G}} \left(M\right)$ be a second-order non-holonomic \linebreak$\overline{G}-$structure. Then, we could define a non-holonomic integrable prolongation of second order in a similar way to integrable $\overline{G}-$structures.\\
In fact, using Eq. (\ref{37}), we can prove that $\overline{\omega}_{\overline{G}}\left(M\right)$ is a non-holonomic integrable prolongation of second-order if, and only if, for all point $x \in M$, there exists a local isomorphism of principal bundles whose isomorphism of Lie groups is the identity map, $\Phi: FU \rightarrow F\overline{U}$, with $x \in U$ such that it induces an isomorphism of principal bundles given by
$$ \Upsilon : \overline{\omega}_{\overline{G}} \left(U\right)  \rightarrow \overline{U} \times \overline{G},$$
where $\Upsilon\left( j_{e_{1} , Z}^{1} \Psi \right) = \left(\phi \left(z\right) , \overline{\Upsilon} \left( j_{e_{1} , Z}^{1} \Psi \right)\right)$ and
$$ \overline{\Upsilon} \left( j_{e_{1} , Z}^{1} \Psi \right) =  j_{e_{1} ,Z_{0}}^{1} \left( F \left( \tau_{-\phi \left( z\right)}\right) \circ \Phi  \circ \Psi \right) ,$$
with $\phi$ the induced map of $\Phi$ over the base manifold and $ \pi_{M}\left(Z\right)=z$.\\
Thus, $\overline{\omega}_{\overline{G}} \left(M\right)$ is a non-holonomic integrable prolongation of second-order if it is locally isomorphic to the trivial principal bundle $\mathbb{R}^{n} \times \overline{G}$ by a more general class of local charts (see Definition \ref{41}).\\\\

\end{remarkth}

\begin{proposition}
Let $\overline{\omega}_{\overline{G}} \left(M\right)$ be a second-order non-holonomic \linebreak$\overline{G}-$structure. If $\overline{\omega}_{\overline{G}} \left(M\right)$ is integrable, then $\overline{\omega}_{\overline{G}} \left(M\right)$ is a non-holonomic integrable prolongation of second-order.
\end{proposition}

Not all non-holonomic integrable prolongation of second-order is integrable (see Proposition \ref{116}).\\

It directly follows that if $\overline{\omega}_{\overline{G}} \left(M\right)$ is a second-order non-holonomic integrable prolongation, then the projected $G-$structure by $\tilde{\pi}^{2}_{1}$ is integrable.\\

\section{Groupoids}
In this section, we will briefly introduce the notion of \textit{Lie groupoid} which is a fundamental concept for the rest of this paper (a good reference on groupoids is \cite{KMG}). A particular groupoid, the \textit{1-jets groupoid}, will be closely related with the frame bundle of the base manifold.

\begin{definition}
\rm
Let $ M$ be a set. A \textit{groupoid} over $M$ is given by a set $\Gamma$ provided with the maps $\alpha,\beta : \Gamma \rightarrow M$ (\textit{source map} and \textit{target map} respectively), $\epsilon: M \rightarrow \Gamma$ (\textit{identities map}), $i: \Gamma \rightarrow \Gamma$ (\textit{inversion map}) and   
$\cdot : \Gamma_{\left(2\right)} \rightarrow \Gamma$ (\textit{composition law}) where for each $k \in \mathbb{N}$,
$$\Gamma_{\left(k\right)} := \{ \left(g_{1}, \hdots , g_{k}\right) \in \Gamma \times \stackrel{k)}{\ldots} \times \Gamma \ : \ \alpha\left(g_{i}\right)=\beta\left(g_{i+1}\right), \ i=1, \hdots , k -1 \},$$
which satisfy the following properties:\\
\begin{itemize}
\item[(1)] $\alpha$ and $\beta$ are surjective and for each $\left(g,h\right) \in \Gamma_{\left(2\right)}$,
$$ \alpha\left(g \cdot h \right)= \alpha\left(h\right), \ \ \ \beta\left(g \cdot h \right) = \beta\left(g\right).$$
\item[(2)] Associative law with the composition law, i.e.,
$$ g \cdot \left(h \cdot k\right) = \left(g \cdot h \right) \cdot k, \ \forall \left(g,h,k\right) \in \Gamma_{\left(3\right)}.$$
\item[(3)] For all $ g \in \Gamma$,
$$ g \cdot \epsilon \left( \alpha\left(g\right)\right) = g = \epsilon \left(\beta \left(g\right)\right)\cdot g .$$
In particular,
$$ \alpha \circ  \epsilon \circ \alpha = \alpha , \ \ \ \beta \circ \epsilon \circ \beta = \beta.$$
Since $\alpha$ and $\beta$ are surjetive we get
$$ \alpha \circ \epsilon = Id_{\Gamma}, \ \ \ \beta \circ \epsilon = Id_{\Gamma}.$$
\item[(4)] For each $g \in \Gamma$,
$$i\left(g\right) \cdot g = \epsilon \left(\alpha\left(g\right)\right) , \ \ \ g \cdot i\left(g\right) = \epsilon \left(\beta\left(g\right)\right).$$
Then,
$$ \alpha \circ i = \beta , \ \ \ \beta \circ i = \alpha.$$
\end{itemize}
These maps will be called \textit{structure maps}. We will denote this groupoid by $ \Gamma \rightrightarrows M$.
\end{definition}
If $\Gamma$ is a groupoid over $M$, then $M$ is also denoted by $\Gamma_{\left(0\right)}$ and it is often identified with the set $\epsilon \left(M\right)$ of identity elements of $\Gamma$. $\Gamma$ is also denoted by $\Gamma_{\left(1\right)}$. The space of sections of the map $\left(\alpha , \beta\right) : \Gamma \rightarrow M \times M$ is denoted by $\Gamma_{\left(\alpha, \beta\right)} \left(\Gamma\right)$.\\\\
Now, we define the morphisms of the category of groupoids.\\
\begin{definition}
\rm

If $\Gamma_{1} \rightrightarrows M_{1}$ and $\Gamma_{2} \rightrightarrows M_{2}$ are two groupoids; then a morphism from $\Gamma_{1} \rightrightarrows M_{1}$ to $\Gamma_{2} \rightrightarrows M_{2}$ consists of two maps $\Phi : \Gamma_{1} \rightarrow \Gamma_{2}$ and $\phi : M_{1} \rightarrow M_{2}$ such that for any $g_{1} \in \Gamma_{1}$
\begin{equation}\label{4}
\alpha_{2} \left( \Phi \left(g_{1}\right)\right) = \phi \left(\alpha_{1} \left(g_{1} \right)\right), \ \ \ \ \ \ \ \beta_{2} \left( \Phi \left(g_{1}\right)\right) = \phi \left(\beta_{1} \left(g_{1} \right)\right),
\end{equation}
where $\alpha_{i}$ and $\beta_{i}$ are the source and the target map of $\Gamma_{i} \rightrightarrows M_{i}$ respectively, for $i=1,2$, and preserves the composition, i.e.,
$$\Phi \left( g_{1} \cdot h_{1} \right) = \Phi \left(g_{1}\right) \cdot \Phi \left(h_{1}\right), \ \forall \left(g_{1} , h_{1} \right) \in \Gamma_{\left(2\right)}.$$
We will denote this morphism as $\left(\Phi , \phi\right)$.
\end{definition}
Observe that, as a consequence, $\Phi$ preserves the identities, i.e., denoting by $\epsilon_{i}$ the section of identities of $\Gamma_{i} \rightrightarrows M_{i}$ for $i=1,2$,
$$\Phi \circ  \epsilon_{1} = \epsilon_{2} \circ \phi .$$
Note that, using equations \ref{4}, $\phi$ is completely determined by $\Phi$.\\\\
Using this definition we define a \textit{subgroupoid} of a groupoid $\Gamma \rightrightarrows M$ as a groupoid $\Gamma' \rightrightarrows M'$ such that $M' \subseteq M$, $\Gamma' \subseteq \Gamma$ and the inclusion map is a morphism of groupoids.

\begin{remarkth}
\rm
There is a more abstract way of defining a groupoid. We can say that a groupoid is a "small" category (the class of objects and the class of morphisms are sets) in which each morphism is invertible.\\
If $ \Gamma \rightrightarrows M$ is the groupoid, then $M$ is the set of objects and $\Gamma$ is the set of morphisms.\\
A groupoid morphism is a functor between these categories which is a more natural definition.
\end{remarkth}
Now, we present some basic examples of groupoids.

\begin{example}\label{5}
\rm
A group is a groupoid over a point. In fact, let $G$ be a group and $e$ the identity element of $G$. Then, $G \rightrightarrows \{e\}$ is a groupoid, where the operation of the groupoid, $\cdot$, is the operation in $G$.
\end{example}
\begin{example}\label{7}

\rm
For any set $M$, the product space $ M \times M$ can be seen as a groupoid over $M$ where the composition is given by:
$$\left(y,z\right)\cdot\left(x,y\right)= \left(x,z\right), \ \forall \left(x,y\right),\left(y,z\right) \in  M \times M.$$ 
This groupoid is called the \textit{pair groupoid on $M$}. Note that, if $\Gamma \rightrightarrows M$ is an arbitrary groupoid over $M$, then the map $\left(\alpha , \beta\right) : \Gamma \rightarrow M \times M$, which is sometimes called the \textit{anchor} of $\Gamma$, is a morphism from $\Gamma \rightrightarrows M$ to the pair groupoid of $M$.
\end{example}

\begin{example}\label{8}
\rm
Let $A$ be a vector bundle over a manifold $M$. Let $\Phi \left(A\right)$ denote the set of all vector space isomorphisms $L_{x,y}: A_{x} \rightarrow A_{y}$ for $x,y \in M$, where for each $z\in M$, $A_{z}$ is the fibre of $A$ over $z$. We can consider $\Phi \left(A\right)$ as a groupoid $\Phi \left(A\right) \rightrightarrows M$ such that, for all $x,y \in M$ and $L_{x,y}\in \Phi \left(A\right)$,
\begin{itemize}
\item[(i)] $\alpha\left(L_{x,y}\right) = x$
\item[(ii)] $\beta\left(L_{x,y}\right) = y$
\item[(iii)] $L_{y,z} \cdot G_{x,y} = L_{y,z} \circ G_{x,y}, \ L_{y,z}: A_{y} \rightarrow A_{z}, \ G_{x,y}: A_{x} \rightarrow A_{y}$
\end{itemize}
This groupoid is called the \textit{frame groupoid on $A$}. As a particular case, when $A$ is the tangent bundle over $M$ we have the \textit{1-jets groupoid on} $M$ which is denoted by $\Pi^{1} \left(M,M\right)$.
\end{example}

Next, we will introduce the notions of orbit and isotropy group of an action.
\begin{definition}
\rm
Let $\Gamma \rightrightarrows M$ be a groupoid with $\alpha$ and $\beta$ the source map and target map, respectively. For each $x \in M$, then
$$\Gamma^{x}_{x}= \beta^{-1}\left(x\right) \cap \alpha^{-1}\left(x\right),$$
is called the i\textit{sotropy group of} $\Gamma$ at $x$. The set
$$\mathcal{O}\left(x\right) = \beta\left(\alpha^{-1}\left(x\right)\right) = \alpha\left(\beta^{-1}\left( x\right)\right),$$
is called the \textit{orbit} of $x$, or \textit{the orbit} of $\Gamma$ through $x$.\\\\
If $\mathcal{O}\left(x\right) =\{x\}$, or equivalently, $\beta^{-1}\left(x\right)=\alpha^{-1} \left(x\right)= \Gamma^{x}_{x}$ then $x$ is called a \textit{fixed point}. \textit{The orbit space of} $\Gamma$ is the space of orbits of $\Gamma $ on $M$, i.e., the quotient space of $M$ by the equivalence relation induced by $\Gamma$: two points of $M$ are equivalent if, and only if, they lie on the same orbit.\\
If $\mathcal{O}\left(x\right) = M$ for all $x \in M$, or equivalently $\left(\alpha,\beta\right) : \Gamma  \rightarrow M \times M$ is a surjective map, the groupoid $\Gamma \rightrightarrows M$ is called \textit{transitive}. If every $x \in M$ is fixed point, then the groupoid $\Gamma \rightrightarrows M$ is called \textit{totally intransitive}. Furthermore, a subset $N$ of $M$ is called \textit{invariant} if it is a union of some orbits.\\
Finally, the preimages of the source map $\alpha$ of a Lie groupoid are called $\alpha-$\textit{fibres}. Those of the target map $\beta$ are called $\beta-$\textit{fibres}.\\
\end{definition}

\begin{definition}\label{9}
\rm
Let $\Gamma \rightrightarrows M$ be a groupoid with $\alpha$ and $\beta$ the source and target map, respectively. We may define the right translation on $g \in \Gamma$ as the map $R_{g} : \alpha^{-1}\left(\beta\left(g\right)\right) \rightarrow \alpha^{-1} \left( \alpha \left(g\right)\right)$, given by
$$ h \mapsto  h \cdot g.$$
We may define the left translation on $g$, $L_{g} : \beta^{-1} \left( \alpha\left(g\right)\right) \rightarrow \beta^{-1} \left(\beta\left(g\right)\right)$ similarly. 
\end{definition}
Note that,
\begin{equation}\label{10} 
Id_{\alpha^{-1}\left(x\right)} = R_{\epsilon \left(x\right)}.
\end{equation}
So, for all $ g \in \Gamma $, the right (left) translation on $g$, $R_{g}$ (resp. $L_{g}$), is a bijective map with inverse $R_{i\left(g\right)}$ (resp. $L_{i\left(g\right)}$), where $i : \Gamma \rightarrow \Gamma$ is the inverse map.\\
One may impose various topological and geometrical structures on a groupoid, depending on the context. We will be mainly interested in Lie groupoids.
\begin{definition}
\rm
A \textit{Lie groupoid} is a groupoid $\Gamma \rightrightarrows M$ such that $\Gamma$ is a smooth manifold, $M$ is a smooth manifold and all the structure maps are smooth. Furthermore, the source and the target map are submersions.\\
A \textit{Lie groupoid morphism} is a groupoid morphism which is differentiable.\\
\end{definition}

\begin{definition}
\rm
Let $\Gamma \rightrightarrows M$ be a Lie groupoid. A \textit{Lie subgroupoid} of $\Gamma \rightrightarrows M$ is a Lie groupoid $\Gamma' \rightrightarrows M'$ such that $\Gamma' $ and $M'$ are submanifolds of $\Gamma$ and $M$ respectively, and the inclusion maps $i_{\Gamma'} : \Gamma' \hookrightarrow \Gamma$ $i_{M'} : M' \hookrightarrow M$ become a morphism of Lie groupoids. $\Gamma' \rightrightarrows M'$ is said to be a \textit{reduced Lie subgroupoid} if it is transitive and $M=M'$.
\end{definition}
Observe that, taking into account that $ \alpha \circ \epsilon = Id_{M} = \beta \circ \epsilon$, then $\epsilon$ is an injective immersion.\\\\
On the other hand, in the case of a Lie groupoid, $R_{g}$ (resp. $L_{g}$) is clearly a diffeomorphism for all $g \in \Gamma$.\\\\
Note also that, for each $k \in \mathbb{N}$, $\Gamma_{\left(k\right)}$ is a pullback space given by $\beta$ and the operation map on $\Gamma_{\left(k-1\right)}$. Thus, by induction, we may prove that $\Gamma_{\left(k\right)}$ is a smooth manifold for all $k \in \mathbb{N}$.\\

\begin{example}\label{11}
\rm
A Lie group is a Lie groupoid over a point. 
\end{example}
\begin{example}
\rm
Let $M$ be a manifold. It is trivial to prove that the pair groupoid on $M$ is a Lie groupoid.
\end{example}
\begin{example}\label{13}
\rm
Let $\pi: P \rightarrow M$ be a principal bundle with structure group $G$. Denote by $\phi :  P  \times G \rightarrow P$ the action of $G$ on $P$.\\
Now, suppose that $\Gamma \rightrightarrows P$ is a Lie groupoid, with $\overline{\phi} :  \Gamma \times G \rightarrow \Gamma$ a free and proper action of $G$ on $\Gamma$ such that, for each $h \in G$, the pair $\left( \overline{\phi}_{h} , \phi_{h}\right)$ is an isomorphism of Lie groupoids.\\
We can construct a Lie groupoid $\Gamma / G \rightrightarrows M$ such that the source map, $\overline{\alpha}$, and the target map, $\overline{\beta}$, are given by
$$ \overline{\beta}\left([ g ]\right) = \pi \left( \beta \left(g\right)\right), \ \ \overline{\alpha}\left([g ] \right) = \pi \left( \alpha\left(g\right)\right),$$
for all $ g \in \Gamma$, $\alpha$ and $\beta$ being the source and the target map on $\Gamma \rightrightarrows P$, respectively, and $[\cdot]$ denotes the equivalence class in the quotient space $\Gamma / G$. These kind of Lie groupoids are called \textit{quotient Lie groupoids by the action of a Lie group}.
\end{example}
There is an interesting particular case of the above example.
\begin{example}\label{14}
\rm
Let $\pi: P \rightarrow M$ be a principal bundle with structure group $G$ and $P \times P \rightrightarrows P$ the pair groupoid. Take $\overline{\phi} :  \left(P \times P\right) \times G \rightarrow P \times P$ the diagonal action of $\phi$, where $\phi :  P \times G \rightarrow P$ is the action of $G$ on $P$.\\
Then it is easy to prove that $\left(\overline{\phi}_{g}, \phi_{g}\right)$ is an isomorphism of Lie groupoids and thus, we may construct the groupoid $\left( P \times P \right) / G \rightrightarrows M$. This groupoid is called \textit{gauge groupoid} and is denoted by $Gauge\left(P\right)$.
\end{example}

\begin{example}\label{99}
\rm
Let $G$ be a Lie group and $M$ be a manifold. The \textit{Trivial Lie groupoid on $M$ with group $G$} is given by $M \times M \times G \rightrightarrows M$ where the structure maps are the following:
\begin{itemize}
\item[(i)] $\alpha\left(x,y,g \right) = x$
\item[(ii)] $\beta\left(x,y,g \right) = y$
\item[(iii)] $\left( y,z , h \right) \cdot \left( x,y , g \right) = \left( x,y,h \cdot g \right)$
\end{itemize}

\end{example}

\begin{example}\label{15}
\rm
Let $A$ be a vector bundle over $M$ then the frame groupoid is a Lie groupoid (see Example \ref{8}). In fact, let $\left(x^{i}\right)$ and $\left(y^{j}\right)$ be local coordinate systems on open sets $U, V \subseteq M$ and $\{\alpha_{p}\}$ and $\{ \beta_{q}\}$ be local basis of sections of $A_{U}$ and $A_{V}$ respectively. The corresponding local coordinates $\left(x^{i} \circ \pi, \alpha^{p}\right)$ and $\left(y^{j} \circ \pi, \beta^{q}\right)$ on $A_{U}$ and $A_{V}$ are given by
\begin{itemize}
\item $ a = \sum_{p} \alpha^{p}\left(a\right) \alpha_{p}\left(x^{i}\left(\pi \left(a\right)\right)\right),\ \forall a \in A_{U}$
\item $ b = \sum_{q} \beta^{q}\left(b\right) \beta_{q}\left(y^{j}\left(\pi \left(b\right)\right)\right), \ \forall b \in A_{V}$
\end{itemize}

Then, we can consider a local coordinate system $\Phi \left(A\right)$
$$ \Phi \left(A_{U,V}\right) : \left(x^{i} , y^{j}, y^{j}_{i}\right),$$
where, $A_{U,V} = \alpha^{-1}\left(U\right) \cap \beta^{-1}\left(V\right)$ and for each $L_{x,y} \in \alpha^{-1}\left(x\right) \cap \beta^{-1}\left(y\right) \subseteq \alpha^{-1}\left(U\right) \cap \beta^{-1}\left(V\right)$, we have
\begin{itemize}
\item $x^{i} \left(L_{x,y}\right) = x^{i} \left(x\right)$
\item $y^{j} \left(L_{x,y}\right) = y^{j} \left( y\right)$
\item $y^{j}_{i}\left( L_{x,y}\right)  = A_{L_{x,y}}$, where $A_{L_{x,y}}$ is the induced matrix of the induced map of $L_{x,y}$ by the local coordinates $\left(x^{i} \circ \pi, \alpha^{p}\right)$ and $\left(y^{j} \circ \pi, \beta^{q}\right)$
\end{itemize}
In particular, if $A=TM$, then the $1-$jets groupoid on $M$, $\Pi^{1}\left(M,M\right)$, is a Lie groupoid and its local coordinates will be denoted as follows
\begin{equation}\label{17}
\Pi^{1}\left(U,V\right) : \left(x^{i} , y^{j}, y^{j}_{i}\right),
\end{equation}
where, for each $ j^{1}_{x,y} \psi \in \Pi^{1}\left(U,V\right)$
\begin{itemize}
\item $x^{i} \left(j^{1}_{x,y} \psi\right) = x^{i} \left(x\right)$
\item $y^{j} \left(j^{1}_{x,y}\psi \right) = y^{j} \left( y\right)$
\item $y^{j}_{i}\left( j^{1}_{x,y}\psi\right)  = \dfrac{\partial \left(y^{j}\circ \psi\right)}{\partial x^{i}_{|x} }$
\end{itemize}
\end{example}

Next, as an important example, we will introduce the \textit{second-order non-holonomic groupoid}.\\\\
Let $M$ be a manifold and $FM$ the frame bundle over $M$. So, we can consider the $1-$jets groupoid on $FM$, $\Pi^{1} \left(FM, FM\right) \rightrightarrows FM$.\\
Thus, we denote by $J^{1}\left(FM \right)$ the subset of $\Pi^{1} \left(FM,FM\right)$ given by the $1-$jets $j^{1}_{X,Y} \Psi$ of local automorphism $\Psi$ of $FM$ such that
$$ \Psi \left( v \cdot g\right) = \Psi \left( v \right) \cdot g, \ \forall v \in Dom\left( \Psi \right), \ \forall g \in Gl \left( n , \mathbb{R} \right).$$
Let $\left(x^{i}\right)$ and $\left(y^{j}\right)$ be local coordinate systems over two open sets $U,V \subseteq M$, the induced coordinate systems over $FM$ are denoted by
$$FU: \left(x^{i} , x^{i}_{j}\right)$$
$$FV: \left(y^{j} , y^{j}_{i}\right).$$
Hence, we can construct induced coordinates over $\Pi^{1}\left(FM,FM\right)$
$$\Pi^{1}\left(FU,FV\right) \triangleq \left( \alpha , \beta \right)^{-1}\left(U , V\right) : \left(\left(x^{i},x^{i}_{j}\right),\left( y^{j} , y^{j}_{i}\right),y^{j}_{,i}, y^{j}_{,ik},y^{j}_{i,k} ,y^{j}_{i,kl}\right) ,$$
where for each $j^{1}_{X,Y} \Psi \in \Pi^{1}\left(FU,FV\right)$, we have
\begin{itemize}
\item $x^{i} \left(j^{1}_{X,Y} \Psi\right) = x^{i} \left(X\right)$
\item $x^{i}_{j} \left(j^{1}_{X,Y} \Psi\right) =x^{i}_{j} \left(X\right)$
\item $y^{j} \left(j^{1}_{X,Y} \Psi\right) = y^{j} \left( \Psi \left(X\right)\right)$
\item $y^{j}_{i}\left( j^{1}_{X,Y} \Psi\right)  = y^{j}_{i} \left( \Psi \left( X\right)\right)$
\item $y^{j}_{,i} \left( j^{1}_{X,Y} \Psi \right) =  \dfrac{\partial \left(y^{j}\circ \Psi\right)}{\partial x^{i}_{|X} }$
\item $y^{j}_{,ik} \left( j^{1}_{X,Y} \Psi \right) =  \dfrac{\partial \left(y^{j}\circ \Psi\right)}{\partial {x^{i}_{k}}_{|X} }$
\item $y^{j}_{i,k} \left( j^{1}_{X,Y} \Psi \right) =  \dfrac{\partial \left(y^{j}_{i}\circ \Psi\right)}{\partial {x^{k}}_{|X} }$
\item $y^{j}_{i,kl} \left( j^{1}_{X,Y} \Psi \right) =  \dfrac{\partial \left(y^{j}_{i}\circ \Psi\right)}{\partial {x^{k}_{l}}_{|X} }$
\end{itemize}
Then, using these coordinates, $J^{1}\left(FM\right)$ can be described as follows:
{\footnotesize $$J^{1}\left(FU,FV\right) \triangleq J^{1}\left(FM\right) \cap \left(\alpha , \beta \right)^{-1}\left(U,V\right) : \left(\left(x^{i},x^{i}_{j}\right),\left( y^{j} , y^{j}_{i}\right),y^{j}_{,i}, 0,y^{j}_{i,k},y^{j}_{i,kl}\right) ,$$}
where
$$ y^{j}_{i,kl} = \left( \sum_{m}y^{j}_{m} \left( x^{-1} \right)^{m}_{k} \right) \delta_{l}^{i}.$$
Thus, $J^{1}\left(FM\right)$ is a submanifold of $\Pi^{1}\left(FM,FM\right)$ and its induced local coordinates will be denoted by
\begin{equation}\label{101}
J^{1}\left(FU,FV\right) : \left(\left(x^{i},x^{i}_{j}\right),\left( y^{j} , y^{j}_{i}\right),y^{j}_{,i},y^{j}_{i,k}\right) .
\end{equation}
Finally, restricting the structure maps we can ensure that $J^{1}\left(FM\right) \rightrightarrows FM$ is a reduced Lie subgroupoid of the $1-$jets groupoid over $FM$.\\\\
Analogously to $F^{2}M$, we may construct $j^{1}\left(FM\right)$ as the set of the $1-$jets of the form $j^{1}_{X,Y} F\psi$, where $\psi: M \rightarrow M$ is a local diffeomorphism. Let $\left( x^{i} \right)$ be a local coordinate system on $M$; then, restricting the induced local coordinates given in Eq. (\ref{101}) to $j^{1} \left( FM \right)$ we have that
$$ y^{j}_{i} = y^{j}_{,l}x^{l}_{i} \ \ \ \ \ ; \ \ \ \ \ y^{j}_{i,k} = y^{j}_{k,i}.$$
We deduce that $j^{1}\left(FM\right) \rightrightarrows FM$ is a reduced Lie subgroupoid of the $1-$jets groupoid over $FM$ and we denoted the coordinates on $j^{1} \left(FM \right)$ by
\begin{equation}\label{140}
j^{1}\left(FU,FV\right) : \left(\left(x^{i},x^{i}_{j}\right),\left( y^{j} , y^{j}_{i}\right),y^{j}_{i,k}\right) , \ \ \ y^{j}_{i,k} = y^{j}_{k,i}.
\end{equation}
Now, we will work with a quotient space of $J^{1}\left(FM\right)$ (resp. $j^{1}\left(FM\right)$) which will be our \textit{non-holonomic groupoid of second order} (resp. \textit{holonomic groupoid of second order}).\\\\
We consider the following right action of $Gl\left(n , \mathbb{R}\right)$ over $J^{1} \left(FM\right)$,

\begin{equation}\label{44}
\begin{array}{rccl}
\Phi : & J^{1} \left(FM\right) \times Gl \left( n , \mathbb{R}\right)  & \rightarrow & J^{1}\left(FM\right) \\
&\left(j^{1}_{X,Y} \Psi , g\right)   &\mapsto &  j^{1}_{X \cdot g , Y \cdot g} \Psi.
\end{array}
\end{equation}

Thus, for each $g \in Gl \left( n , \mathbb{R}\right)$ the pair $\left(\Phi_{g} , R_{g}\right)$ (where $R$ is the natural right action of $Gl\left(n , \mathbb{R}\right)$ over $FM$) is a Lie groupoid automorphism. Therefore, we can consider the quotient Lie groupoid by this action $ \tilde{J}^{1} \left( FM \right) \rightrightarrows M$ which is called \textit{second-order non-holonomic groupoid over $M$}.\\
We will denote the structure maps of $\tilde{J}^{1} \left(FM\right)$ by $\overline{\alpha}$ and $\overline{\beta}$ (source and target maps respectively), $\overline{\epsilon}$ (identities map) and $\overline{i}$ (inversion map). The elements of $\tilde{J}^{1} \left( FM \right)$ are denoted by $j^{1}_{x,y} \Psi$ with $x,y \in M$ and $\overline{\alpha}\left( j^{1}_{x,y} \Psi\right) = x$ and $\overline{\beta}\left(j^{1}_{x,y} \Psi\right) = y$.\\
Then, the induced local coordinates are given by
\begin{equation}\label{45}
\tilde{J}^{1}\left(FU,FV\right) \triangleq \left( \overline{\alpha} , \overline{\beta} \right)^{-1}\left(U , V\right) : \left(\left(x^{i}\right),\left( y^{j} , y^{j}_{i}\right),y^{j}_{,i},y^{j}_{i,k}\right) .
\end{equation}
Considering ${e_{1}}_{x}$ as the $1-$jet through $x \in M$ which satisfies that $x^{i}_{j} \left( {e_{1}}_{x}\right) = \delta^{i}_{j}$ for all $i,j$, for each $j_{x,y}^{1}\Psi \in \tilde{J}^{1}\left(FM\right)$ we have
\begin{itemize}
\item $x^{i} \left(j^{1}_{x,y} \Psi\right) = x^{i} \left(x\right)$
\item $y^{j} \left(j^{1}_{x,y} \Psi\right) = y^{j} \left( y\right)$
\item $y^{j}_{i}\left( j^{1}_{x,y} \Psi\right)  = y^{j}_{i} \left( \Psi \left( {e_{1}}_{x}\right)\right)$.
\item $y^{j}_{,i} \left( j^{1}_{x,y} \Psi \right) =  \dfrac{\partial \left(y^{j}\circ \Psi\right)}{\partial x^{i}_{|{e_{1}}_{x}} }$
\item $y^{j}_{i,k} \left( j^{1}_{x,y} \Psi \right) =  \dfrac{\partial \left(y^{j}_{i}\circ \Psi\right)}{\partial {x^{k}}_{|{e_{1}}_{x}} }$
\end{itemize}

Observe that we can restrict the action $\Phi$ to an action of $Gl \left( n , \mathbb{R} \right)$ over $j^{1} \left(FM\right)$. So, by quotienting, we can build a reduced subgroupoid of $ \tilde{J}^{1} \left( FM \right) \rightrightarrows M$ which is denoted by $ \tilde{j}^{1} \left( FM \right) \rightrightarrows M$ and is called \textit{second-order holonomic groupoid over $M$}. Finally, by restriction, the local coordinates on $j^{1} \left( FM \right)$ are given by
\begin{equation}
\tilde{j}^{1}\left(FU,FV\right) : \left(\left(x^{i}\right),\left( y^{j} , y^{j}_{i}\right),y^{j}_{i,k}\right) , \ \ \ y^{j}_{i,k} = y^{j}_{k,i}.
\end{equation}

\section{Lie Algebroids}
The notion of \textit{Lie algebroid} was introduced by J. Pradines in $1966$ \cite{JPRA} as an infinitesimal version of Lie groupoid and for this reason the author called it \textit{infinitesimal groupoid}. Now, we will recall this concept (we also refer to \cite{KMG}).

\begin{definition}
\rm
A \textit{Lie algebroid over $M$} is a triple $\left(  A \rightarrow M, \sharp , [ \cdot , \cdot ] \right)$, where $\pi : A \rightarrow M$ is a vector bundle together with a vector bundle morphism $\sharp : A \rightarrow TM$, called the \textit{anchor}, and a Lie bracket $[ \cdot  , \cdot ]$ on the space of sections, such that the Leibniz rule holds
\begin{equation}\label{18}
[ \alpha , f \beta ] = f [\alpha , \beta ] + {\sharp}\left(\alpha\right)\left(f\right)\beta,
\end{equation}
for all $\alpha , \beta \in \Gamma \left(A\right)$ and $f \in \mathcal{C}^{\infty} \left(M\right)$.\\
$A$ is \textit{transitive} if $\sharp$ is surjective and \textit{totally intransitive} if $\sharp \equiv 0$. Also, $A$ is said to be \textit{regular} if $\sharp$ has constant rank.
\end{definition}

Looking at $\sharp$ as a $\mathcal{C}^{\infty}\left(M\right)$-module morphism from $\Gamma \left(A\right)$ to $\mathfrak{X} \left(M\right)$, for each section $\alpha \in \Gamma \left(A\right)$ we are going to denote $ \sharp \left(\alpha\right)$ by $\alpha^{\sharp}$. Next, let us show the following fundamental property:
\begin{lemma}\label{19}
If $\left(  A \rightarrow M, \sharp , [ \cdot , \cdot ] \right)$ is a Lie algebroid, then the anchor map is a morphism of Lie algebras, i.e.
\begin{equation}\label{20}
[\alpha , \beta]^{\sharp} = [\alpha^{\sharp}, \beta^{\sharp}] , \ \forall \alpha , \beta \in \Gamma \left(A\right).
\end{equation}

\end{lemma}

\begin{remarkth}
\rm
Eq. (\ref{20}) is often considered as a part of the definition of a Lie algebroid though it is a consequence of the other conditions.\\
\end{remarkth}

An important remark is that the Lie algebra structure on sections is of local type i.e. $[\alpha , \beta ] \left(x\right)$ will depend on $\beta$ (therefore, on $\alpha$ too) around $x$ only, $\forall x \in M$. As a consequence, the restriction of a Lie algebroid over $M$ to a open subset of $M$ is again a Lie algebroid. 

Now, we will give some examples of Lie algebroids

\begin{example}\label{22}
\rm
Any Lie algebra is a Lie algebroid over a single point. Indeed, identifying $\Gamma \left(\mathfrak{g} \right)$ with $\mathfrak{g}$, the Lie bracket on sections is simply the Lie algebra bracket and the anchor map is the trivial one.\\\\
This kind of Lie algebroid is a particular case of the following example.
\end{example}

\begin{example}\label{23}
\rm
Let $\left(A \rightarrow M, \sharp  , [ \cdot , \cdot]\right)$ be a Lie algebroid where $\sharp \equiv 0$. Then, the Lie bracket on $\Gamma \left(A\right)$ is a point-wise Lie bracket, that is, the restriction of $[\cdot , \cdot]$ to the fibres induces a Lie algebra structure on each of them. These kind of Lie algebroids (with $\sharp \equiv 0$) are called \textit{Lie algebra bundles}. Note that the Lie algebra structures on the fibres are not necessary isomorphic to each other.

\end{example}

\begin{example}
\rm
If $M$ is a smooth manifold, then the tangent bundle of $M$, $TM$, is a Lie algebroid: the anchor map is the identity map and the Lie bracket is the usual Lie bracket of vector fields. This is called the \textit{tangent algebroid} of $M$.
\end{example}

\begin{example}\label{25}
\rm
Let $\tau: P \rightarrow M$ be a principal bundle with structure group $G$. Denote by $\phi : G \times P \rightarrow P$ the action of $G$ on $P$. Now, suppose that $\left(A \rightarrow P, \sharp , [\cdot , \cdot]\right)$ is a Lie algebroid, with vector bundle projection $\pi:A \rightarrow P$ and that $\overline{\phi} : G \times A \rightarrow A$ is an action of $G$ on $A$ such that $ \pi$ is a vector bundle action under the action $\overline{\phi}$ where for each $g \in G$, the pair $\left( \overline{\phi}_{g} , \phi_{g}\right)$ satisfies that
\begin{itemize}
\item[(1)] $\sharp \circ  \overline{\phi}_{g} = T \phi_{g} \circ \sharp$.
\item[(2)] $[  \overline{\phi}_{g} \circ \alpha \circ \phi_{g}^{-1},  \overline{\phi}_{g} \circ \beta \circ \phi_{g}^{-1} ] =  \overline{\phi}_{g} \circ [\alpha , \beta] \circ \phi_{g}^{-1}, \ \forall \alpha , \beta \in \Gamma \left(A\right)$.
\end{itemize}
This fact will be equivalent to the fact of that $\left( \overline{\phi}_{g} , \phi_{g}\right)$ is a Lie algebroid isomorphism. Let $\overline{\pi}: A/G \rightarrow M$ be the quotient vector bundle of $\pi$ by the action of $G$. Then, we can construct a Lie algebroid structure on $\overline{\pi}$ such that:
\begin{itemize}
\item[(i)] The anchor map $\overline{\sharp} : A/G \rightarrow TM$ is given by
$$ \overline{\sharp} \left(u\right)= T_{\pi \left(a\right)}\tau \left( \sharp\left(a\right)\right),$$
for all $u \in A/G$ and $a \in A$, where $\overline{\tau}: A \rightarrow A/G $ the quotient projection and $\overline{\tau}\left(a\right)=u$.\\
\item[(ii)] Taking into account $\left( 2 \right)$, the Lie bracket on $\Gamma \left(A\right)$ restricts to $\Gamma \left(A\right)^{G}$, i.e., the $\overline{\phi}-$invariant sections of $\pi$. Then, it is easy to prove that
$$\Gamma \left(A\right)^{G} \cong \Gamma \left(A /G\right).$$
Hence, the induced a Lie algebra structure on $\Gamma\left(A/G\right)$ coincides just with our Lie bracket of the Lie algebroid structure on $\overline{\pi}$.
\end{itemize}
This kind of Lie algebroids are called \textit{quotient Lie algebroids by the action of a Lie group}.
\end{example}
A particular but interesting example of this construction is obtained when we consider the tangent lift of a free and proper action of a Lie group on a manifold.
\begin{example}\label{26}
\rm
Let $\pi: P \rightarrow M$ be a principal bundle with structure group $G$. Denote by $\phi $ the action of $G$ on $P$. Let $\left(TP \rightarrow P, Id_{TP} , [ \cdot , \cdot]\right)$ be the tangent algebroid and $\phi^{T} : G \times TP \rightarrow TP$ be the tangent lift of $\phi$.\\
Then, $\phi^{T}$ satisfies the conditions of Example \ref{25}. Thus, one may consider the quotient Lie algebroid $\left( TP/G \rightarrow M, \overline{\sharp}, \overline{[\cdot ,\cdot]}\right)$ by the action of $G$. This algebroid is called the \textit{Atiyah algebroid associated with the principal bundle} $\pi : P \rightarrow M$.\\
Note that, as we have seen, the space of sections can be considered as the space of invariant vector fields by the action $\phi$ over $M$.
\end{example}

\begin{example}\label{51}
\rm
Let $M$ be a manifold and $\frak g $ be a Lie algebra. We can construct a Lie algebroid structure over the vector bundle $A = TM \oplus \left( M \times \frak g\right) \rightarrow M$ such that
\begin{itemize}
\item[(i)] The anchor $\sharp :TM \oplus \left( M \times \frak g \right) \rightarrow TM$ is the projection.

\item[(ii)] Lie algebra structure over the space of sections is given by:
$$ [ X \oplus f , Y \oplus g ] =[X,Y] \oplus \{ X \left( g \right) - Y \left( f \right)  + [ f,g] \},$$
for all $  X \oplus f , Y \oplus g \in \Gamma \left(A \right)$.
\end{itemize}
This Lie algebroid is called the \textit{Trivial Lie algebroid on $M$ with structure algebra $\frak g$}.

\end{example}

Next, we introduce the definition of a morphism in the category of Lie algebroids. The main problem is that a morphism between vector bundles does not, in general, induce a map between the modules of sections, so it is not immediately clear what should be meant by bracket relation. We will give a direct definition in terms of $\left(\Phi,\phi\right)-$decompositons of sections which is easy to use, and is amenable to categorical methods.\\
Let $\Phi : A' \rightarrow A$, $\phi : M' \rightarrow M$ be a vector bundle morphism between $\pi : A \rightarrow M$ and $\pi' : A' \rightarrow M'$. We know that for each $ \alpha' \in \Gamma \left(A'\right)$, there exists $f_{i} \in \mathcal{C}^{\infty} \left(M'\right)$ and $\alpha_{i} \in \Gamma \left(A\right)$ such that
$$\Phi \circ \alpha' = \sum_{i=1}^{k} f_{i} \left(\alpha_{i} \circ \phi \right).$$
Thus, we are ready to give the definiton of Lie algebroid morphism. 
\begin{definition}
\rm
Let $\left(A \rightarrow M , \sharp, [\cdot ,\cdot ]\right)$, $\left(A' \rightarrow M' , \sharp', [\cdot,\cdot]'\right)$ be Lie algebroids. A morphism of Lie algebroids is a vector bundle morphism $\Phi: A' \rightarrow A$, $\phi : M' \rightarrow M$ such that
\begin{equation}
\sharp \circ \Phi = T\phi \circ \sharp',
\end{equation}
and such that for arbitrary $\alpha', \beta' \in \Gamma \left(A'\right)$ with $\left(\Phi,\phi\right)-$decompositions
$$\Phi \circ \alpha' = \sum_{i=1}^{k} f_{i}\left(\alpha_{i} \circ \phi\right),
\Phi \circ \beta' = \sum_{j=1}^{k} g_{j}\left(\beta_{j} \circ \phi\right),
$$
we have
\begin{small}
\begin{equation}\label{28}
\Phi \circ [\alpha' , \beta' ] = \sum_{i,j=1}^{k}f_{i}g_{j}\left([ \alpha_{i} , \beta_{j} ] \circ \phi \right) + \sum_{j=1}^{k}{\alpha'}^{\sharp'}\left(g_{j}\right)\left(\beta_{j}\circ \phi\right) - \sum_{i=1}^{k}{\beta'}^{\sharp'}\left(f_{i}\right)\left(\alpha_{i}\circ \phi\right).
\end{equation}
\end{small}
In fact, the right-hand side of Eq. (\ref{28}) is independent of the choice of the $\left(\Phi,\phi\right)-$decompositions of $\alpha'$ and $\beta'$.\\
\end{definition}
It is easy to prove that the composition preserves Lie agebroid morphisms and, hence, we can define the category of Lie algebroids.

\begin{remarkth}\label{29}
\rm
In particular, if $\alpha' \sim_{\left(\Phi , \phi\right)} \alpha$ and  $\beta' \sim_{\left(\Phi , \phi\right)} \beta$, then Eq. (\ref{28}) reduces to
$$\Phi \circ [\alpha', \beta'] = [\alpha , \beta]\circ \phi.$$

On the other hand, if $M=M'$ and $\phi = Id_{M}$ then Eq. (\ref{28}) reduces to 
$$ \Phi \circ [\alpha' , \beta'] = [\Phi \circ \alpha' , \Phi \circ \beta' ], \ \forall \alpha' , \beta' \in \Gamma \left(A'\right).$$ 
\end{remarkth}
Next, we will introduce the notion of Lie subalgebroid.
\begin{definition}
\rm
Let $\left(A \rightarrow M, \sharp , [ \cdot , \cdot ]\right)$ be a Lie algebroid. Suppose that $A'$ is an embedded submanifold of $A$ and $M'$ is a inmersed submanifold of $M$ with inclusion maps $i_{A'} : A' \hookrightarrow A$ and $i_{M'} : M' \hookrightarrow M$. $A'$ is called a \textit{Lie subalgebroid of $A$} if $A'$ is a Lie algebroid on $M'$ which is a vector subbundle of $\pi_{|M'}$, where $\pi : A \rightarrow M$ is the projection map of $A$, such that the inclusion is a morphism of Lie algebroids. A \textit{reduced subalgebroid of $A$} is a transitive Lie subalgebroid with $M$ as the base manifold.
\end{definition}
\begin{remarkth}
\rm
Suppose that $M' \subseteq M$ is a closed submanifold then, using the $\left(i_{A'},i_{M'}\right)-$ decomposition and extending functions, it satisfies that for all $\alpha' \in \Gamma \left(A'\right)$ there exists $\alpha \in \Gamma \left(A\right)$ such that
$$ i_{A'} \circ \alpha' = \alpha \circ i_{M'}.$$ 
So, Eq. (\ref{28}) reduces to 
$$ i_{A'} \circ [\alpha' , \beta']_{M'} = [ \alpha ,  \beta ]_{M} \circ i_{M'}, \ \forall \alpha' , \beta' \in \Gamma \left(A'\right).$$ 
\end{remarkth}
\begin{example}
\rm
Let $\left(A \rightarrow M , \sharp , [ \cdot , \cdot ]\right)$ be a Lie algebroid over $M$. Then from Lemma \ref{19} and Remark \ref{29} we deduce the anchor map $\sharp : A \rightarrow TM$ is a Lie algebroid morphism from $A$ to the tangent algebroid of $M$.
\end{example}
\begin{example}
\rm
Let $\phi : M_{1} \rightarrow M_{2}$ be a smooth map. Then $\left(T \phi , \phi \right) $ is a Lie algebroid morphism between the tangent algebroids $TM_{1}$ and $TM_{2}$.
\end{example}
\begin{example}
\rm
Let $\tau: P \rightarrow M$ be a principal bundle with structure group $G$ and $\left(A \rightarrow P, \sharp , [\cdot , \cdot]\right)$ be a Lie algebroid (with vector bundle projection $\pi:A \rightarrow P$) in the conditions of Example \ref{25}. If $\overline{\pi} : A / G \rightarrow M$ is the quotient Lie algebroid by the action of the Lie group $G$ then $\left( \overline{\tau}  , \tau\right)$ is a Lie algebroid morphism. 
\end{example}

Next, as an important example of Lie algebroid, we will explain briefly how to associate a Lie algebroid to a Lie groupoid (which generalize the construction of the Lie algebra associated to a Lie group). Fix $\Gamma \rightrightarrows M$ a Lie groupoid with structure maps $\alpha$, $\beta$, $\epsilon$ and $i$. A vector field $X \in \mathfrak{X} \left(\Gamma\right)$ is said to be a \textit{left-invariant vector field on $\Gamma$} if it satisfies the following two properties:
\begin{itemize}
\item[(a)] $X$ is tangent to $\beta^{-1}\left(x\right)$, for all $x \in M$.
\item[(b)] For each $g \in \Gamma$, the left translation $L_{g}$ preserves $X$.
\end{itemize}
We denote the space of left-invariant vector fields on $\Gamma$ by $\mathfrak{X}_{L} \left(\Gamma\right)$.\\
Similarly to the case of Lie groups, it is clear that the Lie bracket of two left-invariant vector fields is again a left-invariant vector field, say
\begin{equation}\label{52}
[ \mathfrak{X}_{L}\left(\Gamma\right) , \mathfrak{X}_{L}\left(\Gamma\right) ] \subset \mathfrak{X}_{L}\left(\Gamma\right).
\end{equation}
Now, we will describe the associated Lie algebroid to $\Gamma$:
\begin{itemize}
\item[(i)] The vector bundle $ \pi^{\epsilon} :A \Gamma \rightarrow M$ satisfy that for each $x \in M$, the fibre of $A \Gamma$ at $x$ is
$$ A \Gamma_{x} = T_{\epsilon \left( x \right) } \beta^{-1} \left(x \right).$$
\item[(ii)]  Let $\Lambda$ be a section of $A \Gamma$. Then, we can define the left-invariant vector field on $\Gamma$ given by
$$ X_{\Lambda}\left(g\right) = T_{\epsilon \left(\alpha \left(g\right)\right)} L_{g} \left( \Lambda \left(\alpha\left(g\right)\right)\right), \ \forall g \in \Gamma ,$$
i.e., $X_{\Lambda}$ is determinated by the following equality
$$ X_{\Lambda}\left( \epsilon \left(x\right) \right) =  \Lambda \left(x\right), \ \forall x \in M .$$
Conversely, if $X \in \mathfrak{X}_{L} \left(\Gamma\right)$, then $\Lambda^{X} =X \circ \epsilon : M \rightarrow T \Gamma$ induces a section of $A \Gamma$ and, indeed, the correspondence $\Lambda \mapsto X_{\Lambda}$ is an $\mathbb{R}-$bilinear isomorphism from $\Gamma \left(A \Gamma\right)$ to $\mathfrak{X}_{L}\left(\Gamma\right)$ with inverse $X \mapsto \Lambda^{X}$. With this identification  $\Gamma \left(A \Gamma\right)$ inherits a Lie bracket from $ \mathfrak{X}_{L} \left(\Gamma\right)$.\\


\item[(iii)] The anchor map $\sharp$ is defined by
\begin{equation}\label{56}
 \sharp \left( \Lambda \left(x\right)\right) =  T_{\epsilon \left(x\right)} \alpha \left( \Lambda \left(x\right)\right),
\end{equation}

for all $\Lambda \left( x \right) \in A \Gamma_{x}$ and $x \in M$. In fact, for all $\Lambda \in \Gamma \left( A \Gamma \right)$,
\begin{equation}\label{57}
\Lambda^{\sharp}  =  T \alpha \circ \Lambda.
\end{equation}
\end{itemize}
The Lie algebroid $\left( A \Gamma \rightarrow M , \sharp , [ \cdot , \cdot ]\right)$ is called the \textit{Lie algebroid associated to the Lie groupoid} $\Gamma \rightrightarrows M$, and sometimes denoted by $A \Gamma$.\\

\begin{remarkth}\label{58}
\rm
Let $\Gamma \rightrightarrows M$ be a Lie groupoid. For any $x \in M$, the associated Lie algebra to the isotropy Lie group $\Gamma^{x}_{x}$, $A \left(\Gamma^{x}_{x}\right)$ is isomorphic to the isotropy Lie algebra through $x$, i.e.,
\begin{equation}\label{59}
 A \left(\Gamma^{x}_{x}\right) \cong Ker \left( \sharp_{x}\right),
\end{equation}
where $\sharp_{x}$ is the restriction of $\sharp$ to the fibre at $x$ of $A \Gamma$.
\end{remarkth}

\begin{theorem}
There is a natural functor from the category of Lie groupoids to the category of Lie algebroids.

\begin{proof}
We will given an sketch of the proof (a detailed proof can found in \cite{PJHKM}).\\
We already have given the definition of the correspondence between objects ($\Gamma \rightrightarrows M \ \rightarrow A\Gamma$) and we will obtain the correspondence between morphisms.\\
Let $\left( \Phi , \phi\right) : \Gamma_{1} \rightrightarrows M_{1} \rightarrow \Gamma_{2} \rightrightarrows M_{2}$ be a Lie groupoid morphism, with $\Phi : \Gamma_{1} \rightarrow \Gamma_{2}$ and $\phi : M_{1} \rightarrow M_{2}$. Then, $\left( \Phi, \phi\right)$ induces a morphism of Lie algebroids from $A\Gamma_{1}$ to $A \Gamma_{2}$ given by $\left(\Phi_{*}, \phi\right)$ such that for each $v_{x} \in T_{\epsilon_{1}\left(x\right)} \beta_{1}^{-1} \left( x \right),$
\begin{equation}\label{60}
\Phi_{*} \left( v_{x} \right) \triangleq T_{\epsilon_{1}\left(x\right)} \Phi_{x} \left( v_{x} \right)
\end{equation}
where $\Phi_{x} : \beta_{1}^{-1}\left(x\right) \rightarrow \beta_{2}^{-1}\left(x\right)$ is the restriction of $\Phi$ to $\beta_{1}^{-1}\left(x\right)$ for each $x \in M$.
\end{proof}
\end{theorem}
This morphism induced by a morphism $\left(\Phi , \phi\right)$ of Lie groupoids over the associated Lie algebroids will be denoted by $A \Phi $.\\

Now, we will give some examples of the above general construction.
\begin{example}
\rm
Let $M$ be a smooth manifold and $M \times M  \rightrightarrows M$ be the pair groupoid (see Example \ref{7}). Then, the vector bundle $A \left( M \times M \right)$ can be seen as the tangent bundle $\pi_{M}: TM \rightarrow M$. It follows that the associated Lie algebroid to $M \times M \rightrightarrows M$ is the tangent algebroid.
\end{example}

\begin{example}
\rm
Let $M$ be a manifold and $G$ be a Lie group. Consider the trivial Lie groupoid on $M$ with group $G$ (see Example \ref{99}). Then, the associated Lie algebroid is the trivial Lie algebroid on $M$ with structure algebra $\frak g $ (see Example \ref{51}), i.e., $ TM \oplus \left( M \times \frak g\right) \rightarrow M$.\\
Notice that, using the Lie's third fundamental theorem, every trivial Lie algebroid on $M$ with structure algebra $\mathfrak{g}$ is the induced Lie algebroid of a trivial groupoid on $M$ with a group $G$.
\end{example}

\begin{example}\label{98}
\rm
Let $\pi: P \rightarrow M$ be a principal bundle with structure group $G$. Denote by $\phi :  P  \times G \rightarrow P$ the action of $G$ on $P$.\\
Now, suppose that $\Gamma \rightrightarrows P$ is a Lie groupoid, with $\overline{\phi} :  \Gamma \times G \rightarrow \Gamma$ a free and proper action of $G$ on $\Gamma$ such that, for each $g \in G$, the pair $\left( \overline{\phi}_{g} , \phi_{g}\right)$ is an isomorphism of Lie groupoids. So, we may construct the quotient Lie groupoid by the action of a Lie group, $\Gamma / G \rightrightarrows M$ (see Example \ref{13}).\\
Then, by construction, we may identify $A \left( \Gamma / G \right)$ with the quotient Lie algebroids by the action of a Lie group, $A \Gamma /G$ (see Example \ref{25}).
\end{example}
As a particular case, we may give the following interesting example.
\begin{example}\label{61}

\rm
Let $\pi: P \rightarrow M$ be a principal bundle with structure group $G$ and $gauge\left(P\right)$ be the gauge groupoid (see Example \ref{14}). Then, the associated Lie algebroid to $gauge \left(P\right)$ is the Atiyah algebroid associated with the principal bundle $\pi : P \rightarrow M$ (see Example \ref{26}).
\end{example}

\begin{example}\label{62}
\rm
Let $\Phi \left(A\right) \rightrightarrows M$ be the frame groupoid. Then $A \Phi \left(A\right)$ is called \textit{frame algebroid} (see Example \ref{15}). As a particular case $A \Pi^{1}\left(M,M\right)$ is called \textit{$1-$jets algebroid}.\\
Let $\left(x^{i}\right)$ be a local coordinate system defined on some open subset $U \subseteq M$, using Eq. (\ref{17}) we can consider local coordinates on $A \Pi^{1}\left(M,M\right)$ as follows
\begin{equation}\label{63}
A\Pi^{1}\left(U,U\right) : \left(\left(x^{i} , x^{i}, \delta^{i}_{j}\right),  v^{j}, 0 , v^{j}_{i}\right) \cong \left(x^{i} , v^{j}, v^{j}_{i}\right).
\end{equation}

\end{example}

\begin{example}\label{64}
\rm
Let $\tilde{J}^{1} \left(FM\right)$ be the second-order non-holonomic \linebreak{groupoid} over a manifold $M$ and $\left(x^{i}\right)$ be a local coordinate system on an open set $U \subseteq M$. Using Eq. (\ref{45}) we can construct induced local coordinates over $A \tilde{J}^{1}\left(FM\right)$ as follows:
\begin{small}
\begin{equation}\label{65}
A\tilde{J}^{1}\left(FU\right) : \left(\left(\left(x^{i}\right),\left( x^{i} , \delta^{i}_{j}\right),\delta^{i}_{j},0\right),v^{i}, v^{i}_{j},0 ,  v^{i}_{,j}, v^{i}_{j,k}\right) \cong  \left(x^{i},v^{i}, v^{i}_{j}, v^{i}_{,j}, v^{i}_{j,k}\right).
\end{equation}
\end{small}
$A \tilde{J}^{1}\left(FM\right)$ is called the \textit{second-order non-holonomic algebroid over $M$}. We will denote the anchor of this Lie algebroid by $\overline{\sharp}$.
\end{example}

\begin{example}\label{141}
\rm
Let $\tilde{j}^{1} \left(FM\right)$ be the second-order holonomic \linebreak{groupoid} over a manifold $M$ and $\left(x^{i}\right)$ be a local coordinate system on an open set $U \subseteq M$. Using the above example we can construct induced local coordinates over $A \tilde{j}^{1}\left(FM\right)$ as follows:
\begin{small}
\begin{equation}\label{142}
A\tilde{j}^{1}\left(FU\right) :   \left(x^{i},v^{i}, v^{i}_{j}, v^{i}_{j,k}\right), \ \ \ v^{i}_{j,k} = v^{i}_{k,j}
\end{equation}
\end{small}
$A \tilde{j}^{1}\left(FM\right)$ is called the \textit{second-order holonomic algebroid over $M$}.
\end{example}


\section{Derivation algebroid}
In this section, we will extend the notion of exponential map for Lie groups to the context of Lie groupoids and algebroids. Next, we are going to use the exponential map in order to get another way of interpreting the associated Lie algebroid of $\Pi^{1}\left(M,M\right)$. A more detailed construction of this algebroid can be found in \cite{KMG} (see also \cite{YKKM}).\\

Let $\Gamma \rightrightarrows M$ be a Lie groupoid and $\Lambda \in \Gamma \left( A \Gamma\right)$ be a section of $A \Gamma  $. We consider the left-invariant vector field associated to $\Lambda$, $X_{\Lambda} \in \mathfrak{X}_{L} \left( \Gamma \right)$, $\{ \varphi_{t} : \mathcal{U}_{t} \rightarrow \mathcal{U}_{-t} \}$ the local flow of $X_{\Lambda}$ and $\{ \psi_{t} : U_{t} \rightarrow U_{-t} \}$ the local flow of $\Lambda^{\sharp}$. Then, we have

\begin{itemize}
\item[(i)]  $ \beta \circ \varphi_{t} = \beta, \  \forall t$
\item[(ii)]  $\alpha \circ \varphi_{t} = \psi_{t} \circ \alpha, \  \forall t$
\end{itemize}

Then, we can define a differentiable map $ Exp_{t} \Lambda  : U_{t} \rightarrow \Gamma$ in the following way
$$ Exp_{t} \Lambda  \left( x\right) = i \left( \xi \right) \cdot  \varphi_{t} \left( \xi\right), \ \forall x \in  U_{t},$$
being $ \xi \in \mathcal{U}_{t} \cap \alpha^{-1} \left(x\right)$.\\
Observe that, for all $t$
\begin{itemize}
\item[(i)]  $\beta \circ Exp_{t} \Lambda  = Id_{U_{t}}.$
\item[(ii)] $\alpha \circ Exp_{t} \Lambda  = \psi_{t} \circ \alpha.$
\end{itemize}
Therefore, we have the following result.
\begin{proposition}\label{103}
Let $\Gamma \rightrightarrows M$ be a Lie groupoid, $ W \subseteq M$ be an open subset, and take $\Lambda \in \Gamma_{W} \left( A \Gamma \right) $. Then, there is a differentiable map $ Exp  \Lambda : \mathcal{D} \rightarrow \Gamma$ given by
$$ \left(t, x\right) \mapsto Exp_{t} \Lambda \left(x\right),$$
where $\mathcal{D}$ is the domain of the flow of $\Lambda^{\sharp}$, such that
\begin{itemize}
\item[(i)] $\dfrac{\partial}{\partial t_{|t=0}} \left( Exp_{t} \Lambda \right) = \Lambda.$
\item[(ii)] $Exp_{0}\Lambda = \epsilon_{U}.$
\item[(iii)] $\{ \alpha \circ Exp_{t} \Lambda : U \rightarrow U_{-t}\}$, is a local $1-$parameter group of transformations of $\Lambda^{\sharp}$
\end{itemize}
This map is the called \textit{exponential of $\Lambda$}.
\end{proposition}

Now, as a second step, we will have to construct a new Lie algebroid. To do this, first we have to introduce some notions.

\begin{definition}
\rm
Let $\pi : A \rightarrow M$ be a vector bundle. A \textit{derivation on $A$} is a  $\mathbb{R}-$linear map $D : \Gamma \left(A\right) \rightarrow \Gamma \left(A\right)$ with a vector field $X \in \frak X \left(M\right)$ such that for each $f \in \mathcal{C}^{\infty}\left(M\right)$ and $\Lambda \in \Gamma \left(A\right)$,
$$ D \left(f \Lambda\right) = f D \left( \Lambda \right) + X\left(f\right) \Lambda.$$
We call $X$ the \textit{base vector field of $D$}. So, a derivation on $A$ is characterized by two geometrical objects, $D$ and $X$.
\end{definition}
Let us give some examples of these objects.

\begin{example}
\rm
Let $\left(  A \rightarrow M, \sharp , [ \cdot , \cdot ] \right)$ be a Lie algebroid. For each $\Lambda \in \Gamma \left(A\right)$, the map
$$
\begin{array}{rccl}
[\Lambda , \cdot] : & \Gamma \left(A\right)  & \rightarrow & \Gamma \left(A\right) \\
&\Theta  &\mapsto &  [\Lambda , \Theta].
\end{array}
$$
is a derivation on $A$ with base vector field $\Lambda^{\sharp}$.
\end{example}

\begin{example}\label{47}
\rm
Let $\pi : A \rightarrow M$ be a vector bundle and $\nabla : \frak X \left(M\right) \times \Gamma \left(A\right) \rightarrow \Gamma \left(A\right)$ be a \textit{covariant derivative}, i.e., $\nabla$ is a $\mathbb{R}-$bilinear map such that,
\begin{itemize}
\item[(1)] It is $\mathcal{C}^{\infty} \left(M\right)-$linear in the first variable.
\item[(2)] For all $X \in \frak X \left( M \right)$, $ \Theta \in \Gamma \left(A\right)$ and $f \in \mathcal{C}^{\infty} \left(M\right)$,
\begin{equation}
\nabla_{X} f \Theta = f \nabla_{X}\Theta + X \left( f\right) \Theta.
\end{equation}
\end{itemize}

Then, any vector field $X \in \frak X \left(M\right)$ generates a derivation on $A$, $\nabla_{X}$, (with base vector field $X$) fixing the first coordinate again, i.e.,
$$ \nabla_{X} : \Gamma \left(A\right) \rightarrow \Gamma \left(A\right),$$
such that
$$ \nabla_{X}\left(\Theta\right) = \nabla_{X}\Theta, \ \forall \Theta \in \Gamma \left(A\right).$$

\end{example}

Now, the space of derivations on $A$ can be considered as the space of sections of a vector bundle $\mathfrak{D}\left(A\right)$ on $M$. We can endow this vector bundle with a Lie algebroid structure.
\begin{itemize}
\item Let $D_{1} , D_{2}$ be derivations on $A$, we can define $[D_{1} , D_{2}]$ as the commutator, i.e.,
$$[D_{1} , D_{2}] = D_{1} \circ D_{2} - D_{2}\circ D_{1}.$$
A simple computation shows that the commutator of two derivations is again a derivation, indeed, the base vector field of $[D_{1} , D_{2}]$ is given by 
\begin{equation}
[X_{1} , X_{2}],
\end{equation}
where $X_{1}$ and $X_{2}$ are the base vector fields of $D_{1}$ and $D_{2}$ respectively.
\item Let $D$ be a derivation on $A$, then $D^{\sharp}$ is its base vector field.
\end{itemize}
Thus, with this structure $\mathfrak{D}\left(A\right)$ is a transitive Lie algebroid called the \textit{Lie algebroid of derivations on $A$}. The space of sections of $\mathfrak{D}\left(A\right)$, the derivations on $A$, will be denoted by $Der\left(A\right)$.

Note that in this Lie algebroid the linear sections of $\sharp$ are \linebreak$\mathcal{C}^{\infty}\left(M\right)-$linear maps from $\frak X \left(M\right)$ to $Der\left(A\right)$. So, the space of linear sections of $\sharp$ is, indeed, the space of covariant derivatives on $M$. Conversely, it is easy to see that a covariant derivative $\nabla$ is a section (in the category of Lie algebroids) of $\sharp$ if, and only if, $\nabla$ is flat, i.e., for all $X, Y \in \frak X \left( M \right)$
\begin{equation}\label{150}
 R \left( X , Y \right) = \nabla_{[X,Y]} - \nabla_{X} \nabla_{Y} + \nabla_{Y} \nabla_{X} = 0.
\end{equation}
Finally, as a last step, it is turn to relate this algebroid with the frame algebroid. Let $\Phi \left(A\right) \rightrightarrows M$ be the frame groupoid of a vector bundle $A \rightarrow M$. Consider $\Lambda \in \Gamma \left( A\Phi \left(A\right) \right)$ and its exponential map $ Exp_{t} \Lambda : U_{t} \rightarrow \Phi \left(A\right) $. Then, we can define a (local) linear map $Exp_{t} \Lambda^{*} : \Gamma \left( A \right) \rightarrow \Gamma \left( A \right)$ satisfying
$$ \{Exp_{t} \Lambda^{*} \left( X\right) \} \left( x \right)  =  Exp_{t} \Lambda \left( x\right) \left( X \left(\left( \alpha \circ Exp_{t} \Lambda\right)\left(x\right)\right)\right),$$
for each $X \in \Gamma \left(A\right)$ and $x \in U_{t}$. 
Thus, we can give the following result,

\begin{theorem}\label{48}
Let $A$ be a vector bundle over $M$. We can consider a map $\mathcal{D} : \Gamma \left(A \Phi \left(A\right)\right) \rightarrow Der\left(A\right)$ given by
$$ \mathcal{D}\left( \Lambda \right) \triangleq D^{\Lambda}  = \dfrac{\partial}{\partial t_{|t=0}} \left( Exp_{t}\Lambda^{*}\right),$$
which define a Lie algebroid isomorphism $\mathcal{D} : A \Phi \left(A\right) \rightarrow  \mathfrak{D} \left(A\right)$ over the identity map on $M$.
\end{theorem}
More detailed, for each $ X \in \Gamma \left( A \right)$ and $x \in M$ we have
\begin{eqnarray*}
D^{\Lambda}X \left( x \right) &=& \dfrac{\partial}{\partial t_{|t=0}} \left( Exp_{t}\Lambda  \left(x \right)   \left( X \left( \alpha \circ Exp_{t} \Lambda \right) \left(x \right)\right)\right).
\end{eqnarray*}
Notice that, for all $f \in \mathcal{C}^{\infty} \left( M \right)$
\begin{small}
\begin{eqnarray*}
D^{\Lambda}f X \left( x \right) &=& \dfrac{\partial}{\partial t_{|t=0}} \left( Exp_{t}\Lambda  \left( x \right)   \left( f\left( \left( \alpha \circ Exp_{t} \Lambda \right) \left(x \right) \right)X \left(\left( \alpha \circ Exp_{t} \Lambda \right) \left( x \right)\right) \right)\right)\\
&=& \dfrac{\partial}{\partial t_{|t=0}} \left( f\left(\left( \alpha \circ Exp_{t} \Lambda \right) \left(x \right)\right) Exp_{t}\Lambda  \left(x \right)   \left( X \left(\left( \alpha \circ Exp_{t} \Lambda \right) \left(x \right)\right)\right)\right)\\
 &=& \Lambda^{\sharp} \left( x \right) \left( f \right) X \left( x \right) + f \left( x \right) D^{\Lambda}X \left( x \right).
\end{eqnarray*}

\end{small}
This theorem gives us another way of interpreting the $1-$jets Lie algebroid. In fact, the $1-$jets Lie algebroid $A \Pi^{1} \left( M , M \right)$ is isomorphic to the algebroid of derivations on $TM$.

Notice that using this isomorphism, we can consider a one-to-one map from linear sections of $\sharp$ in $ A \Pi^{1}\left(M,M\right)$ to covariant derivatives over $M$. Thus, having a section $\Lambda$ of $\sharp$ in $A \Pi^{1}\left(M,M\right)$ we will denote its associated covariant derivative by $\nabla^{\Lambda}$. Furthermore, $\Lambda$ is a morphism of Lie algebroids if, and only if, $\nabla^{\Lambda}$ is flat.

Let $\left( x^{i} \right)$ and $\left(y^{j}\right)$ be coordinate systems on $M$, then we can induce local coordinates over $\Pi^{1}\left(M,M\right)$ and $A \Pi^{1}\left(M,M\right)$, i.e.,
\begin{itemize}
\item $\Pi^{1}\left(M,M\right) : \left(x^{i}, y^{j} , y^{j}_{i}\right).$

\item $A\Pi^{1}\left(M,M\right) : \left(x^{i} , v^{j}, v^{j}_{i}\right).$
\end{itemize}
With these (local) coordinates we can give a result which can help us to understand the shape of the above theorem.
\begin{lemma}\label{76}
Let $M$ be a manifold and $\Lambda $ be a section of the $1-$jets algebroid with local expression
$$\Lambda \left(x^{i}\right) = \left(x^{i} , \Lambda^{j} , \Lambda^{j}_{i}\right).$$
The matrix $\Lambda^{j}_{i}$ is (locally) the associated matrix to $D^{\Lambda}$, i.e.,
$$D^{\Lambda} \left( \dfrac{\partial}{\partial x^{i}} \right)= \sum_{j} \Lambda^{j}_{i}  \dfrac{\partial}{\partial x^{j}},$$
and the base vector field of $D^{\Lambda}$ is $\Lambda^{\sharp}$ which is given locally by $\left(x^{i} ,  \Lambda^{j} \right)$.
\begin{proof}
Let $\Lambda \in \Gamma \left(A\Pi^{1}\left(M,M\right)\right)$ be a section of the $1-$jets algebroid and $X_{\Lambda}$ its associated left-invariant vector field over $\Pi^{1}\left(M,M\right)$. Considering the flow of $X_{\Lambda}$, $\{\varphi_{t} : \mathcal{U}_{t} \rightarrow \mathcal{U}_{-t}\}$ we have by construction that
$$ Exp_{t} \Lambda \left( x\right) = \xi^{-1} \cdot \varphi_{t} \left( \xi\right) , \ \forall x \in  \alpha \left( \mathcal{U}_{t}\right),$$
where $ \xi \in \mathcal{U}_{t} \cap \alpha^{-1} \left(x\right)$.\\
Now, let us take local coordinate systems $\left(x^{i}\right)$ and $\left(y^{j}\right)$ and its induced local coordinates over $\Lambda$, then
$$\Lambda \left(x^{i}\right) = \left(x^{i} , \Lambda^{j} , \Lambda^{j}_{i}\right).$$
Thus, the associated left-invariant vector field is (locally) as follows
$$X_{\Lambda} \left(x^{i},y^{j}, y^{j}_{i} \right) = \left(\left(x^{i},y^{j}, y^{j}_{i} \right) ,\Lambda^{j} ,0, y^{j}_{l}\cdot \Lambda^{l}_{i}\right).$$
Therefore, its flow
$$\varphi_{t} \left(x^{i},y^{j}, y^{j}_{i} \right) = \left( \psi_{t} \left(x^{i}\right) ,y^{j} , y^{j}_{l} \cdot \overline{\varphi}_{t}\left(x^{i}\right)\right),$$
satisfies that
\begin{itemize}
\item[(i)] $\psi_{t}$ is the flow of $\Lambda^{\sharp}$.
\item[(ii)] $ \dfrac{\partial}{\partial t_{|t=0}}  \left(\overline{\varphi}_{t} \left(x^{i}\right)\right) = \Lambda^{j}_{i}$.
\end{itemize}
So,
$$ Exp_{t} \Lambda \left( x^{i} \right) = \left( \psi_{t} \left(x^{i}\right) , x^{i} , \overline{\varphi}_{t}\left(x^{i}\right)\right).$$
Then, 
\begin{eqnarray*}
Exp_{t}\Lambda^{*} \left(\dfrac{\partial}{\partial x^{k}_{|x^{i}}} \right)&=& \left( \psi_{t} \left(x^{i}\right) , x^{i} , \overline{\varphi}_{t}\left(x^{i}\right)\right) \left( \dfrac{\partial}{\partial x^{k}_{| \psi_{t} \left( x^{i} \right)}}\right)\\
&=& \overline{\varphi}_{t}\left(x^{i}\right)  \dfrac{\partial}{\partial x^{k}_{|x^{i}}}
\end{eqnarray*}
Hence,
$$D^{\Lambda}  \dfrac{\partial}{\partial x^{k}} = \sum_{j} \Lambda^{j}_{k}  \dfrac{\partial}{\partial x^{j}},$$
i.e., the matrix $\Lambda^{j}_{i}$ is (locally) the associated matrix to $D^{\Lambda}$. 
\end{proof}
\end{lemma}
Let $\Lambda$ be a linear section of $\sharp$ in $A \Pi^{1}\left(M,M\right)$. Then, $\mathcal{D}$ induces a covariant derivative on $M$, $\nabla^{\Lambda}$. Thus, for each $\left(x^{i}\right)$ local coordinate system on $M$
$$\Lambda \left(x^{i} , \dfrac{\partial}{\partial x^{j}}\right) = \left(x^{i} , \dfrac{\partial}{\partial x^{j}}, \Lambda^{j}_{i}\right),$$
where $\Lambda^{j}_{i}$ depends on $\dfrac{\partial}{\partial x^{j}}$. Taking into account that $\Lambda^{j}_{i}$ is linear in the second coordinate we will change the notation as follows
\begin{equation}\label{145}
\Lambda^{j}_{i} \left(x^{l} , \dfrac{\partial}{\partial x^{k}}\right) \triangleq \Lambda^{j}_{i,k} \left( x^{l} \right).
\end{equation}
Therefore, locally $\Lambda$ will be written in the following way
$$\Lambda \left(x^{i} , \dfrac{\partial}{\partial x^{j}}\right) = \left(x^{i} , \dfrac{\partial}{\partial x^{j}}, \Lambda^{k}_{i,j}\right).$$
and thus
$$ \nabla^{\Lambda}_{\dfrac{\partial}{\partial x^{j}}} \dfrac{\partial}{\partial x^{i}} = D^{\Lambda \left( \dfrac{\partial}{\partial x^{j}}\right)} \dfrac{\partial}{\partial x^{i}} =\sum_{k} \Lambda^{k}_{i,j}  \dfrac{\partial}{\partial x^{k}},$$
where $\Lambda \left( \dfrac{\partial}{\partial x^{j}}\right)$ is the (local) section of $A\Pi^{1}\left(M,M\right)$ given by
$$ \Lambda \left( \dfrac{\partial}{\partial x^{j}}\right) \left(x\right) = \Lambda\left(x\right) \left( \dfrac{\partial}{\partial x^{j}_{|x}}\right).$$

So, $\Lambda^{k}_{i,j}$ are just the Christoffel symbols of $\nabla^{\Lambda}$.\\


Thus, our goal is to use this isomorphism to give another interpretation of the second-order non-holonomic algebroid. Then, as a first approximation to the second-order non-holonomic case we are going to restrict the isomorphism constructed in Theorem \ref{48} to a particular case.\\

Consider the $1-$jets groupoid on $FM$, $\Pi^{1} \left(FM, FM\right) \rightrightarrows FM$ and $J^{1} \left( FM \right) \rightrightarrows FM$ the Lie subgroupoid of all $1-$jets of local automorphisms on $FM$.\\
Let $\left(x^{i}\right)$ and $\left(y^{j}\right)$ be local coordinate systems over two opens $U,V \subseteq M$; then, the induced coordinate systems over $FM$ are denoted by
$$FU: \left(x^{i} , x^{i}_{j}\right)$$
$$FV: \left(y^{j} , y^{j}_{i}\right).$$
Hence, (see Eq. (\ref{101})) we can construct induced coordinates over $J^{1}\left(FM\right)$,
$$J^{1}\left(FU,FV\right) : \left(\left(x^{i},x^{i}_{j}\right),\left( y^{j} , y^{j}_{i}\right),y^{j}_{,i},y^{j}_{i,k}\right) .$$

So, we can consider its associated Lie algebroid $AJ^{1}\left(FM\right)$ as a reduced subalgebroid of the $1-$jets algebroid $A \Pi^{1}\left(FM,FM\right)$ and, hence, its induced coordinates will be
\begin{equation}\label{133}
AJ^{1}\left(FU\right) : \left(\left(x^{i},x^{i}_{j}\right),\left( v^{j} , v^{j}_{i}\right),v^{j}_{,i}, 0,v^{j}_{i,k},v^{j}_{i,kl}\right) \cong 
\end{equation}
$$ \cong \left(\left(x^{i},x^{i}_{j}\right),\left( v^{j} , v^{j}_{i}\right),v^{j}_{,i},v^{j}_{i,k}\right),$$
where,
$$ v^{j}_{i,kl} = \left(\sum_{m}v^{j}_{m} \left( x^{-1}\right)^{m}_{k}\right)\delta_{l}^{i}.$$
In this way, we can restrict the isomorphism given in the Theorem \ref{48} to get another isomophism between Lie algebroids, $AJ^{1}\left(FM\right) \rightarrow \mathfrak{D}^{1}\left(FM\right)$ where $\mathfrak{D}^{1}\left(FM\right)$ is the resulting Lie algebroid from the restriction of the isomorphism.\\
Let $\Lambda$ be a section of $A J^{1}\left(FM\right)$ such that (locally)
\begin{equation}\label{134}
\Lambda \left( x^{i},  x^{i}_{j}\right) = \left(\left(x^{i},x^{i}_{j}\right),\left( \Lambda^{j} , \Lambda^{j}_{i}\right),\Lambda^{j}_{,i},\Lambda^{j}_{i,k}\right)
\end{equation}
Then, the associated derivation is characterized by the following identities
\begin{itemize}
\item[(i)] $D^{\Lambda}  \dfrac{\partial}{\partial x^{i}} = \sum_{j} \Lambda^{j}_{,i} \dfrac{\partial}{\partial x^{j}} + \sum_{j,k} \Lambda^{j}_{k,i} \dfrac{\partial}{\partial x^{j}_{k}}.$

\item[(ii)] $D^{\Lambda}  \dfrac{\partial}{\partial x^{i}_{j}} = \sum_{k,m} \left(\Lambda^{k}_{m} \left(x^{-1}\right)^{m}_{i} \right)\dfrac{\partial}{\partial x^{k}_{j}}.$
\end{itemize}
So, conditions $\left(i\right)$ and $\left(ii\right)$ characterize the sections of Lie algebroid $\mathfrak{D}^{1} \left(FM\right)$. This space will be denoted by $Der^{1}\left(FM\right)$.\\
\begin{remarkth}\label{107}
\rm
We can characterize $Der^{1}\left(FM\right)$ in the following way. Let $\{ \left( \Phi_{t} , \psi_{t} \right) \}$ be the flow of $D^{\Lambda}$. Then,
$$ D^{\Lambda} X = \dfrac{\partial \Phi_{t}^{*} X}{\partial t_{|0}}, \ \forall X \in \mathfrak{X} \left( FM \right).$$
Hence, by uniqueness,
$$ Exp_{t} \Lambda^{*}   \left( X \right) = \Phi^{*}_{t} X,$$
i.e.,
$$ Exp_{t} \Lambda \left( x \right) = \Phi_{-t}.$$
Thus, we can say that $D^{\Lambda} \in Der^{1}\left(FM\right)$ if, and only if, its flow is the tangent map of an automorphism of frame bundles (over the identity map) at each fibre.

\end{remarkth}

Finally, we will work with the second-order non-holonomic algebroid $A \tilde{J}^{1} \left( FM \right)$ (see Example \ref{64}). As we know (see Example \ref{98}), $A\tilde{J}^{1}\left(FM\right)$ can be seen as the quotient Lie algebroid by the induced action of $\Phi$ over $AJ^{1}\left(FM\right)$.\\
In this way we can consider a relation in $ \mathfrak{D}^{1}\left(FM\right)$ given by the restriction of the isomorphism defined in the before section $ \mathcal{D}: AJ^{1}\left(FM\right) \rightarrow \mathfrak{D}^{1}\left(FM\right)$ and the relation in $AJ^{1}\left(FM\right)$, i.e.,
$$ \mathcal{D}\left(a\right) \sim \mathcal{D}\left(b\right) \Leftrightarrow a \sim b, \ \forall a,b \in AJ^{1}\left(FM\right).$$
The new quotience space is denoted by $ \tilde{\mathfrak{D}}^{1} \left(FM\right)$ and it is obvious that this space inherit the Lie algebroid structure from $A\tilde{J}^{1} \left(FM\right)$. In fact, considering $\tilde{\mathcal{D}}: A\tilde{J}^{1} \left(FM\right) \rightarrow \tilde{\mathfrak{D}}^{1} \left(FM\right)$ the map which commutes with the projections, the Lie algebroid structure over $\tilde{\mathfrak{D}}^{1} \left(FM\right)$ is the unique Lie algebroid structure such that $\tilde{\mathcal{D}}$ is a Lie algebroid isomorphism over the identity map on $M$. This Lie algebroid will be called \textit{second-order non-holonomic algebroid of derivations on $TM$}.\\ 
Let $\left(x^{i}\right)$ be a local coordinate system on an open set $U \subseteq M$. Using Eq. (\ref{65}) we can construct induced local coordinates over $A \tilde{J}^{1}\left(FM\right)$ as follows:
$$A\tilde{J}^{1}\left(FU\right)  \cong  \left(x^{i},v^{i}, v^{i}_{j}, v^{i}_{,j}, v^{i}_{j,k}\right).$$
Thus, the non-holonomic second-order derivation algebroid is characterized by the following equalities:

\begin{itemize}
\item[(i)] $D \dfrac{\partial}{\partial x^{i}} = \sum_{j} f^{j}_{,i} \dfrac{\partial}{\partial x^{j}} + \sum_{j,k} f^{j}_{k,i} \dfrac{\partial}{\partial x^{j}_{k}}$

\item[(ii)] $D  \dfrac{\partial}{\partial x^{i}_{j}} = \sum_{k} f^{k}_{i}\dfrac{\partial}{\partial x^{k}_{j}}$
\end{itemize}
\noindent where the local fuctions $f^{j}_{,i}, f^{j}_{k,i}$ and $f^{k}_{i}$ do not depend on $x^{i}_{j}$.\\
Observe that, we could restrict $\mathcal{D}$ to $Aj^{1} \left( FM \right)$ and we obtain a Lie subalgebroid of $ \mathfrak{D}^{1} \left(FM\right)$ which is denoted by $ \mathfrak{d}^{1} \left(FM\right)$. Proceeding in the same way as in the case of $\tilde{J}^{1} \left( FM \right)$, we obtain a reduced Lie subalgebroid of $\tilde{\mathfrak{D}}^{1}  \left( FM \right)$. This Lie algebroid is denoted by $\tilde{\mathfrak{d}}^{1} \left( FM \right)$ and it is called \textit{second-order holonomic algebroid of derivations on $TM$}. Obviously, this subalgebroid is isomorphic to holonomic algebroid of second order $A\tilde{j}^{1} \left( FM \right)$ by restricting $\tilde{\mathcal{D}}$.\\

The second-order holonomic algebroid of derivations on $TM$ is characterized by the above equalities satisfying additionally that
$$ f^{j}_{i} = f^{j}_{,i} \ \ \ \ \ ; \ \ \ \ \ f^{j}_{i,k} = f^{j}_{k,i}.$$

\begin{remarkth}\label{114}
\rm
Denote by $\Gamma \left( AJ^{1}\left(FM\right)\right)^{G}$ the set of $\left(A \Phi , R \right)-$invariant sections of $AJ^{1}\left(FM\right)$, i.e., for all $\Lambda \in \Gamma \left(AJ^{1}\left(FM\right)\right)^{G}$ and $g \in Gl\left(n, \mathbb{R}\right)$, the diagram

\vspace{0.5cm}
\begin{picture}(375,50)(50,40)
\put(180,20){\makebox(0,0){$FM $}}
\put(230,25){$\Lambda$}               \put(200,20){\vector(1,0){80}}
\put(310,20){\makebox(0,0){$AJ^{1}\left(FM\right)$}}
\put(145,50){$R_{g}$}                  \put(180,70){\vector(0,-1){40}}
\put(310,50){$A \Phi_{g}$}                  \put(300,70){\vector(0,-1){40}}
\put(180,80){\makebox(0,0){$  FM$}}
\put(230,85){$\Lambda$}               \put(200,80){\vector(1,0){80}}
\put(310,80){\makebox(0,0){$AJ^{1}\left(FM\right)$}}
\end{picture}

\vspace{35pt}
\noindent is commutative, namely
$$ \Lambda \left(X \cdot g\right) = T_{\epsilon \left(X\right) } \Phi_{g}^{X} \left(\Lambda \left(X\right)\right), \ \forall X \in FM, \ \forall g \in Gl\left(n, \mathbb{R}\right),$$
where $\Phi^{X}_{g} : \beta^{-1} \left(X\right) \rightarrow \beta^{-1} \left( X \cdot g \right)$ is the restriction of $\Phi_{g}$ to  $\beta^{-1} \left(X\right)$.\\

As we have seen in Example \ref{25}, the space of $\left(A \Phi , R \right)$-invariant sections of $AJ^{1}\left(FM\right)$ is isomorphic to $\Gamma \left( A \tilde{J}^{1} \left( FM \right) \right)$.\\
Next, take $\Lambda \in \Gamma \left( A J^{1} \left( FM \right) \right)^{G}$. Then, $Exp_{t} \Lambda$ satisfies
\begin{equation}\label{106}
Exp_{t} \Lambda \circ R_{g} = \Phi_{g} \circ Exp_{t} \Lambda, \  \forall g \in Gl \left( n , \mathbb{R} \right).
\end{equation}
Thus, let $Exp_{t}\Lambda^{*} : \frak X \left( FM  \right) \rightarrow \frak X \left( FM  \right)$ be the induced linear map over $\frak X \left( FM  \right)$. Let be $\overline{X} \in \frak X \left( FM  \right)$ then, we have
$$\{ Exp_{t}\Lambda^{*} \circ T R_{g}^{*} \left( \overline{X}  \right) \} \left( X \right)$$
$$= Exp_{t}\Lambda \left( X  \right) \left(  \left(T_{\left( \alpha \circ Exp_{t}\Lambda \right)\left( X\right)\cdot g^{-1}}R_{g} \left( \overline{X} \left( \left( \alpha \circ Exp_{t}\Lambda \right)\left( X\right) \cdot g^{-1}\right)\right)\right) \right)$$
$$=T_{ X \cdot g^{-1} }R_{g} \left( Exp_{t}\Lambda \left( X \cdot g^{-1} \right)  \left( \overline{X} \left( \left( \alpha \circ Exp_{t}\Lambda \right)\left( X \cdot g^{-1}\right)\right)\right)\right)$$
$$= \{ T R_{g}^{*} \circ Exp_{t}\Lambda^{*} \left( \overline{X}  \right) \} \left( X \right),$$
for all $ g \in Gl \left( n , \mathbb{R} \right)$.\\♥
Conversely, suppose that $\Lambda \in \Gamma \left(A J^{1}\left(FM\right)\right)$ satisfies the above equality. Then, in a similar way, we can prove that $\Lambda$ is $\left(A \Phi , R \right)-$invariant.
Deriving this equality, we have that it is equivalent to
\begin{equation}\label{108}
 D^{\Lambda} \circ TR^{*}_{g} = TR_{g}^{*} \circ D^{\Lambda}.
\end{equation}
So, we have proved that the space $ \tilde{\mathfrak{D}}^{1} \left(FM \right)$ can be seen as the derivations in $\mathfrak{D}^{1} \left(FM \right)$ which commute with $TR_{g}^{*}$ for all $g \in Gl \left( n , \mathbb{R} \right)$.\\
Analogously, the second-order holonomic algebroid of derivations \linebreak$ \tilde{\mathfrak{d}}^{1} \left(FM \right)$ can be seen as the derivations in $\mathfrak{d}^{1} \left(FM \right)$ which commute with $TR_{g}^{*}$ for all $g \in Gl \left( n , \mathbb{R} \right)$.\\

Observe that, if $D^{\Lambda}$ satisfies Eq. (\ref{108}) then, its base vector field $\Lambda^{\sharp} \in \mathfrak{X} \left( FM \right)$ is right-invariant, .i.e.,
$$ TR_{g}^{*} \Lambda^{\sharp} = \Lambda^{\sharp}, \ \forall g \in  Gl \left( n , \mathbb{R} \right),$$
or, equivalently,
$$ T_{Z}R_{g} \left( \Lambda^{\sharp} \left( Z \right) \right) = \Lambda^{\sharp} \left( Z \cdot g \right),$$
for all $Z \in FM$. Thus, $\Lambda^{\sharp}$ is $\pi_{M}-$related with a (unique) vector field over $M$.\\
Let $\left( x^{i} \right)$ be local coordinates on $M$ and $ \tilde{\Lambda}$ be a section of $ A \tilde{J}^{1} \left( FM \right)$ which satisfies that
$$\tilde{\Lambda} \left( x^{i} \right) = \left(x^{i},\Lambda^{i}, \Lambda^{i}_{j}, \Lambda^{i}_{,j}, \Lambda^{i}_{j,k} \right).$$
Then, its associated $\left(A \Phi , R \right)$-invariant section $\Lambda$ of $ A J^{1} \left( FM \right)$ is given by
$$  \Lambda \left( x^{i}, x^{i}_{j} \right) = \left(x^{i}, x^{i}_{j} ,\Lambda^{i}, \Lambda^{i}_{l}x^{l}_{j}, \Lambda^{i}_{,j}, \Lambda^{i}_{j,k} \right).$$
Hence,
$$\Lambda^{\sharp} \left( x^{i}, x^{i}_{j} \right) = \left(x^{i}, x^{i}_{j} ,\Lambda^{i}, \Lambda^{i}_{l}x^{l}_{j}\right),$$
and its $\pi_{M}-$related vector field on $M$ is
$$\tilde{\Lambda}^{\overline{\sharp}} \left( x^{i} \right) = \left(x^{i} ,\Lambda^{i} \right).$$

Finally, suppose that $\tilde{\Lambda}$ is a linear section of $\overline{\sharp}$. Then, we can consider the map $\Lambda : \frak X \left( M  \right) \rightarrow \Gamma \left( AJ^{1}\left(FM\right)\right)^{G}$ such that for each vector field $X \in \frak X \left( M  \right)$, $\Lambda \left( X \right)$ is the associated $\left(A \Phi , R \right)$-invariant section of $ A J^{1} \left( FM \right)$ to $\tilde{\Lambda} \left( X \right)$.\\
Then, for all $f \in \mathcal{C}^{\infty} \left( M \right)$ and $X \in \frak X \left( M  \right)$, 
$$\Lambda \left(f X \right) = \left( f \circ \pi_{M} \right) \Lambda \left( X \right).$$
Hence, considering the associated derivation to $\Lambda \left( X \right)$ we obtain the following map

\begin{equation}\label{137}
\nabla^{\Lambda} : \frak X \left( M \right) \times \frak X \left( F M \right) \rightarrow \frak X \left( FM \right),
\end{equation}
which satisfies
\begin{itemize}
\item[(i)] For all $f \in \mathcal{C}^{\infty} \left( M \right)$, $X \in \frak X \left( M  \right)$ and $ \tilde{Y} \in \frak X \left( FM \right)$ we have
$$ \nabla^{\Lambda}_{fX}\tilde{Y} = \left( f \circ \pi_{M} \right) \nabla^{\Lambda}_{X}\tilde{Y}.$$
\item[(ii)] For all $F \in \mathcal{C}^{\infty} \left( FM \right)$, $X \in \frak X \left( M  \right)$ and $ \tilde{Y} \in \frak X \left( FM \right)$ we have
$$ \nabla^{\Lambda}_{X}F\tilde{Y} = F \nabla^{\Lambda}_{X}\tilde{Y} + \Lambda\left( X \right)^{\sharp} \left( F \right) \tilde{Y}.$$
\item[(iii)] For all $X \in \frak X \left( M  \right)$ the base vector field of $\nabla^{\Lambda}_{X}$ is $\Lambda \left( X \right)^{\sharp}$ which is $\pi_{M}-$related to $X$.
\item[(iv)] For all $g \in Gl \left( n , \mathbb{R} \right)$ and $X \in \frak X \left( M  \right)$,
$$\nabla_{X}^{\Lambda} \circ TR^{*}_{g} = TR_{g}^{*} \circ \nabla_{X}^{\Lambda}.$$
\item[(v)] For all $X \in \frak X \left( M \right)$ the flow of $\nabla^{\Lambda}_{X}$ is the tangent map of an automorphism of frame bundles (over the identity map) at each fibre.

\end{itemize}
These kind of objects will be called \textit{second-order non-holonomic covariant derivatives on $M$}.\\
Roughly speaking, the isomorphism $\tilde{\mathcal{D}}$ gives us a way to interpret a linear section of $\overline{\sharp}$ as a map which turn a vector field $X \in \frak X \left( M \right)$ into a $TR^{*}_{g}-$invariant derivation over $TFM$ with a base vector field which projects over $X$. Note that, in this case, this map is not exactly a covariant derivative but it has a similar shape.
\end{remarkth}

\section{Cosserat Media: Uniformity and Homogeneity}

A \textit{body} $\mathcal{B}$ is a three-dimensional differentiable manifold which can be covered with just one chart. An embedding $\phi : \mathcal{B} \rightarrow \mathbb{R}^{3}$ is called a \textit{configuration of} $\mathcal{B}$ and its $1-$jet $j_{x,\phi \left(x\right)}^{1} \phi$ at $x \in \mathcal{B}$ is called an \textit{infinitesimal configuration at $x$}. We usually identify the body with any one of its configurations, say $\phi_{0}$, called \textit{reference configuration}. Given any arbitrary configuration $\phi$, the change of configurations $\kappa = \phi \circ \phi_{0}^{-1}$ is called a \textit{deformation}, and its $1-$jet $j_{\phi_{0}\left(x\right) , \phi \left(x\right)}^{1} \kappa$ is called an \textit{infinitesimal configuration at $\phi_{0}\left(x\right)$}.\\
For elastic bodies, the material is completely characterized by one function $W$ which depends, at each point $x \in \mathcal{B}$, on the gradient of the deformation evaluated at the point. Thus, $W$ is defined (see \cite{MELZA}) as a differentiable map
$$ W : Gl\left(3 , \mathbb{R}\right) \times \mathcal{B} \rightarrow V,$$
where $V$ is a real vector space. Another equivalent way of considering $W$ is as a differentiable map
$$ W : \Pi^{1} \left( \mathcal{B}, \mathcal{B}\right) \rightarrow V,$$
which does not depend on the final point, i.e., for all $x,y,z \in \mathcal{B}$
$$W \left( j_{x,y}^{1} \phi\right) = W \left( j_{x,z}^{1} \left( \tau_{z-y} \circ \phi\right)\right), \ \forall j_{x,y}^{1}\phi \in \Pi^{1} \left( \mathcal{B}, \mathcal{B}\right),$$
where $\tau_{v}$ is the translation map on $\mathbb{R}^{3}$ by the vector $v$.
The picture describing a Cosserat medium is a little bit more difficult of describing. A \textit{Cosserat medium} is the linear frame bundle $F \mathcal{B}$ of a body $\mathcal{B}$. $\mathcal{B}$ is usually called the \textit{macromedium} or \textit{underlying body}. With some abuse of notation, we shall call $\mathcal{B}$ the \textit{Cosserat continuum}.\\
A \textit{configuration} of a Cosserat medium $F\mathcal{B}$ is an embedding $\Psi : F \mathcal{B} \rightarrow F \mathbb{R}^{3}$ of principal bundles such that the induced Lie group morphism $\tilde{\psi} : Gl \left( 3 , \mathbb{R}\right) \rightarrow Gl\left(3, \mathbb{R}\right)$ is the identity map. Hence $\Psi$ satisfies
$$ \Psi \left( \tilde{X} \cdot g \right) = \Psi \left(\tilde{X}\right) \cdot g, \ \forall \tilde{X} \in F \mathcal{B}, \ \forall g \in Gl \left( 3 , \mathbb{R} \right).$$
Also, $\Psi$ induces an embedding $\psi : \mathcal{B} \rightarrow \mathbb{R}^{3}$ verifying
$$ \pi_{\mathbb{R}^{3}} \circ \Psi = \psi \circ \pi_{\mathcal{B}}.$$
In particular, $\psi$ is a configuration of the macromedium $\mathcal{B}$.\\
Notice that the subbundle $\Psi \left( F \mathcal{B}\right)$ of $F \mathbb{R}^{3}$ is just the frame bundle of $\psi \left( \mathcal{B}\right)$, i.e.,
$$\Psi \left( F \mathcal{B}\right) = F \psi \left( \mathcal{B}\right).$$
Since we are dealing with equivariants embedding, we can consider equivalence classes of the $1-$jets $j^{1}_{\tilde{X} , \Psi \left( \tilde{X}\right)} \Psi$ according to the action (\ref{44}). So, the equivalence class of an $1-$jet $j^{1}_{\tilde{X} , \Psi \left( \tilde{X}\right)} \Psi$, which is denoted by $j^{1}_{x,\psi \left(x\right)} \Psi$ like in the non-holonomic groupoid of second order, is called \textit{infinitesimal configuration at $x$}. We usually identify the Cosserat medium with any one on its cofigurations, say $\Psi_{0} : F \mathcal{B} \rightarrow F \mathbb{R}^{3}$ and we denote by $\psi_{0}$ the induced map of $\Psi_{0}$. $\Psi_{0}$ is called \textit{reference configuration}. Given any configuration $\Psi$, the change of configuration $\tilde{\kappa} = \Psi \circ \Psi_{0}^{-1}$ is called a \textit{deformation}, and its class of $1-$jets $j^{1}_{\psi_{0} \left(x\right) , \psi \left(x\right)} \tilde{\kappa}$ is called an \textit{infinitesimal deformation} at $\psi_{0}\left(x\right)$. Notice that the induced map of $\tilde{\kappa}$, $\kappa = \psi \circ \psi_{0}^{-1}$, is a deformation on the body $\mathcal{B}$.\\
From now on we make the following identification: $F \mathcal{B} \cong F \psi_{0} \left( \mathcal{B} \right)$.\\
Our assumption is that the material is completely characterized by one differentiable function $W : \tilde{J}^{1} \left( F \mathcal{B} \right) \rightarrow V$ which does not depend on the end point. In this case, this condition can be translated by the following equality
\begin{equation}\label{94}
W \left( j_{x,y}^{1} \Psi\right) = W \left( j_{x,z}^{1} \left( F \tau_{z-y} \circ \Psi\right)\right), \ \forall j_{x,y}^{1}\Psi \in \tilde{J}^{1} \left( F \mathcal{B}\right).
\end{equation}
This function measures, for instance, the stored energy per unit mass and, again, we will call this function \textit{response functional} or \textit{mechanical response}.\\
Notice that, using that $W$ does not depend on the final point, we can define $W$ over $\tilde{J}^{1} \left( F\mathcal{B} , F\mathbb{R}^{3}\right)$, which is the open subset of $\tilde{J}^{1} \left( F\mathbb{R}^{3} \right)$ given by $\left(\overline{\alpha} , \overline{\beta} \right)^{-1} \left( \mathcal{B} \times \mathbb{R}^{3} \right) $.\\

Now, suppose that an infinitesimal neighbourhood of the material around the point $Y$ can be turned into a neighbourhood of $X$ such that the transformation cannot be detected by any mechanical experiment. If this condition is satisfied with every point $X$ of $ F\mathcal{B}$, the body is said \textit{uniform}. We can express this physical property in a geometric way as follows.

\begin{definition}
\rm
A Cosserat continuum $\mathcal{B}$ is said to be \textit{uniform} if for each two points $x,y \in \mathcal{B}$ there exists a local principal bundle isomorphism over the identity map on $Gl\left(3 , \mathbb{R}\right)$, $\Psi$, from an open neighbourhood $FU \subseteq F\mathcal{B}$ over $x$ to an open neighbourhood $FV \subseteq F\mathcal{B}$ over $y$ such that $\psi \left(x\right) =y$ and
\begin{equation}\label{66}
W \left( j^{1}_{y, \kappa \left(y\right)} \tilde{\kappa} \circ j^{1}_{x,y} \Psi \right) = W \left( j^{1}_{y, \kappa \left(y\right)} \tilde{\kappa}\right),
\end{equation}
for all infinitesimal deformation $j^{1}_{y , \kappa \left(y\right)} \tilde{\kappa}$.
\end{definition}
This kind of maps are relevant for the sequel and we will endow these maps with a groupoid structure. For each two points we will denote by $\overline{G} \left(x,y\right)$ the collection of all $1-$jets $j_{x,y}^{1}\Psi$ which satisfy Eq. (\ref{66}). So, the set $\overline{\Omega} \left( \mathcal{B}\right) = \cup_{x,y \in \mathcal{B}} \overline{G}\left(x,y\right)$ can be considered as a groupoid over $\mathcal{B}$ which is, indeed, a subgroupoid of the second-order non-holonomic groupoid $\tilde{J}^{1}\left(F \mathcal{B}\right)$. We will denote $\overline{\alpha}^{-1} \left(x\right)$ (resp. $\overline{\beta}^{-1} \left(x\right)$) by $\overline{\Omega}_{x}\left( \mathcal{B}\right)$ (resp. $\overline{\Omega}^{x}\left( \mathcal{B}\right)$).\\

\begin{definition}
\rm
Given a material point $x \in \mathcal{B}$ a \textit{material simmetry} at $x$ is a class of $1-$jets $j_{x, x}^{1}\Psi$, where $\Psi$ is a local automorphism at $x$ over the identity map on $Gl\left( 3 , \mathbb{R}\right)$, which satisfies Eq. (\ref{66}).
\end{definition}
We denote by $\overline{G}\left(x\right)$ the set of all material simmetries which is, indeed, the isotropy group of $\overline{\Omega} \left(\mathcal{B}\right)$ at $x$. So, the following result is obvious.
\begin{proposition}
Let $\mathcal{B}$ be a Cosserat contiuum. $\mathcal{B}$ is uniform if, and only if, $\overline{\Omega} \left( \mathcal{B}\right)$ is a reduced subgroupoid of $\tilde{J}^{1} \left( F \mathcal{B}\right)$.
\end{proposition}
Notice that, at general, we cannot ensure that $\overline{\Omega} \left( \mathcal{B}\right) \subseteq \tilde{J}^{1} \left( F\mathcal{B}  \right)$ is a Lie subgroupoid. Our assumption is that $\overline{\Omega} \left( \mathcal{B}\right)$ is in fact a Lie subgroupoid and, in this case, $\overline{\Omega} \left(\mathcal{B}\right)$ is said to be the \textit{second-order non-holonomic material groupoid of $\mathcal{B}$}.\\

As we have seen, a Cosserat medium is uniform if the function $W$ does not depend on the point $x$. In addition, a body is said to be \textit{homogeneous} if we can choose a global section of the second-order non-holonomic material groupoid which is constant on the body, more precisely:

\begin{definition}\label{151}
\rm
A Cosserat medium $\mathcal{B}$ is said to be \textit{homogeneous} if it admits a global deformation $\tilde{\kappa}$ which induces a global section of $\left(\overline{\alpha} , \overline{\beta}\right)$ in $\overline{\Omega} \left( \mathcal{B}\right)$, $\overline{\mathcal{P}}$, i.e., for each $x , y \in \mathcal{B}$
$$ \overline{\mathcal{P}}\left(x, y\right) = j^{1}_{x,y} \left(\tilde{\kappa}^{-1} \circ F\tau_{\kappa\left(y\right) - \kappa \left(x\right)} \circ \tilde{\kappa}\right),$$
where $\tau_{\kappa\left(y\right) - \kappa \left(x\right)}: \mathbb{R}^{3} \rightarrow \mathbb{R}^{3}$ denotes the translation on $\mathbb{R}^{3}$ by the vector $\kappa\left(y\right) - \kappa \left(x\right)$. $\mathcal{B}$ is said to be \textit{locally homogeneous} if there exists a covering of $\mathcal{B}$ by homogeneous open sets.
\end{definition}


Now, suppose that $\mathcal{B}$ is homogeneous. Then, if we take the global coordinates $\left(x^{i}\right) $ on $\mathcal{B}$ given by the induced diffeomorphism $\kappa$, we deduce that $\overline{\mathcal{P}}$ is expressed by
\begin{equation}\label{67}
\overline{\mathcal{P}}\left(x^{i},y^{j}\right) = \left(\left(x^{i},y^{j}, P^{j}_{i}\right) , \delta^{j}_{i} , \dfrac{\partial P^{j}_{i} }{\partial x^{k}} + \dfrac{\partial P^{j}_{i} }{\partial y^{k}} \right),
\end{equation}
If $\mathcal{B}$ is locally homogeneous we can cover $\mathcal{B}$ by local coordinate systems $\left(x^{i}\right)$ which generate (local) sections of $\left( \overline{\alpha} , \overline{\beta} \right)$ in $\overline{\Omega} \left( \mathcal{B} \right)$ satisfying Eq. (\ref{67}).\\

Next, we want to give some equivalent definitions and relate this definition with the usual one, which is given for second-order non-holonomic $\overline{G}-$structures (see \cite{MAREMDL}).

\section{Integrability}
Now, we will introduce the notion of \textit{integrability} of reduced Lie subgroupoids of the second-order non-holonomic groupoid.\\
In order to do that, we will proceed in a similar way to $\overline{F}^{2}M$. Thus, there exists a canonical Lie groupoid isomorphism over the identity on $\mathbb{R}^{n}$, $L: \tilde{J}^{1}\left(F \mathbb{R}^{n}\right) \cong \mathbb{R}^{n} \times \mathbb{R}^{n} \times \overline{G}^{2}\left(n\right)$, where $\mathbb{R}^{n} \times \mathbb{R}^{n} \times \overline{G}^{2}\left(n\right)$ is the trivial Lie groupoid of $\overline{G}^{2}\left(n\right)$ over $\mathbb{R}^{n}$ defined by
$$L\left( j^{1}_{x,y} \Psi\right) = \left(x,y, j^{1}_{0,0} \left( F \tau_{-y} \circ \Psi \circ F \tau_{x} \right) \right), \ \forall x,y \in \mathbb{R}^{n},$$
where $\tau_{-y}$ and $\tau_{x}$ denote the translations on $\mathbb{R}^{n}$ by the vectors $-y$ and $x$ respectively. Thus, if $\overline{G}$ is a Lie subgroup of $\overline{G}^{2} \left(n\right)$, we can transport $\mathbb{R}^{n} \times \mathbb{R}^{n} \times \overline{G}$ by this isomorphism to obtain a reduced Lie subgroupoid of $\tilde{J}^{1}\left(F \mathbb{R}^{n}\right)$. This reduced Lie subgroupoid of $\tilde{J}^{1}\left(F \mathbb{R}^{n}\right)$ will be called \textit{standard flat} subgroupoid of $\tilde{J}^{1}\left(F \mathbb{R}^{n}\right)$ over $\overline{G}$.\\

Let $U, V \subseteq M$ be two open subsets of $M$. We denote by $\tilde{J}^{1}\left(FU,FV\right)$ the open subset of $\tilde{J}^{1}\left(FM\right)$ defined by $\left(\overline{\alpha} , \overline{\beta} \right)^{-1}\left(U \times V\right)$. Note that if $U=V$, $\tilde{J}^{1}\left(FU, FU \right)$ is in fact the second-order non-holonomic groupoid over $U$, i.e., $\tilde{J}^{1}\left(FU, FU \right) = \tilde{J}^{1}\left(FU\right)$. Furthermore, we will think about $\tilde{J}^{1}\left(FU,FV\right)$ as the restriction of the Lie groupoid $\tilde{J}^{1}\left( FM\right)$ equipped with the restriction of the structure maps (this could not be a Lie groupoid). We will also use this notation for subgroupoids of $\tilde{J}^{1} \left(FM\right)$.\\

Next, we will introduce the notion of \textit{integrability of a reduced Lie subgroupoid}.\\

\begin{definition}
\rm
Let $\tilde{J}^{1}_{\overline{G}}\left(FM\right)$ be a reduced Lie subgroupoid of $\tilde{J}^{1}\left(F M\right)$. $\tilde{J}^{1}_{\overline{G}}\left(FM\right)$ is \textit{integrable} if it is locally isomorphic to the trivial Lie groupoid $\mathbb{R}^{n} \times \mathbb{R}^{n} \times \overline{G}$ for some Lie subgroup $\overline{G}$ of $\overline{G}^{2}\left(n\right)$.
\end{definition}

Note that $\tilde{J}^{1}_{\overline{G}}\left(FM\right)$ is "locally isomorphic" to $\mathbb{R}^{n} \times \mathbb{R}^{n} \times \overline{G} \rightrightarrows \mathbb{R}^{n}$ if for all $x\in M$ there exist an open set $U \subseteq M$ with $x \in U$ and a local chart, $\psi_{U}: U \rightarrow \overline{U}$, which induces a Lie groupoid isomorphism,
\begin{equation}\label{68}
 \Psi_{U} : \tilde{J}^{1}_{\overline{G}}\left(FU\right) \rightarrow \overline{U} \times \overline{U} \times \overline{G},
\end{equation}
such that $ \Psi_{U} = \left( \psi_{U} \circ \overline{\alpha}, \psi_{U} \circ \overline{\beta} , \overline{\Psi}_{U}\right)$, where for each $j_{x,y}^{1} \Psi \in \tilde{J}^{1}_{\overline{G}}\left(FU\right)$
$$ \overline{\Psi}_{U} \left( j_{x,y}^{1} \Psi \right) = j_{0,0}^{1}\left( F\left(\tau_{-\psi_{U}\left(y\right)} \circ \psi_{U}\right) \circ \Psi \circ F\left(\psi_{U}^{-1} \circ \tau_{\psi_{U}\left(x\right)}\right)\right).$$
So, we can claim that $\tilde{J}^{1}_{\overline{G}} \left(FM\right)$ is locally isomorphic to $\mathbb{R}^{n} \times \mathbb{R}^{n} \times \overline{G}$ if we can cover $M$ by local charts $\left(\psi_{U} , U\right)$ such that induce Lie groupoid isomorphisms from $\tilde{J}^{1}_{\overline{G}} \left(FU\right)$ to the restrictions of the standard flat over $\overline{G}$ to $\overline{U}$.\\

\begin{remarkth}\label{69}
\rm
Let $\tilde{J}^{1}_{\overline{G}} \left(FM\right)$ be an integrable subgroupoid of $\tilde{J}^{1} \left(FM\right)$, i.e., locally isomorphic to $\mathbb{R}^{n} \times \mathbb{R}^{n} \times \overline{G}$ with $\overline{G}$ a Lie subgroup of $\overline{G}^{1} \left( n \right)$. Suppose that there exists another Lie subgroup of $\overline{G}^{2}\left(n\right)$, $\tilde{G}$, such that $\tilde{J}^{1}_{\overline{G}} \left(FM\right)$ is locally isomorphic to $\mathbb{R}^{n} \times \mathbb{R}^{n} \times \tilde{G}$. Then, it is easy to see that $G$ and $\tilde{G}$ are conjugated subgroups of $\overline{G}^{2}\left(n\right)$, i.e., there exists a frame at $0$, $\overline{g} \in \overline{G}^{2}\left(n\right)$, such that
$$  \tilde{G} = \overline{g}^{-1} \cdot \overline{G} \cdot \overline{g}.$$
However, in this case the converse is not true.

\end{remarkth}

Now, there is a special kind of reduced subgroupoids of $\tilde{J}^{1} \left(FM\right)$ which will play an important role in the following. A trivial reduced subgroupoid of $\tilde{J}^{1} \left(FM\right)$ or \textit{parallelism of} $\tilde{J}^{1} \left(FM\right)$ is a reduced subgroupoid of $\tilde{J}^{1} \left(FM\right)$, $\tilde{J}^{1}_{\overline{e}} \left(FM\right)$, such that for each $x,y \in M$ there exists a unique $1-$jet $j^{1}_{x,y}\Psi \in \tilde{J}^{1}_{\overline{e}} \left(FM\right)$, or equivalently, the restriction of $\left(\overline{\alpha} , \overline{\beta}\right)$ to $\tilde{J}^{1}_{\overline{e}} \left(FM\right)$ is a Lie groupoid isomorphism.\\
So, having a trivial reduced subgroupoid $\tilde{J}^{1}_{\overline{e}} \left( FM\right)$ of $\tilde{J}^{1} \left(FM\right)$ we can consider a global section of $\left(\overline{\alpha} , \overline{\beta}\right)$, $\overline{\mathcal{P}}: M \times M \rightarrow \tilde{J}^{1} \left(FM\right)$, such that $\overline{\mathcal{P}}\left(x,y\right)$ is the unique $1-$jet from $x$ to $y$ which is in $\tilde{J}^{1}_{\overline{e}} \left(FM\right)$, i.e., $\overline{\mathcal{P}} = \left(\overline{\alpha} , \overline{\beta}\right)_{|\tilde{J}^{1}_{\overline{e}}\left(FM\right)}^{-1}$. Conversely, every global section of $\left(\overline{\alpha} , \overline{\beta}\right)$ can be seen as a parallelism of $\tilde{J}^{1} \left(FM\right)$ (we are understanding "section" as section in the category of Lie groupoids, i.e., Lie groupoid morphism from the pair groupoid $M \times M$ to $\tilde{J}^{1} \left( FM\right)$ which is a section of the morphism $\left(\overline{\alpha} , \overline{\beta}\right)$). Using this, we can also speak about \textit{integrable sections of} $\left(\overline{\alpha} , \overline{\beta}\right)$.\\

Next, we will express a necessary result to interpret the integrability in another equivalent way. The proof is very technical and is omitted.\\

\begin{proposition}\label{70}

Let $\tilde{J}^{1}_{\overline{G}} \left(FM\right)$ be a reduced Lie subgroupoid of $\tilde{J}^{1} \left(FM\right)$. $\tilde{J}^{1}_{\overline{G}} \left(FM\right)$ is integrable if, and only if, for all $x,y \in M$ there exist two open sets $U , V \subseteq M$ with $x \in U$ and $y \in V$ and two local charts $\psi_{U} : U \rightarrow \overline{U} $ and $\psi_{V} : V  \rightarrow \overline{V}$ which induce a diffeomorphism
\begin{equation}\label{71}
\begin{array}{rccl}
\Psi_{U,V} : & \tilde{J}^{1}_{\overline{G}}\left(FU,FV\right) & \rightarrow & \overline{U} \times \overline{V} \times \overline{G} \\
& j^{1}_{x,y}\Psi  &\mapsto &  \left(\psi_{U}\left(x\right), \psi_{V}  \left(y\right), \overline{\Psi}_{U,V}\left(j^{1}_{x,y} \Psi \right)\right).
\end{array}
\end{equation}
where,
$$ \overline{\Psi}_{U,V}\left(j^{1}_{x,y} \Psi \right) = j_{0,0}^{1} \left(F\left( \tau_{-\psi_{V} \left(y\right)} \circ \psi_{V} \right)\circ \Psi \circ F\left(\psi_{U}^{-1} \circ \tau_{\psi_{U}\left(x\right)}\right)\right).$$\\

\end{proposition}

Let $\overline{\mathcal{P}} : M \times M \rightarrow \tilde{J}^{1} \left(FM\right)$ be a section of $\left(\overline{\alpha} , \overline{\beta}\right)$. Using this result we can claim that $\overline{\mathcal{P}}$ is integrable if, and only if, for each $x,y \in M$
\begin{equation}\label{73}
\overline{\mathcal{P}} \left(x,y\right) = j^{1}_{x,y} \left(F \left(\psi_{V}^{-1} \circ \tau_{\psi_{V} \left(y\right) - \psi_{U} \left(x\right)} \circ \psi_{U} \right)\right),
\end{equation}
for some two local charts $\left( \psi_{U} , U\right) , \left(\psi_{V} , V\right)$ on $M$ through $x$ and $y$ respectively.\\
Equivalently, using the local coordinates given in Eq. (\ref{45}), $\overline{\mathcal{P}}$ can be written locally as follows,
\begin{equation}\label{74}
\overline{\mathcal{P}}\left(x^{i}, y^{j} \right) = \left(\left(x^{i} , y^{j} , \delta^{j}_{i}\right),\delta^{j}_{i} , 0 \right),
\end{equation}

Let $\tilde{J}^{1}_{\overline{G}} \left(FM\right)$ be a reduced subgroupoid of $\tilde{J}^{1} \left(FM\right)$ and $\overline{Z}_{0} \in \overline{F}^{2} M$ be a second-order non-holonomic frame at $z_{0} \in M$. Then, we define
\begin{equation}\label{75}
\overline{G}  := \{ \overline{Z}_{0}^{-1} \cdot \overline{g} \cdot \overline{Z}_{0} \ / \ \overline{g} \in \tilde{J}^{1}_{\overline{G}} \left(z_{0}\right) \} = \overline{Z}_{0}^{-1} \cdot  \tilde{J}^{1}_{\overline{G}} \left(z_{0}\right) \cdot \overline{Z}_{0},
\end{equation}
where $\tilde{J}^{1}_{\overline{G}} \left(z_{0}\right)$ is the isotropy group of $\tilde{J}^{1}_{\overline{G}} \left(FM\right)$ at $z_{0}$. Therefore, $\overline{G}$ is a Lie subgroup of $\overline{G}^{2}\left(n\right)$. This Lie group will be called \textit{associated Lie group} to $\tilde{J}^{1}_{\overline{G}} \left(FM\right)$.\\
Note that, contrary to the case of non-holonomic $\overline{G}-$structures of second order, we do not have a unique Lie group $\overline{G}$. In fact, let $\tilde{Z}_{0}$ be a non-holonomic frame of second order at $\tilde{z}_{0}$ and $\tilde{G}$ be the associated Lie group; then, if we take $\overline{L}_{z_{0}, \tilde{z}_{0}} \in \tilde{J}^{1}_{\overline{G}}\left(z_{0} , \tilde{z}_{0}\right)$ we have
$$ \overline{G} = [\tilde{Z}_{0}^{-1}\cdot \overline{L}_{z_{0}, \tilde{z}_{0}} \cdot \overline{Z}_{0}]^{-1} \cdot \tilde{G} \cdot[\tilde{Z}_{0}^{-1}\cdot \overline{L}_{z_{0}, \tilde{z}_{0}} \cdot \overline{Z}_{0}],$$
i.e., $\overline{G}$ and $\tilde{G}$ are conjugated subgroups of $\overline{G}^{2}\left(n\right)$. Notice that, by construction, if $\tilde{J}^{1}_{\overline{G}}\left( FM \right)$ is integrable by $\overline{G}$ (i.e. locally isomorphic to the Lie trivial Lie groupoid of $\overline{G}$ over $\mathbb{R}^{n}$), $\overline{G}$ can be constructed using Eq. (\ref{75}).

\begin{proposition}\label{104}
A reduced subgroupoid $\tilde{J}^{1}_{\overline{G}}\left(FM\right)$ of $\tilde{J}^{1}\left(FM\right)$ is integrable if, and only if, for each point $x  \in M$ there exists a (local) coordinate system $\left(x^{i}\right)$ on an open set $U \subseteq M$ with $x \in U$ such that the local section, 
\begin{equation}\label{77}
\overline{\mathcal{P}}\left(x^{i}, y^{j}\right) = \left(\left(x^{i}, y^{j} ,\delta^{j}_{i}\right),\delta^{j}_{i} , 0 \right),
\end{equation}
takes values into $\tilde{J}^{1}_{\overline{G}}\left(FM\right)$.

\begin{proof}
The proof is similar to the proof of Proposition \ref{31}. Let us consider an (local) integrable section $\overline{\mathcal{P}} : U \times U \rightarrow \tilde{J}^{1}_{\overline{G}}\left(FU\right)$ of $\left(\overline{\alpha} , \overline{\beta}\right)$ given by
$$ \overline{\mathcal{P}}\left(x,y\right) = j_{x,y}^{1}\left(F\left(\psi_{U}^{-1} \circ \tau_{\psi_{U}\left(y\right) - \psi_{U}\left(x\right)} \circ \psi_{U}\right)\right),$$
where $\psi_{U} : U \rightarrow \overline{U}$ is a local chart. Now, Let $z_{0} \in U$ be a point at $U$, we choose $\overline{Z}_{0} \triangleq j_{0,z_{0}}^{1}\left(F\left(\psi_{U}^{-1} \circ \tau_{\psi_{U}\left(z_{0}\right)}\right)\right) \in \overline{F}^{2} U$ be a non-holonomic second-order frame at $z_{0}$ and $\overline{G}$ be the Lie subgroup of $\overline{G}^{2}\left(n\right)$ satisfying Eq. (\ref{75}). Then, the rest of the proof is easy.
\end{proof}

\end{proposition}

Now, we want to define the notion of \textit{second-order non-holonomic prolongation} in $\tilde{J}^{1}\left(FM\right)$. In order to do this, we will define the projections $\overline{\Pi}_{1}^{2}$ and $\tilde{\Pi}_{1}^{2}$ which will be closely related with the maps $\tilde{\pi}^{2}_{1}$ and $\overline{\pi}^{2}_{1}$ (see section 2) by the Equalities given in Remark \ref{120}. Thus, we define

$$
\begin{array}{rccl}
\overline{\Pi}_{1}^{2} : & \tilde{J}^{1}\left(FM\right)  & \rightarrow & \Pi^{1}\left(M,M\right)\\
& j_{x,y}^{1}\Psi  &\mapsto &  \Psi\left(X\right)[X^{-1}]
\end{array}
$$
where $X \in FM$ is a frame at $x$. It is easy to show that $\overline{\Pi}_{1}^{2}$ is well-defined and, locally,
$$\overline{\Pi}_{1}^{2} \left(\left(x^{i}\right),\left( y^{j} , y^{j}_{i}\right),y^{j}_{,i},y^{j}_{i,k}\right) = \left(x^{i}, y^{j} , y^{j}_{i}\right).$$
On the other hand we consider
$$
\begin{array}{rccl}
\tilde{\Pi}_{1}^{2} : & \tilde{J}^{1}\left(FM\right)  & \rightarrow & \Pi^{1}\left(M,M\right)\\
& j_{x,y}^{1}\Psi  &\mapsto &  j_{x,y}^{1}\psi
\end{array}
$$
where $\psi$ is the induced map of $\Psi$ over $M$. Then, locally
$$\tilde{\Pi}_{1}^{2} \left(\left(x^{i}\right),\left( y^{j} , y^{j}_{i}\right),y^{j}_{,i},y^{j}_{i,k}\right) = \left(x^{i}, y^{j} , y^{j}_{,i}\right).$$
Notice that $\overline{\Pi}_{1}^{2}$ and $\tilde{\Pi}_{1}^{2}$ are, indeed, Lie groupoid morphims over the identity map on $M$. Then, let $\overline{\mathcal{P}} : M \times M \rightarrow \tilde{J}^{1}\left(FM\right)$ be a section of $\left( \overline{\alpha}, \overline{\beta}\right)$ in $\tilde{J}^{1}\left(FM\right)$ the projections $\mathcal{P} = \overline{\Pi}_{1}^{2} \circ \overline{\mathcal{P}}$ and $\mathcal{Q} = \tilde{\Pi}_{1}^{2} \circ \overline{\mathcal{P}}$ are sections of $\left(\alpha , \beta \right)$ in $\Pi^{1}\left(M,M\right)$.\\

Next, we will invert this process and, to do this, we will get inspired from Remark \ref{36}. Let $\mathcal{P}, \mathcal{Q} : M \times M \rightarrow \Pi^{1}\left(M,M\right)$ be two sections of $\left(\alpha , \beta\right)$ in $\Pi^{1}\left(M,M\right)$ such that
$$ \mathcal{Q}\left(x,y\right) = j_{x,y}^{1}\psi_{xy}, \ \forall x , y \in M.$$
Thus, we constuct the following map

$$
\begin{array}{rccl}
\overline{\mathcal{P}  \circ \psi_{xy}} : & FU & \rightarrow & FV \\
& j^{1}_{0,f\left(0\right)}f  &\mapsto &  \mathcal{P}\left(f\left(0\right), \psi_{xy} \left(f\left(0\right)\right)\right) \cdot j^{1}_{0,f\left(0\right)}f
\end{array}
$$
where $\psi_{xy}: U \rightarrow V$. Analogously to Remark \ref{36}, $\overline{ \mathcal{P} \circ \psi_{xy}}$ is a local principal bundle isomorphism with $\psi_{xy}: U \rightarrow V$ as its induced map over $M$. In fact, the inverse is given by
$$ j^{1}_{0,g\left(0\right)}g  \mapsto   [\mathcal{P}\left(\psi_{xy}^{-1}\left(g\left(0\right)\right),g\left(0\right)\right)]^{-1} \cdot j^{1}_{0,g\left(0\right)}g.$$
Furthermore, let $\left(x^{i}\right)$ be a local coordinate system on an open set $U \subseteq M$ and $\left(x^{i}, x^{i}_{j}\right)$ its induced coordinates over $FU$, we have
$$ \overline{\mathcal{P} \circ \psi_{xy}} \left(x^{i}, x^{i}_{j} \right) = \left(\psi_{xy}\left(x^{i}\right), P^{j}_{l}\left(x^{i}, \psi_{xy}\left(x^{i}\right)\right) x^{l}_{i} \right),$$
where for each another local coordinate system $\left(y^{j}\right)$ on an open set $V \subseteq M$
$$ \mathcal{P}\left(x^{i}, y^{j} \right) = \left(x^{i}, y^{j}, P^{j}_{i}\left(x^{i}, y^{j}\right)\right).$$
Thus, we define

\begin{equation}\label{78}
\begin{array}{rccl}
\mathcal{P}^{1}\left(\mathcal{Q}\right)  : & M \times M & \rightarrow & \tilde{J}^{1}\left(FM\right) \\
& \left(x,y\right)  &\mapsto &  j_{x,y}^{1} \left( \overline{\mathcal{P} \circ \psi_{xy}}\right)
\end{array}
\end{equation}
where we are considering the equivalence class in $\tilde{J}^{1}\left(FM\right)$. Notice that $\mathcal{P}^{1}\left(\mathcal{Q}\right)$ does not depend on $\psi_{xy}$ because of $\mathcal{Q}$ does not depend on $\psi_{xy}$. $\mathcal{P}^{1}\left(\mathcal{Q}\right)$ will be called \textit{second-order non-holonomic prolongation of $\mathcal{P}$ and $\mathcal{Q}$} and satisfies that
\begin{itemize}
\item[(i)] For all $x,y \in M$ and $j_{0,x}^{1}f \in F M$,
$$ \overline{\Pi}^{2}_{1} \circ \mathcal{P}^{1}\left(\mathcal{Q}\right) \left(x,y\right) = [ \mathcal{P}\left(x, \psi_{xy} \left(x\right)\right) \cdot j_{0,x}^{1} f ] \cdot \left(j_{0,x}^{1}f\right)^{-1} = \mathcal{P}\left(x,y\right).$$

\item[(ii)]  For all $x,y \in M$,
$$ \tilde{\Pi}^{2}_{1} \circ \mathcal{P}^{1}\left(\mathcal{Q}\right) \left(x,y\right)  = j_{x,y}^{1}\psi_{xy} = \mathcal{Q}\left(x,y\right).$$
\end{itemize}
In fact, let be $\left(x^{i}\right)$ and $\left(y^{j}\right)$ local coordinate systems on open sets $U,V \subseteq M$ and $\left(x^{i}, x^{i}_{j}\right)$ and $\left(y^{j}, y^{j}_{i}\right)$ its induced coordinates over $FU$ and $FV$ respectively, then we have
$$\mathcal{P}^{1}\left(\mathcal{Q}\right) \left(x^{i}, y^{j} \right) =  \left(\left(x^{i} , y^{j} , P^{j}_{i}\right), Q^{j}_{i}, R^{j}_{i,k} \right),$$
where
$$ \mathcal{P}\left(x^{i}, y^{j} \right) = \left(x^{i}, y^{j}, P^{j}_{i}\right), \ \ \ \ \ \mathcal{Q}\left(x^{i}, y^{j} \right) = \left(x^{i}, y^{j}, Q^{j}_{i}\right).$$
Furthermore, for each $k=1, \cdots, n$
$$\dfrac{\partial \left(P^{j}_{i} \circ \left( Id_{U}, \psi_{xy}\right)\right)}{\partial x^{k}_{|x} }= d{ P^{j}_{i}}_{|\left(x,y\right)}\circ \dfrac{\partial \left( Id_{U}, \psi_{xy}\right)}{\partial x^{k}_{|x} } $$
$$= d{ P^{j}_{i}}_{|\left(x,y\right)}\left( x^{k} ,  \left(Q_{k}^{1}\left(x,y\right), \cdots ,Q_{k}^{n}\left(x,y\right)\right)\right) $$
$$ =d{\left(P^{j}_{i}\right)_{y}}_{|x} \left(x^{k}\right) + d{\left(P^{j}_{i}\right)_{x}}_{|y} \left(Q_{k}^{1}\left(x,y\right), \cdots ,Q_{k}^{n}\left(x,y\right)\right)$$
$$ =\dfrac{\partial \left(P^{j}_{i}\right)_{y}}{\partial x^{k}_{|x} } + d{\left(P^{j}_{i}\right)_{x}}_{|y} \left(Q_{k}^{1}\left(x,y\right), \cdots ,Q_{k}^{n}\left(x,y\right)\right),$$
where we are fixing the first (resp. the second) coordinate when we write $\left(P^{j}_{i}\right)_{x}$ (resp. $\left(P^{j}_{i}\right)_{y}$). Then, by definition of induced coordinates, $R^{j}_{i,k}$ is given by
$$R^{j}_{i,k}\left(x,y\right) = \dfrac{\partial \left(P^{j}_{i}\right)_{y}}{\partial x^{k}_{|x} } + \sum_{l} Q_{k}^{l}\left(x,y\right) \dfrac{\partial \left(P^{j}_{i}\right)_{x}}{\partial y^{l}_{|y} }.$$
We will denote this expression by
$$R^{j}_{i,k} = \dfrac{\partial P^{j}_{i}}{\partial x^{k} } + \sum_{l} Q_{k}^{l}\dfrac{\partial P^{j}_{i}}{\partial y^{l} }.$$
Then, any section $\overline{\mathcal{P}}  \left(x^{i}, y^{j} \right) =  \left(\left(x^{i} , y^{j} , P^{j}_{i}\right), Q^{j}_{i}, R^{j}_{i,k} \right) $ of $\left(\overline{\alpha}, \overline{\beta}\right)$ in \linebreak $\tilde{J}^{1}\left(FM\right)$ which projects onto $\mathcal{P}  \left(x^{i}, y^{j} \right) =  \left(x^{i} , y^{j} , P^{j}_{i}\right)$ and $\mathcal{Q}   \left(x^{i}, y^{j} \right) =  \left(x^{i} , y^{j} , Q^{j}_{i} \right)$ via $\overline{\Pi}^{2}_{1}$ and $\tilde{\Pi}^{2}_{1}$ respectively is a prolongation if, and only if,
\begin{equation}\label{105}
R^{j}_{i,k} = \dfrac{\partial P^{j}_{i}}{\partial x^{k} } + \sum_{l} Q_{k}^{l}\dfrac{\partial P^{j}_{i}}{\partial y^{l} }.
\end{equation}
Thus, we have established the notion of prolongation in the second-order non-holonomic groupoid. Then, we can give the following definition.

\begin{definition}
\rm
Let $\mathcal{P}^{1}\left(\mathcal{Q}\right)$ be a non-holonomic prolongation of second order in $\tilde{J}^{1}\left(FM\right)$. $\mathcal{P}^{1}\left(\mathcal{Q}\right)$ is said to be \textit{integrable} in $\tilde{J}^{1}\left(FM\right)$ if $\mathcal{Q}$ is an integrable section of $\Pi^{1}\left(M,M\right)$.
\end{definition}
Notice that, using the introduced coordinates, an integrable prolongation can be seen locally as follows
\begin{eqnarray*}\label{79}
\mathcal{P}^{1}\left(\mathcal{Q}\right) \left(x^{i}, y^{j} \right) &=&  \left(\left(x^{i} , y^{j} , P^{j}_{i}\right), \delta^{j}_{i}, \dfrac{\partial P^{j}_{i}}{\partial x^{k}} + \dfrac{\partial P^{j}_{i} }{\partial y^{k} }\right)
\end{eqnarray*}

Thus, as in the case of second-order non-holonomic frame bundle, we have two remarkable sections: integrable sections and integrable prolongations. So, it is easy to give the following result (similar to Proposition \ref{116}) which will help us to understand why the integrable prolongation are not necessarily integrable as sections.

\begin{proposition}\label{80}
Let $\overline{\mathcal{P}}$ be a section of $\left(\overline{\alpha} , \overline{\beta}\right)$ in $\tilde{J}^{1}\left(FM\right)$. $\overline{\mathcal{P}}$ is integrable if, and only if, $\overline{\mathcal{P}} = \mathcal{P}^{1}\left(\mathcal{Q}\right)$ is a second-order non-holonomic integrable prolongation and $\mathcal{P} = \mathcal{Q}$. In particular, a second-order non-holonomic integrable prolongation $\mathcal{P}^{1}\left(\mathcal{Q}\right)$ is integrable if, and only if, $\mathcal{P}^{1} \left( \mathcal{Q} \right)$ takes values in $\tilde{j}^{1} \left( FM \right)$.
\end{proposition}
Thus, second-order non-holonomic integrable prolongations can be seen as a natural generalization of integrable sections of $\left(\overline{\alpha} , \overline{\beta}\right)$ (see Proposition \ref{81}).\\

Notice that, analogously to Eq. (\ref{37}), we can prove the following result.

\begin{proposition}\label{81}
Let $\overline{\mathcal{P}}$ be a section of $\left(\overline{\alpha} , \overline{\beta}\right)$ in $\tilde{J}^{1}\left(FM\right)$. $\overline{\mathcal{P}}$ is a second-order non-holonomic integrable prolongation if, and only if, for all $x_{0},y_{0} \in M$ there exist two open sets $U,V \subseteq M$ with $x_{0} \in U $ and $y_{0} \in V$ and two local principal bundle isomorphisms $\Psi : FV \rightarrow F \overline{V}$ and $\Phi : FU \rightarrow F \overline{U}$ such that
\begin{equation}\label{82}
\overline{\mathcal{P}}\left(x, y\right) = j^{1}_{x,y} \left(\Psi^{-1} \circ F\tau_{\psi\left(y\right) - \phi \left(x\right)} \circ \Phi\right), \ \forall \left(x,y\right) \in U \times V.
\end{equation}
\end{proposition}


Now, we will extend this concept to reduced subgroupoids.
\begin{definition}\label{83}
\rm
Let $\tilde{J}^{1}_{\overline{G}} \left(FM\right)$ be a reduced subgroupoid of $\tilde{J}^{1}\left(FM\right)$. $\tilde{J}^{1}_{\overline{G}} \left(FM\right)$ is an \textit{integrable prolongation} if can be covered $M$ with local integrable prolongations which take values in $\tilde{J}^{1}_{\overline{G}} \left(FM\right)$.

\end{definition}
\begin{proposition}
Let $\tilde{J}^{1}_{\overline{G}} \left(FM\right)$ be an integrable prolongation. $\tilde{J}^{1}_{\overline{G}} \left(FM\right)$ is integrable if, and only if, $\tilde{J}^{1}_{\overline{G}} \left(FM\right)$ is contained in $\tilde{j}^{1} \left(FM\right)$.
\end{proposition}
Notice that, Definition \ref{83} can be expressed as follows: For any point $x \in M$ there exists a local coordinate system $\left(x^{i}\right)$ over an open set $U \subseteq M$ which contains $x$ such that there is a local section
\begin{equation}\label{84}
\mathcal{P}^{1}\left(\mathcal{Q}\right) \left(x^{i}, y^{j} \right) =  \left(\left(x^{i} , y^{j} , P^{j}_{i}\right), \delta^{j}_{i}, \dfrac{\partial P^{j}_{i} }{\partial x^{k}} + \dfrac{\partial P^{j}_{i} }{\partial y^{k} }\right),
\end{equation}
which takes values in $\tilde{J}^{1}_{\overline{G}} \left(FM\right)$.\\

\begin{remarkth}\label{135}
\rm
Analogously to Remark \ref{39}, let $\tilde{J}^{1}_{\overline{G}} \left(FM\right)$ be a reduced subgroupoid of $\tilde{J}^{1} \left(FM\right)$. We can prove that $\tilde{J}^{1}_{\overline{G}} \left(FM\right)$ is an integrable prolongation if, and only if, for each point $x \in M$, there exists a local isomorphism of principal bundles, $\Psi_{U}: FU \rightarrow F\overline{U}$, with $x \in U$ such that induces an isomorphism of Lie groupoids given by
$$ \Upsilon_{U} : \tilde{J}^{1}_{\overline{G}} \left(FU\right)  \rightarrow \overline{U} \times \overline{U} \times \overline{G},$$
where $\Upsilon_{U}\left( j_{x,y}^{1} H \right) = \left(\psi_{U} \left(x\right) , \psi_{U} \left(y\right) , \overline{\Upsilon}_{U} \left( j_{x,y}^{1} H \right)\right)$ and
$$ \overline{\Upsilon}_{U} \left( j_{x,y}^{1} H \right) =  j_{0,0}^{1} \left(F\left(\tau_{-\psi_{U} \left(y\right)}\right) \circ \Psi_{U} \circ H \circ \Psi_{U}^{-1} \circ F\left(\tau_{\psi_{U} \left(x\right)}\right)\right) ,$$
with $\psi_{U}$ is the induced map of $\Psi_{U}$ over the base manifold.

\end{remarkth}

So, in a similar way to Proposition \ref{70}, we may prove the following:

\begin{proposition}\label{147}

Let $\tilde{J}^{1}_{\overline{G}} \left(FM\right)$ be a reduced Lie subgroupoid of \linebreak$\tilde{J}^{1} \left(FM\right)$. $\tilde{J}^{1}_{\overline{G}} \left(FM\right)$ is an integrable prolongation if, and only if, for all $x,y \in M$ there exist two open sets $U , V \subseteq M$ with $x \in U$ and $y \in V$ and two local isomorphisms $\Psi_{U} : FU \rightarrow F\overline{U} $ and $\Psi_{V} : FV  \rightarrow F\overline{V}$ which induce the following isomorphism of Lie groupoids:
\begin{equation}\label{148}
 \Upsilon_{U,V} : \tilde{J}^{1}_{\overline{G}} \left(FU, FV\right)  \rightarrow \overline{U} \times \overline{V} \times \overline{G},
 \end{equation}
where $\Upsilon_{U,V}\left( j_{x,y}^{1} H \right) = \left(\psi_{U} \left(x\right) , \psi_{V} \left(y\right) , \overline{\Upsilon}_{U,V} \left( j_{x,y}^{1} H \right)\right)$ and
$$ \overline{\Upsilon}_{U,V} \left( j_{x,y}^{1} H \right) =  j_{0,0}^{1} \left(F\left(\tau_{-\psi_{V} \left(y\right)}\right) \circ \Psi_{V} \circ H \circ \Psi_{U}^{-1} \circ F\left(\tau_{\psi_{U} \left(x\right)}\right)\right) ,$$
with $\psi_{U}$ and $\psi_{V}$ are the induced map of $\Psi_{U}$ and $\Psi_{V}$ over the base manifold respectively.
\end{proposition}

Hence, we have that: $\tilde{J}^{1}_{\overline{G}} \left( FM \right)$ is an integrable prolongation if, and only if, for any two point $x, y\in M$ there exist two local coordinate systems $\left(x^{i}\right)$ and $\left( y^{j} \right)$ over $x$ and $y$ respectively such that there is a local section
\begin{equation}\label{84}
\mathcal{P}^{1}\left(\mathcal{Q}\right) \left(x^{i}, y^{j} \right) =  \left(\left(x^{i} , y^{j} , P^{j}_{i}\right), \delta^{j}_{i}, \dfrac{\partial P^{j}_{i} }{\partial x^{k}} + \dfrac{\partial P^{j}_{i} }{\partial y^{k} }\right),
\end{equation}
which takes values in $\tilde{J}^{1}_{\overline{G}} \left(FM\right)$.\\



Now, we will translate these results to the associated Lie algebroid. Thus, we will express the notions of integrability over the second-order non-holonomic algebroid over an $n-$manifold $M$. We will begin defining the notion of \textit{integrability of a reduced Lie subalgebroid}. In order to do that, we will denote by $\overline{\mathfrak{g}}^{2}\left(n\right)$ the associated Lie algebra of $\overline{G}^{2}\left(n\right)$.\\

\begin{definition}
\rm
Let $A\tilde{J}^{1}_{\overline{G}}\left(FM\right)$ be a reduced Lie subalgebroid of $A\tilde{J}^{1}\left( F M \right)$. $A\tilde{J}^{1}_{\overline{G}}\left( F M \right)$ is \textit{integrable by $\overline{G}$} if it is locally isomorphic to the trivial algebroid $T \mathbb{R}^{n} \oplus \left( \mathbb{R}^{n} \times \overline{\mathfrak{g}} \right)$ where $\overline{\mathfrak{g}}$ is the Lie subalgebra of $\overline{\mathfrak{g}}^{2}\left( n \right)$.
\end{definition}
We will consider $\overline{G}$ the (unique) Lie subgroup of $\overline{G}^{2}\left( n \right)$ whose associated Lie algebra is $\overline{\mathfrak{g}}$.\\
Note that $A\tilde{J}^{1}_{\overline{G}}\left(FM\right)$ is locally isomorphic to $T \mathbb{R}^{n} \oplus \left( \mathbb{R}^{n} \times \overline{\mathfrak{g}} \right)$  if for all $x\in M$ there exists an open set $U \subseteq M$ with $x \in U$ and a local chart, $\psi_{U}: U \rightarrow \overline{U}$, which induces an isomorphism of Lie algebroids,
\begin{equation}\label{86}
 A\Psi_{U} : A\tilde{J}^{1}_{\overline{G}}\left(FU\right) \rightarrow T \overline{U} \oplus \left( \overline{U} \times \overline{\mathfrak{g}} \right),
\end{equation}
where $A \Psi_{U}$ is the induced map of the isomophism of Lie groupoids $\Psi_{U}$ which is given by
$$\Psi_{ U} : \tilde{J}^{1}_{\overline{G}} \left(FU\right) \rightarrow \overline{U} \times \overline{U} \times \overline{G},$$
such that $ \Psi_{U} = \left( \psi_{U} \circ \overline{\alpha}, \psi_{U} \circ \overline{\beta} , \overline{\Psi}_{U}\right)$, where for each $j_{x,y}^{1} \Psi \in \tilde{J}^{1}_{\overline{G}}\left(FU\right)$
\begin{equation}\label{87}
\overline{\Psi}_{U} \left( j_{x,y}^{1} \Psi \right) =j_{0,0}^{1}\left( F\left(\tau_{-\psi_{U}\left(y \right)} \circ \psi_{U} \right)\circ \Psi \circ F\left(\psi_{U}^{-1} \circ \tau_{\psi_{U}\left(x\right)}\right)\right),
\end{equation}
for some Lie subgroupoid $\tilde{J}^{1}_{\overline{G}} \left(FU\right)$ of $\tilde{J}^{1}\left(FU\right)$.\\

So, for each open $U \subseteq M$, $A \tilde{J}^{1}_{\overline{G}} \left(FU\right)$ is integrable by a Lie subgroupoid $\tilde{J}^{1}_{\overline{G}} \left(FU\right)$ of $\tilde{J}^{1} \left(FU\right)$. Using the uniqueness of integrating immersed (source-connected) subgroupoids (see for example \cite{IMJMLS}), \linebreak $A\tilde{J}^{1}_{\overline{G}} \left(FM\right)$ is integrable by a Lie subgroupoid of $\tilde{J}^{1} \left(FM\right)$ which will be denoted by $\tilde{J}^{1}_{\overline{G}} \left(FM\right)$. Obviously, $A\tilde{J}^{1}_{\overline{G}} \left(FM\right)$ is integrable if, and only if, $\tilde{J}^{1}_{\overline{G}} \left(FM\right)$ is integrable.\\\\

Analogously to the case of the $1-$jets groupoid, a \textit{parallelism} of $A\tilde{J}^{1}\left( FM \right)$ is an associated Lie algebroid of a parallelism of $\tilde{J}^{1}\left( FM\right)$. Hence, using the Lie's second fundamental theorem, a parallelism is a section of $\overline{\sharp}$, where $\overline{\sharp}$ is the anchor of $A\tilde{J}^{1}\left(FM\right)$ (understanding "section" as section in the category of Lie algebroids, i.e., Lie algebroid morphism from the tangent algebroid $TM$ to $A\tilde{J}^{1}\left(FM\right)$ which is a section of the morphism $\overline{\sharp}$), and conversely. In this way, we will also speak about \textit{integrable sections of} $\overline{\sharp}$.\\

Let $\left(x^{i}\right)$ be a local coordinate system defined on some open subset $U \subseteq M$, then, we will use the local coordinate system defined in Eq. (\ref{65}),
\begin{equation}\label{88}
A\tilde{J}^{1}\left(FU\right) :  \left(x^{i},v^{i}, v^{i}_{j}, v^{i}_{,j}, v^{i}_{j,k}\right)
\end{equation}
which are, indeed, induced coordinates by the functor $A$ of local coordinates on $\tilde{J}^{1}\left(FU\right)$.\\

Notice that each integrable section of $\left(\overline{\alpha} , \overline{\beta} \right)$ in $\tilde{J}^{1}\left(FM\right)$, $\overline{\mathcal{P}}$, is a Lie groupoid morphism. Hence, $\overline{\mathcal{P}}$ induces a Lie algebroid morphism $A\overline{\mathcal{P}} : TM \rightarrow A\tilde{J}^{1}\left(FM\right)$ which is a section of $\overline{\sharp}$ and is given (see Eq. (\ref{60})) by
\begin{equation}\label{112}
 A\overline{\mathcal{P}}\left(v_{x}\right) = T_{x}\overline{\mathcal{P}}_{x} \left( v_{x}\right), \ \forall v_{x} \in T_{x} M,
\end{equation}
where $\overline{\mathcal{P}}_{x} : M \rightarrow {\tilde{J}^{1}}{x}\left( FM\right)$ satisfies that
$$ \overline{\mathcal{P}}_{x}\left(y\right)=\overline{\mathcal{P}}\left(y,x\right), \forall x,y \in M.$$
So, taking into account Eq. (\ref{74}), locally, 
$$ \overline{\mathcal{P}} \left(x^{i} , y^{j} \right) = \left(\left(x^{i} , y^{j} , \delta^{j}_{i}\right),\delta^{j}_{i} , 0 \right),$$
we have that each integrable section can be written locally as follows
$$ A \overline{\mathcal{P}} \left(x^{i}, \dfrac{\partial}{\partial x^{i}} \right) = \left( \left(x^{i} , \dfrac{\partial}{\partial x^{i}} , 0\right), 0 , 0 \right).$$

Now, using Proposition \ref{104}, we have the following analogous proposition.
\begin{proposition}\label{89}
A reduced subalgebroid $A\tilde{J}^{1}_{\overline{G}}\left( FM\right)$ of $A\tilde{J}^{1}\left( FM\right)$ is integrable by $\overline{G}$ if, and only if, there exist local integrable sections of $\overline{\sharp}$ covering $M$ which takes values on $A\tilde{J}^{1}_{\overline{G}}\left( FM\right)$.
\end{proposition}
Equivalently, for each point $ x \in M$ there exists a local coordinate system $\left(x^{i}\right)$ over an open set $U \subseteq M$ with $x \in U$ such that the local sections
$$\Lambda \left(x^{i}, \dfrac{\partial}{\partial x^{i}} \right) = \left( \left(x^{i} , \dfrac{\partial}{\partial x^{i}} , 0\right), 0 , 0 \right),$$
takes values in $A \tilde{J}^{1}_{\overline{G}}\left(FM\right)$.\\

Next, we will have to introduce the notion of prolongation over the induced Lie algebroid $A \tilde{J}^{1} \left( FM \right)$. In this way, taking into account that $\overline{\Pi}^{2}_{1}$ and $\tilde{\Pi}^{2}_{1}$ are morphisms of Lie groupoids we can consider the induced morphisms of Lie algebroids $A\overline{\Pi}^{2}_{1} , A\tilde{\Pi}^{2}_{1}: A \tilde{J}^{1} \left( FM \right) \rightarrow A\Pi^{1} \left( M , M \right)$.\\
Thus, it is easy to induce the construction of the second-order non-holonomic prolongation over $ A \tilde{J}^{1} \left( FM \right)$. Given two section of $\sharp$
$$ A \mathcal{P}, A \mathcal{Q} : TM \rightarrow A \Pi^{1} \left( M,M \right),$$
we define the \textit{second-order non-holonomic prolongation associated to $A \mathcal{P}$ and $A \mathcal{Q}$} as follows,
$$ A \mathcal{P}^{1} \left( A \mathcal{Q}\right) \triangleq  A  \left( \mathcal{P}^{1} \left( \mathcal{Q} \right)\right) .$$
Then, $A \mathcal{P}^{1} \left( A \mathcal{Q}\right)$ projects via $A\overline{\Pi}^{2}_{1} $ (resp. $A\tilde{\Pi}^{2}_{1}$) over $A \mathcal{P}$ (resp. $A \mathcal{Q}$).\\
Using that the functor $A$ preserves integrability (see \cite{VMJIMM}), $A \mathcal{P}^{1} \left( A \mathcal{Q}\right)$ is said to be \textit{integrable} if $A \mathcal{Q}$ is integrable (equivalently $\mathcal{Q}$ is integrable). Therefore, if $A \mathcal{P}^{1} \left( A \mathcal{Q}\right)$ takes values in $A \tilde{j}^{1} \left( FM \right)$, $A \mathcal{P}^{1} \left( A \mathcal{Q}\right)$ is an integrable prolongation if, and only if, it is integrable.\\
Finally, we can introduce the following definition.
\begin{definition}\label{109}
\rm
Let $ A \tilde{J}^{1}_{\overline{G}} \left( FM \right)$ be a Lie subalgebroid of $ A \tilde{J}^{1} \left( FM \right)$. $ A \tilde{J}^{1}_{\overline{G}} \left( FM \right)$ is a non-holonomic integrable prolongation of second-order if we can cover $M$ by local non-holonomic integrable prolongations of second order which take values in $ A \tilde{J}^{1}_{\overline{G}} \left( FM \right)$.
\end{definition}

\begin{remarkth}\label{152}
\rm
Thus, $A \tilde{J}^{1}_{\overline{G}} \left( FM \right)$ is a non-holonomic integrable prolongation of second-order if, and only if, $\tilde{J}^{1}_{\overline{G}} \left( FM \right)$ is a non-holonomic integrable prolongation of second-order. Notice that, if $\tilde{J}^{1}_{\overline{G}} \left( FM \right)$ is a non-holonomic integrable prolongation of second-order then, we can cover $M$ by open sets $U$ and second-order non-holonomic integrable prolongations $\mathcal{P}^{1}\left(\mathcal{Q}\right) : U \times U \rightarrow \tilde{J}^{1}_{\overline{G}}\left(FU\right)$. However, we cannot take $A\mathcal{P}^{1}\left(\mathcal{Q}\right)$ because these sections do not have to be morphisms of Lie groupoids.\\
To solve this, we fix $ z_{0} \in M$ and define
$$ \mathcal{P}^{1}\left(\mathcal{Q}\right)^{z_{0}} \left( x , y \right) = \mathcal{P}^{1}\left(\mathcal{Q}\right) \left( z_{0} , y \right) \cdot [ \mathcal{P}^{1}\left(\mathcal{Q}\right) \left( z_{0} , x \right)]^{-1}, \ \forall x,y \in U.$$
Then, these family of sections are morphisms of Lie groupoids and non-holonomic integrable prolongations of second-order.\\
\end{remarkth}

Now, express this condition locally. Let $\overline{\mathcal{P}}: M \times M \rightarrow \tilde{J}^{1} \left( FM \right)$ be a section of $\left( \overline{\alpha} , \overline{\beta} \right)$ in $\tilde{J}^{1}\left( FM \right)$ and $\left( x^{i} \right)$ be a local coordinate system on $M$ such that
$$ \overline{\mathcal{P}} \left( x^{i} ,  y^{j} \right) = \left( \left( x^{i}, y^{j} , P^{j}_{i} \right) , Q^{j}_{i} , R^{j}_{i,k} \right).$$
Then,
$$ A \overline{\mathcal{P}} \left( x^{i} ,  \dfrac{\partial}{\partial x^{l}} \right) = \left( \left( x^{i}, \dfrac{\partial}{\partial x^{l}} , \dfrac{\partial  P^{j}_{i} }{\partial x^{l}}  \right) , \dfrac{\partial  Q^{j}_{i} }{\partial x^{l}} , \dfrac{\partial  R^{j}_{i,k} }{\partial x^{l}} \right),$$
where we are deriving fixing the first coordinate (see Eq. (\ref{112})).\\
In this way, take two section of $\sharp$, $ A \mathcal{P}$ and $ A \mathcal{Q}$, in $A \Pi^{1} \left( M , M \right)$ such that
$$ A \mathcal{P} \left(  x^{i} ,  \dfrac{\partial}{\partial x^{l}} \right) = \left( x^{i}, \dfrac{\partial}{\partial x^{l}} , \dfrac{\partial  P^{j}_{i} }{\partial x^{l}}  \right) \ ; \ \ \  A \mathcal{Q} \left( x^{i} ,  \dfrac{\partial}{\partial x^{l}} \right) = \left( x^{i}, \dfrac{\partial}{\partial x^{l}} , \dfrac{\partial  Q^{j}_{i} }{\partial x^{l}}  \right).$$
Hence,
\begin{equation}\label{110}
A \mathcal{P}^{1} \left( A \mathcal{Q} \right) \left( x^{i} ,  \dfrac{\partial}{\partial x^{l}} \right) = \left( \left( x^{i}, \dfrac{\partial}{\partial x^{l}} , \dfrac{\partial  P^{j}_{i} }{\partial x^{l}}  \right) ,  \dfrac{\partial  Q^{j}_{i} }{\partial x^{l}} , R^{j}_{i,kl}  \right),
\end{equation}
where
$$R^{j}_{i,kl}= \dfrac{\partial^{2} P^{j}_{i}}{\partial x^{l}\partial x^{k} } + \sum_{m} \dfrac{\partial Q_{k}^{m}}{\partial x^{l}} \dfrac{\partial P^{j}_{i}}{\partial y^{m} } +  Q_{k}^{m} \dfrac{\partial^{2} P^{j}_{i}}{\partial x^{l}\partial y^{m} }.$$

To understand why we obtain this local expression we have to take into account that we are fixing the second coordinate to do the induced map $A \mathcal{P}^{1} \left( A \mathcal{Q} \right)$.
Finally, using Eq. (\ref{110}), $A \mathcal{P}^{1} \left( A \mathcal{Q} \right) $ is integrable if, and only if, it can be locally expressed as follows
 
\begin{equation}\label{111}
A \mathcal{P}^{1} \left( A \mathcal{Q} \right) \left( x^{i} ,  \dfrac{\partial}{\partial x^{l}} \right) = \left( \left( x^{i}, \dfrac{\partial}{\partial x^{l}} , \dfrac{\partial  P^{j}_{i} }{\partial x^{l}}  \right) ,  0 ,   \dfrac{\partial^{2} P^{j}_{i}}{\partial x^{l}\partial x^{k} } +   \dfrac{\partial^{2} P^{j}_{i}}{\partial x^{l}\partial y^{k} } \right).
\end{equation}

So, we can rewrite Definition \ref{109} in the following way: Let \linebreak $ A \tilde{J}^{1}_{\overline{G}} \left( FM \right)$ be a Lie subalgebroid of $ A \tilde{J}^{1} \left( FM \right)$. $ A \tilde{J}^{1}_{\overline{G}} \left( FM \right)$ is a non-holonomic integrable prolongation of second-order if for each $x \in M$ there exists a local coordinate system $\left( x^{i} \right)$ over $x$ such that the local section of $\overline{\sharp}$,
$$A \mathcal{P}^{1} \left( A \mathcal{Q} \right) \left( x^{i} ,  \dfrac{\partial}{\partial x^{l}} \right) = \left( \left( x^{i}, \dfrac{\partial}{\partial x^{l}} , \dfrac{\partial  P^{j}_{i} }{\partial x^{l}}  \right) ,  0 ,   \dfrac{\partial^{2} P^{j}_{i}}{\partial x^{l}\partial x^{k} } +   \dfrac{\partial^{2} P^{j}_{i}}{\partial x^{l}\partial y^{k} } \right),$$
take values in $ A \tilde{J}^{1}_{\overline{G}} \left( FM \right)$.\\\\

\section{Characterization of homogeneity}
In this section we will use the general development made in the previous section to interpret the (local) homogeneity of Cosserat media in many different ways. So, let $F \mathcal{B}$ be a Cosserat medium with $\phi_{0}: \mathcal{B} \rightarrow \mathbb{R}^{3}$ as reference configuration and $W : \tilde{J}^{1} \left( F\mathcal{B}  \right) \rightarrow V$ as mechanical response. Consider $\overline{\Omega} \left( \mathcal{B} \right)$ the corresponding non-holonomic material groupoid of second order. Then, $\mathcal{B}$ is locally homogeneous if, and only if, for each point $x \in \mathcal{B}$ there exists an open subset $U \subseteq \mathcal{B}$ with $x \in U$ and a local deformation $\tilde{\kappa}$ over $U$ such that the (local) section $\overline{\mathcal{P}}: U \times U \rightarrow \tilde{J}^{1} \left( F\mathcal{B} \right)$ given by
$$ \overline{\mathcal{P}}\left(z,y \right) = j^{1}_{z,y} \left(\tilde{\kappa}^{-1} \circ F\tau_{\kappa\left(y\right) - \kappa \left(z\right)} \circ \tilde{\kappa}\right),$$
where $\tau_{\kappa\left(y\right) - \kappa \left(z\right)}: \mathbb{R}^{3} \rightarrow \mathbb{R}^{3}$ denotes the translation on $\mathbb{R}^{3}$ by the vector $\kappa\left(y\right) - \kappa \left(z\right)$ takes values in $\overline{\Omega} \left( \mathcal{B} \right)$ (see Definition \ref{151}).\\
So, using Proposition \ref{81}, we immediately have
\begin{proposition}\label{85}
Let $\mathcal{B}$ be a Cosserat medium. If $\mathcal{B}$ is homogeneous then $\overline{\Omega} \left( \mathcal{B}\right)$ is a second-order non-holonomic integrable prolongation. In fact, $\overline{\Omega} \left( \mathcal{B}\right)$ is a second-order non-holonomic integrable prolongation if, and only if, $\mathcal{B}$ is locally homogeneous.
\end{proposition}
Considering the local coordinates $\left(x^{i}\right) $ given by the deformations $\tilde{\kappa}$ satistying the deformation conditions, we deduce that $\mathcal{B}$ is locally homogeneous if, and only if, $\Omega \left( \mathcal{B} \right)$ can be locally covered by (local) sections of $\left( \overline{\alpha} , \overline{\beta} \right)$ in $\overline{\Omega} \left( \mathcal{B}\right)$ as follows:
$$\overline{\mathcal{P}}\left(x^{i},y^{j}\right) = \left(\left(x^{i},y^{j}, P^{j}_{i}\right) , \delta^{j}_{i} , \dfrac{\partial P^{j}_{i} }{\partial x^{k}} + \dfrac{\partial P^{j}_{i} }{\partial y^{k}} \right).$$

Next, let us consider the induced subalgebroid of the second-order non-holonomic material groupoid, $A \overline{\Omega} \left( \mathcal{B} \right)$. This Lie algebroid will be called \textit{second-order non-holonomic material algebroid of $\mathcal{B}$}.\\

Take $\Theta \in \Gamma  \left( A \overline{\Omega} \left( \mathcal{B}\right)\right)$. So, the flow of the left-invariant vector field $X_{\Theta}$, $\{ \varphi^{\Theta}_{t} \}$, can be restricted to $\overline{\Omega}\left( \mathcal{B} \right)$.\\
Hence, for any infinitesimal deformation $ g $, we have
$$ W\left( \varphi^{\Theta}_{t} \left( g \cdot \overline{\epsilon} \left( \overline{\alpha} \left(g \right) \right) \right) \right) = W \left(g \right).$$
Indeed, this equality is equivalent to the following one
\begin{equation}\label{92}
W \left( \varphi^{\Theta}_{t} \left( g \right) \right) = W \left(g \right), \forall g \in \tilde{J}^{1} \left( F \mathcal{B}  \right).
\end{equation}
Thus, for each $g \in \tilde{J}^{1} \left( F \mathcal{B} \right)$, we deduce
$$T W \left( X_{\Theta} \left( g \right)\right) = \dfrac{\partial}{\partial t_{|0}}\left(W \left(\varphi^{\Theta}_{t}\left( g \right) \right)\right)= \dfrac{\partial}{\partial t_{|0}}\left(W \left( g \right) \right) = 0.$$
Therefore,
\begin{equation}\label{93}
T W \left( X_{\Theta}\right) = 0.
\end{equation}

Conversely, suppose that Eq. (\ref{93}) is satisfied. Then,
$$\dfrac{\partial}{\partial t_{|s}}\left(W \left(\varphi^{\Theta}_{t}\left( g \right) \right)\right) = 0, \ \forall g \in \tilde{J}^{1} \left( F \mathcal{B} \right), \ \forall s.$$
Thus, taking into account that
$$ W \left(\varphi^{\Theta}_{0}\left(g \right) \right) = W \left(g \right),$$
we have
$$W \left(\varphi^{\Theta}_{t}\left(g \right) \right) = W \left(g \right),$$
i.e.,
$$ \Theta \in \Gamma \left( A \overline{\Omega} \left( \mathcal{B} \right)\right).$$

In this way, the second-order non-holonomic material algebroid can be defined without using the non-holonomic material groupoid of second order. Thus, we can characterize the homogeneity and uniformity using the material Lie algebroid. Taking into account that the fact of being an ``integrable prolongation" can be equivalently defined over the associated Lie algebroid (see Remark \ref{152}) we get the following result:

\begin{proposition}\label{95}
Let $\mathcal{B}$ be a Cosserat continuum. If $\mathcal{B}$ is homogeneous, then, $A\overline{\Omega} \left( \mathcal{B} \right)$ is an integrable non-holonomic prolongation of second order. Conversely, $A\overline{\Omega} \left( \mathcal{B} \right)$ is an integrable non-holonomic prolongation of second order implies that $\mathcal{B}$ is locally homogeneous.

\end{proposition}

Using the local expression (\ref{111}), this result can be expressed locally as follows.

\begin{proposition}\label{96}
Let $\mathcal{B}$ be a Cosserat continuum. $\mathcal{B}$ is locally homogeneous if, and only if, for each point $x  \in \mathcal{B}$ there exists a local coordinate system $\left(x^{i}\right)$ over $U \subseteq \mathcal{B}$ with $x \in U$ such that the local section of $\overline{\sharp}$, 
$$A \mathcal{P}^{1} \left( A \mathcal{Q} \right) \left( x^{i} ,  \dfrac{\partial}{\partial x^{l}} \right) = \left( \left( x^{i}, \dfrac{\partial}{\partial x^{l}} , \dfrac{\partial  P^{j}_{i} }{\partial x^{l}}  \right) ,  0 ,   \dfrac{\partial^{2} P^{j}_{i}}{\partial x^{l}\partial x^{k} } +   \dfrac{\partial^{2} P^{j}_{i}}{\partial x^{l}\partial y^{k} } \right),$$
takes values in $A \overline{\Omega}\left( \mathcal{B}\right)$.
\end{proposition}

Finally, we will use Theorem \ref{48} to give another characterization of the homogeneity. Indeed, let be $\tilde{\Lambda} : T\mathcal{B} \rightarrow A\tilde{J}^{1} \left( F\mathcal{B} \right)$ a linear section of $\overline{\sharp}$. Then, using Remark \ref{114}, $\tilde{\Lambda}$ can be seen as a map
$$ \nabla^{\Lambda} : \mathfrak{X} \left( \mathcal{B} \right) \times \mathfrak{X} \left( F \mathcal{B} \right) \rightarrow \mathfrak{X} \left( F \mathcal{B} \right),$$
where, for all $\left( X , \tilde{Y}\right) \in \mathfrak{X} \left( \mathcal{B} \right) \times \mathfrak{X} \left( F\mathcal{B} \right)$, $f \in \mathcal{C}^{\infty} \left( \mathcal{B} \right)$ and $F   \in \mathcal{C}^{\infty} \left( F\mathcal{B} \right)$ satisfies that
\begin{itemize}
\item[(i)] $ \nabla^{\Lambda}_{fX}\tilde{Y} = \left( f \circ \pi_{\mathcal{B}} \right) \nabla^{\Lambda}_{X}\tilde{Y}.$
\item[(ii)] $ \nabla^{\Lambda}_{X}F\tilde{Y} = F \nabla^{\Lambda}_{X}\tilde{Y} + \Lambda\left( X \right)^{\sharp} \left( F \right) \tilde{Y}.$
\item[(iii)] The base vector field of $\nabla^{\Lambda}_{X}$ is $\Lambda \left( X \right)^{\sharp}$ which is $\pi_{\mathcal{B}}-$related to $X$.
\item[(iv)] For all $g \in Gl \left( 3, \mathbb{R} \right)$,
$$\nabla_{X}^{\Lambda} \circ TR^{*}_{g} = TR_{g}^{*} \circ \nabla_{X}^{\Lambda}.$$
\item[(v)] The flow of $\nabla^{\Lambda}_{X}$ is the tangent map of an automorphism of frame bundles (over the identity map) at each fibre.

\end{itemize}
Let $\left( x^{i}\right)$ be a local coordinate system on $\mathcal{B}$ such that
$$ \tilde{\Lambda} \left( x^{i} ,  \dfrac{\partial}{\partial x^{l}} \right) = \left(  x^{i} ,  \dfrac{\partial}{\partial x^{l}} , \Lambda^{j}_{il} , \Lambda^{j}_{,il} , \Lambda^{j}_{i,kl} \right).$$
Then, $\nabla^{\Lambda}$ is locally characterized as follows
\begin{itemize}
\item[(i)] $\nabla^{\Lambda}_{ \dfrac{\partial}{\partial x^{j}}} \dfrac{\partial}{\partial x^{i}} = \sum_{k} \Lambda^{k}_{,ij} \dfrac{\partial}{\partial x^{k}} + \sum_{k,l} \Lambda^{k}_{l,ij} \dfrac{\partial}{\partial x^{k}_{l}}$
\item[(ii)] $\nabla^{\Lambda}_{\dfrac{\partial}{\partial x^{k}}}  \dfrac{\partial}{\partial x^{i}_{j}} = \sum_{l} \Lambda^{l}_{ik}\dfrac{\partial}{\partial x^{l}_{j}}$
\end{itemize}
In this way, if $\tilde{\Lambda} = A \mathcal{P}^{1} \left( A \mathcal{Q} \right)$ is a non-holonomic prolongation of second order we have that
\begin{itemize}
\item[(i)] $\nabla^{A \mathcal{P}^{1} \left( A \mathcal{Q} \right)}_{ \dfrac{\partial}{\partial x^{j}}} \dfrac{\partial}{\partial x^{i}} = \sum_{k} \dfrac{\partial  Q^{k}_{i} }{\partial x^{j}}  \dfrac{\partial}{\partial x^{k}} + \sum_{k,l} R^{k}_{l,ij} \dfrac{\partial}{\partial x^{k}_{l}}$
\item[(ii)] $\nabla^{A \mathcal{P}^{1} \left( A \mathcal{Q} \right)}_{\dfrac{\partial}{\partial x^{k}}}  \dfrac{\partial}{\partial x^{i}_{j}} = \sum_{l} \dfrac{\partial  P^{l}_{i} }{\partial x^{k}} \dfrac{\partial}{\partial x^{l}_{j}},$
\end{itemize}
where 
$$R^{k}_{l,ij}= \dfrac{\partial^{2} P^{k}_{l}}{\partial x^{j}\partial x^{i} } + \sum_{m} \dfrac{\partial Q_{i}^{m}}{\partial x^{j}} \dfrac{\partial P^{k}_{l}}{\partial x^{m} } +  Q_{i}^{m} \dfrac{\partial^{2} P^{k}_{l}}{\partial x^{j}\partial x^{m} }.$$
Hence, $\Lambda$ is an integrable non-holonomic prolongation of second order if an only if
\begin{itemize}
\item[(i)] $\nabla^{A \mathcal{P}^{1} \left( A \mathcal{Q} \right)}_{ \dfrac{\partial}{\partial x^{j}}} \dfrac{\partial}{\partial x^{i}} =  \sum_{k,l} R^{k}_{l,ij} \dfrac{\partial}{\partial x^{k}_{l}}$
\item[(ii)] $\nabla^{A \mathcal{P}^{1} \left( A \mathcal{Q} \right)}_{\dfrac{\partial}{\partial x^{k}}}  \dfrac{\partial}{\partial x^{i}_{j}} = \sum_{l} \dfrac{\partial  P^{l}_{i} }{\partial x^{k}} \dfrac{\partial}{\partial x^{l}_{j}},$
\end{itemize}
where 
$$R^{k}_{l,ij}= \dfrac{\partial^{2} P^{k}_{l}}{\partial x^{j}\partial x^{i} } + \dfrac{\partial^{2} P^{k}_{l}}{\partial x^{j}\partial x^{i} }.$$\\
Using this we can give the following result:

\begin{proposition}\label{115}
Let $\mathcal{B}$ be a Cosserat continuum. $\mathcal{B}$ is locally homogeneous if, and only if, for each point $x  \in \mathcal{B}$ there exists a local coordinate system $\left(x^{i}\right)$ over $U \subseteq \mathcal{B}$ with $x \in U$ such that the local non-holonomic covariant derivative of second order $\nabla$ satisfies
\begin{itemize}
\item[(i)] $\nabla_{ \dfrac{\partial}{\partial x^{j}}} \dfrac{\partial}{\partial x^{i}} = \sum_{k,l} R^{k}_{l,ij} \dfrac{\partial}{\partial x^{k}_{l}}$
\item[(ii)] $\nabla_{\dfrac{\partial}{\partial x^{k}}}  \dfrac{\partial}{\partial x^{i}_{j}} = \sum_{l} \dfrac{\partial  P^{l}_{i} }{\partial x^{k}} \dfrac{\partial}{\partial x^{l}_{j}}$
\end{itemize}
where 
$$R^{k}_{l,ij}=  \dfrac{\partial^{2} P^{k}_{l}}{\partial x^{j}\partial x^{i} } + \dfrac{\partial^{2} P^{k}_{l}}{\partial x^{j}\partial x^{i} },$$
takes values in $\mathcal{D} \left( A \overline{\Omega}\left( \mathcal{B}\right) \right)$.
\end{proposition}
Let $\tilde{\Lambda} : T \mathcal{B} \rightarrow A\tilde{J}^{1} \left( F\mathcal{B} \right)$ be a linear section of $\overline{\sharp}$ and $\nabla^{\tilde{\Lambda}}$ its associated covariant derivative. We can construct $\left( \Lambda_{1} \right)^{1} \left( \Lambda_{2} \right)$ where
$$ \Lambda_{1} \triangleq A\overline{\Pi}^{2}_{1} \circ \tilde{\Lambda} \ \ \ \ ; \ \ \ \  \Lambda_{2} \triangleq A\tilde{\Pi}^{2}_{1} \circ \tilde{\Lambda}.$$
Then, $\left( \Lambda_{1} \right)^{1} \left( \Lambda_{2} \right)$ has an associated map $\nabla^{\left( \Lambda_{1} \right)^{1} \left( \Lambda_{2} \right)} :  \mathfrak{X} \left( \mathcal{B} \right) \times \mathfrak{X} \left( F \mathcal{B} \right) \rightarrow \mathfrak{X} \left( F \mathcal{B} \right)$ which satisfies $\left( i \right)$, $\left( i i \right)$, $\left( i ii\right)$, $\left( i v \right)$ and $\left( v \right)$. So, $\tilde{\Lambda}$ is a prolongation if, and only if,
$$ \nabla^{\Lambda} = \nabla^{\left( \Lambda_{1} \right)^{1} \left( \Lambda_{2} \right)} .$$
On the other hand, using that for all $X \in \frak X \left( M \right)$, $\nabla^{\Lambda}_{X}$ is $TR_{g}^{*}-$invariant we have that $\nabla^{\Lambda}_{X}$ preserves right-invariant vector fields on $FM$.\\
Then, we can project $\nabla^{\Lambda}$ onto a covariant derivative on $M$, $\nabla^{1}$, in the following way: Let $X,Y$ be two vector fields on $M$ and $Y^{c}$ the complete lift of $Y$ over $FM$ (see \cite{LACMDL}). Then, $Y^{c}$ is right-invariant which implies that $ \nabla^{\Lambda}_{X} Y^{c}$ is right-invariant. So, $\nabla^{\Lambda}_{X} Y^{c}$ projects onto a unique vector field on $M$. This vector field will be $\nabla^{1}_{X}Y$. It is straightforward to prove that $\nabla^{1}$ is a covariant derivative over $M$; indeed, let $\left( x^{i} \right)$ be a local coordinate system on $\mathcal{B}$ such that
$$ \tilde{\Lambda} \left( x^{i} ,  \dfrac{\partial}{\partial x^{l}} \right) = \left(  x^{i} ,  \dfrac{\partial}{\partial x^{l}} , \Lambda^{j}_{il} , \Lambda^{j}_{,il} , \Lambda^{j}_{i,kl} \right).$$
Then, $\nabla^{1}$ satisfies that
\begin{equation}
\nabla^{1}_{ \dfrac{\partial}{\partial x^{j}}} \dfrac{\partial}{\partial x^{i}} = \sum_{k} \Lambda^{k}_{,ij} \dfrac{\partial}{\partial x^{k}}.
\end{equation}
Hence, suppose that $\tilde{\Lambda}$ is a non-holonomic prolongation of second order. $\tilde{\Lambda}$ is an integrable prolongation if, and only if, $\nabla^{1}$ is locally trivial, i.e., the Christoffel symbols are zero.
There is an alternative way to construct $\nabla^{1}$: Using Theorem \ref{48}, we can construct a covariant derivative on $M$, $\nabla^{\Lambda_{2}}$, such that
$$ \nabla^{\Lambda_{2}} = \nabla^{1}.$$
To summarize, we have introduced a new frame (groupoids and Lie algebroid) to study Cosserat media. In this frame, we have been able to express the homogeneity in many different (but equivalent) ways: Over the non-holonomic material groupoid of second order, over the associated Lie algebroid (which can be contructed without using the material groupoid) and over the Lie algebroid of derivations. Finally, using the interpretation over the algebroid of derivations, we have developed a method to know if a covariant derivative is a non-holonomic integrable prolongation without using coordinates.
\section{Homogeneity with non-holonomic $\overline{G}-$structures of second order}
Now, we will introduce the definition of homogeneity used in \cite{MAREMDL} where the authors discuss second-order non-holonomic $\overline{G}-$structures. We will fix $\mathcal{B}$ a Cosserat continuum.\\
Let be $\overline{Z}_{0} = j^{1}_{e_{1},Z_{0}} \Phi \in \overline{F}^{2} \mathcal{B}$ a non-holonomic frame of second order at $z_{0}$. Define a non-holonomic $\overline{G}_{0}-$structure of second order $\overline{\omega}_{\overline{G}_{0}} \left( \mathcal{B}\right)$ on $\mathcal{B}$ (which contains $\overline{Z}_{0}$) as follows:

$$\overline{\omega}_{\overline{G}_{0}} \left( \mathcal{B}\right) \triangleq \overline{\Omega}_{z_{0}} \left( \mathcal{B}\right) \cdot \overline{Z}_{0}.$$

Notice that the principal bundles $\overline{\Omega}_{z_{0}} \left( \mathcal{B}\right)$ and $\overline{\omega}_{\overline{G}_{0}} \left(\mathcal{B}\right)$ are obviously isomorphic (as principal bundles). Indeed, the structure group of $\overline{\omega}_{\overline{G}_{0}} \left( \mathcal{B}\right)$ is given by
$$ \overline{G}_{0} \triangleq \overline{Z}_{0}^{-1} \cdot \overline{G}\left(z_{0} \right) \cdot \overline{Z}_{0},$$
where $\overline{G} \left( z_{0} \right)$ is the isotropy group of $\overline{\Omega} \left( \mathcal{B} \right)$ at $z_{0}$, i.e., the family of all material symmetries at $z_{0}$.\\
A local section of $\overline{\omega}_{\overline{G}_{0}} \left( \mathcal{B}\right)$ will be called \textit{local uniform reference}. A global section of $\overline{\omega}_{\overline{G}_{0}} \left( \mathcal{B}\right)$ will be called \textit{global uniform reference}. We call \textit{reference crystal} to any frame $\overline{Z}_{0} \in \overline{F}^{2} \mathcal{B}$ at $z_{0}$.\\

Now, the canonical projection of the second-order non-holonomic $\overline{G}_{0}-$structure $\overline{\omega}_{\overline{G}_{0}} \left( \mathcal{B}\right)$ is a $G_{0}-$structure denoted by $\omega_{G_{0}} \left( \mathcal{B}\right)$.

\begin{remarkth}\label{118}
\rm
\begin{itemize}
\item[(1)] If we change the point $z_{0}$ to another point $z_{1}$ then we can obtain the same second-order non-holonomic \\$\overline{G}_{0}-$structure. We only have to take a frame $\overline{Z}_{1}$ as the composition of $\overline{Z}_{0}$ with a $j^{1}_{z_{0},z_{1}} \Psi \in \overline{G}\left(z_{0}, z_{1}\right)$.
\item[(2)] We have fixed a configuration $\Phi_{0}$. Suppose that $\Phi_{1}$ is another reference configuration such that the change of configuration (or deformation) is given by $\tilde{\kappa} = \Phi_{1}^{-1} \circ \Phi_{0}$. Transporting the reference crystal $\overline{Z}_{0}$ via $F\tilde{\kappa}$ we get another reference crystal such that the second-order non-holonomic $\overline{G}_{0}-$structures are isomorphic.
\item[(3)] Finally suppose that we have another crystal reference $\overline{Z}'_{0}$ at $z_{0}$. Hence, the new second-order non-holonomic $\overline{G}_{0}'-$structure, $\overline{\omega}_{\overline{G}_{0}'} \left(\mathcal{B}\right)$, is conjugated of $\overline{\omega}_{\overline{G}_{0}} \left( \mathcal{B}\right)$, namely,
$$\overline{G}_{0}' =\overline{g} \cdot \overline{G}_{0} \cdot \overline{g}^{-1}, \ \ \ \overline{\omega}_{\overline{G}_{0}'} \left(\mathcal{B}\right) = \overline{\omega}_{\overline{G}_{0}} \left( \mathcal{B}\right) \cdot \overline{g},$$
for $\overline{g} = \overline{Z}_{0}^{-1} \cdot \overline{Z}'_{0}$.
\end{itemize}
\end{remarkth}

In this way, the definition of homogeneity is the following,

\begin{definition}
\rm
A Cosserat continuum $\mathcal{B}$ is said to be \textit{homogeneous} with respect to the crystal reference $\overline{Z}_{0}$ if it admits a global deformation $\overline{\kappa}$ such that $\overline{\kappa}^{-1}$ induces a uniform reference $\overline{P}$, i.e., for each $x \in \mathcal{B}$
$$ \overline{P}\left(x\right) = j^{1}_{0,x} \left(\overline{\kappa}^{-1} \circ F\tau_{\kappa\left(x\right)}\right),$$
where $\tau_{\kappa \left(x\right)}: \mathbb{R}^{3} \rightarrow \mathbb{R}^{3}$ denotes the translation on $\mathbb{R}^{3}$ by the vector $\kappa \left(x\right)$ and $\kappa$ is the induced map of $\overline{\kappa}$ over $\mathcal{B}$. $\mathcal{B}$ is said to be \textit{locally homogeneous} if every $x \in \mathcal{B}$ has a neighbourhood which is homogeneous.
\end{definition}
Using Eq. (\ref{37}) it is easy to prove the following result:

\begin{proposition}\label{119}
If $\mathcal{B}$ is homogeneous with respect to $\overline{Z}_{0}$ then $\overline{\omega}_{\overline{G}_{0}} \left(\mathcal{B}\right)$ is a non-holonomic integrable prolongation of second order. Conversely, $\overline{\omega}_{\overline{G}_{0}} \left(\mathcal{B}\right)$ is a second-order non-holonomic integrable prolongation implies that $\mathcal{B}$ is locally homogeneous with respect to $\overline{Z}_{0}$.
\end{proposition}

Notice that, this result shows us that the homogeneity does not depend on the point and reference configuration but depends on the reference crystal (see Remark \ref{118}).\\
It is important to recall that our definition of homogeneity does not depend on a crystal reference. So, these definitions cannot be equivalent (as a difference with simple media). However, there are a close relation. In fact, $\mathcal{B}$ is homogeneous (resp. locally homogeneous) if, and only if, there exists a crystal reference $\overline{Z}_{0}$ such that $\mathcal{B}$ is homogeneous (resp. locally homogeneous) with respect to $\overline{Z}_{0}$.\\
To prove this, we will begin defining the following map

$$
\begin{array}{rccl}
\overline{\mathcal{G}} : & \Gamma \left( \overline{F}^{2}M\right)  & \rightarrow & \Gamma_{\left(\overline{\alpha} , \overline{\beta}\right)}\left(\tilde{J}^{1}\left(FM\right)\right)\\
&\overline{P}   &\mapsto &  \overline{\mathcal{G}} \overline{P},
\end{array}
$$
such that
$$\overline{\mathcal{G}} \overline{P} \left(x,y\right) = \overline{P}\left(y\right) \cdot [\overline{P}\left(x\right)]^{-1}, \ \forall x,y \in M,$$
where we are considering the equivalence class in $\tilde{J}^{1}\left(FM\right)$.\\

\begin{remarkth}\label{120}
\rm
Notice that the following equalities relate $\overline{\Pi}^{2}_{1}$ (resp. $\tilde{\Pi}^{2}_{1}$) with $\overline{\pi}^{2}_{1}$ (resp. $\tilde{\pi}^{2}_{1}$)
\begin{itemize}
\item[(i)] $\overline{\Pi}^{2}_{1} \circ \overline{\mathcal{G}} = \mathcal{G} \circ \overline{\pi}^{2}_{1}$

\item[(ii)] $\tilde{\Pi}^{2}_{1} \circ \overline{\mathcal{G}} = \mathcal{G} \circ \tilde{\pi}^{2}_{1}$

\end{itemize}
where $\mathcal{G} : \Gamma \left( FM\right) \rightarrow \Gamma_{\left(\alpha , \beta\right)}\left(\Pi^{1}\left(M,M\right)\right)$ has been defined in \cite{VMJIMM} as follows:
$$\mathcal{G} P \left(x,y\right) = P\left(y\right) \cdot [P\left(x\right)]^{-1}, \ \forall x,y \in M.$$
\end{remarkth}

Before working with second-order non-holonomic prolongations, we are interested in knowing when an element of $\Gamma_{\left(\overline{\alpha} , \overline{\beta}\right)}\left(\tilde{J}^{1}\left(FM\right)\right)$ can be inverted by $\overline{\mathcal{G}}$. First, we consider $\overline{P} \in \Gamma \left( \overline{F}^{2}M\right)$; then for all $x,y,z \in M$
\begin{equation}\label{121}
\overline{\mathcal{G}} \overline{P} \left(y,z\right) \cdot  \overline{\mathcal{G}} \overline{P} \left(x,y\right)  = \overline{\mathcal{G}} \overline{P} \left(x,z\right),
\end{equation}
i.e., $\overline{\mathcal{G}}\overline{P}$ is a morphism of Lie groupoids over the identity map on $M$ from the pair groupoid $M \times M$ to $\tilde{J}^{1}\left(FM\right)$. Therefore, not all element of $\Gamma_{\left(\overline{\alpha} , \overline{\beta}\right)}\left(\tilde{J}^{1}\left(FM\right)\right)$ can be inverted by $\overline{\mathcal{G}}$ but we can prove the following result:

\begin{proposition}\label{122}
Let $\overline{\mathcal{P}}$ be a section of $\tilde{J}^{1}\left(FM\right)$. There exists a section of $\overline{F}^{2}  M$ such that
\begin{equation}\label{144}
\overline{\mathcal{G}} \overline{P} = \overline{\mathcal{P}},
\end{equation}
if, and only if, $\overline{\mathcal{P}}$ is a morphism of Lie groupoids over the identity map from the pair groupoid $M \times M$ to $\tilde{J}^{1}\left(FM\right)$.
\begin{proof}
We have already proved that Eq. (\ref{144}) implies that $\overline{\mathcal{P}}$ is a morphism of Lie groupoids over the identity map from the pair groupoid $M \times M$ to $\tilde{J}^{1}\left(FM\right)$. Conversely, if Eq. (\ref{121}) is satisfied we can define $\overline{P} \in \Gamma \left(\overline{F}^{2}M\right)$ as follows
$$ \overline{P}\left(x\right) =\overline{\mathcal{P}}\left(z_{0},x\right) \cdot \overline{Z}_{0},$$
where $\overline{Z}_{0} \in \overline{F}^{2}M$ with $\overline{\pi}^{2} \left( \overline{Z}_{0} \right)=z_{0}$ is fixed. Then, using Eq. (\ref{121}), we have
$$ \overline{\mathcal{G}} \overline{P} = \overline{\mathcal{P}}.$$
\end{proof}
\end{proposition}
However, there is not a unique $\overline{P}$ such that $ \overline{\mathcal{G}} \overline{P} = \overline{\mathcal{P}}$. In fact, let $\overline{P}$ and $\overline{Q}$ be sections of $\overline{F}^{2}M$ such that
$$ \overline{\mathcal{G}} \overline{P} = \overline{\mathcal{G}} \overline{Q}.$$
Then, there exists $\overline{g} \in \overline{G}^{2}\left(n\right)$ such that
\begin{equation}\label{127}
\overline{P} = \overline{Q} \cdot \overline{g},
\end{equation}
where, we are choosing representatives of the equivalence class to do the jet composition.\\
Notice that the sections of $\left(\overline{\alpha} , \overline{\beta}\right)$ which are morphisms of Lie groupoids over the identity map from the pair groupoid $M \times M$ to $\tilde{J}^{1}\left(FM\right)$ are, precisely, the parallelisms.\\\\

Now, lets see that $\overline{\mathcal{G}} : \Gamma \left( \overline{F}^{2}M \right)   \rightarrow  \Gamma_{\left(\overline{\alpha}, \overline{\beta}\right)} \left(\tilde{J}^{1}\left(FM\right)\right)$ preserves prolongations. In fact, let $P^{1}\left(Q\right) \in \Gamma \left( \overline{F}^{2}M\right)$ be a second-order non-holonomic prolongation of $P$ and $Q$, then 
\begin{equation}\label{117}
\overline{\mathcal{G}} P^{1}\left(Q\right) = \mathcal{G}P^{1} \left( \mathcal{G}Q\right),
\end{equation}
To prove the last equality, we use that
$$ P^{1}\left(Q\right) \left(x\right) = j_{e_{1}, P\left(x\right)}^{1} \left(\overline{P \circ \psi_{x}}\right),$$
where $Q\left(x\right) = j_{0,x}^{1} \psi_{x}$ (see Remark \ref{36}). Then,
$$ j_{0,y}^{1}\left( \overline{P  \circ \psi_{y}} \right) \cdot j_{x,0}^{1} \left( \overline{P  \circ \psi_{x}}\right)^{-1} = j_{x,y}^{1}\left( \overline{\mathcal{G}P \circ \left(\psi_{y}\circ  \psi_{x}^{-1}\right)}\right),$$
and, as we know
$$\mathcal{G} Q\left(x,y\right) = j_{x,y}^{1} \left(\psi_{y} \circ \psi_{x}^{-1}\right).$$

Then, taking into account that $\mathcal{G}$ preserves integrability (see \cite{VMJIMM}), we can assume that $\overline{\mathcal{G}}$ preseves integrable sections and integrable prolongations of $\overline{F}^{2}M$.\\
Conversely, we want to study if we can invert integrable sections (resp. non-holonomic integrable prolongations of second-order) in $\tilde{J}^{1} \left(FM\right)$. Notice that both kinds of sections can be written as second-order non-holonomic prolongations and, in this way, we will study when we can invert non-holonomic prolongations of second-order.\\
So, let $\mathcal{P}^{1}\left( \mathcal{Q}\right)$ be a second-order non-holonomic prolongation in $\tilde{J}^{1} \left( FM\right)$. Using Eq. (\ref{117}) and Remark \ref{120}, if we can invert $\mathcal{P}^{1}\left( \mathcal{Q}\right)$ then, there exist $P , Q \in \Gamma \left( FM\right)$ such that
$$\overline{\mathcal{G}} P^{1}\left(Q\right) = \mathcal{P}^{1} \left( \mathcal{Q}\right).$$
Therefore, analogously to Proposition \ref{122}, $\mathcal{P}$ and $\mathcal{Q}$ have to be Lie groupoid morphisms from the pair groupoid $M \times M$ to $\Pi^{1}\left(M,M\right)$.

\begin{proposition}\label{123}
Let $\mathcal{P}^{1} \left( \mathcal{Q}\right)$ be a second-order non-holonomic prolongation in $\tilde{J}^{1}\left(FM\right)$. There exists a second-order non-holonomic prolongation in $\overline{F}^{2} M$ such that
$$ \overline{\mathcal{G}} P^{1}\left(Q\right) = \mathcal{P}^{1} \left( \mathcal{Q}\right),$$
if, and only if, $\mathcal{P}$ and $\mathcal{Q}$ are morphisms of Lie groupoids from the pair groupoid $M \times M$ to $\Pi^{1}\left(M , M\right)$.
\end{proposition}
Now, notice that, by construction, every integral section of $ \Pi^{1} \left( M , M\right)$ is a morphism of Lie groupoids from the pair groupoid $M \times M$ to $\Pi^{1}\left(M , M\right)$. So we can state the following result:

\begin{corollary}\label{124}
Let $\mathcal{P}^{1} \left( \mathcal{Q}\right)$ be a second-order non-holonomic prolongation in $\tilde{J}^{1}\left(FM\right)$.
\begin{itemize}
\item[(i)] If $\mathcal{P}^{1} \left( \mathcal{Q}\right)$ is integrable then, there exists an integrable section of $\overline{F}^{2}M$, $P^{1}\left(Q\right)$, such that
$$\overline{\mathcal{G}} P^{1}\left(Q\right) = \mathcal{P}^{1} \left( \mathcal{Q}\right).$$
\item[(ii)] If $\mathcal{P}^{1} \left( \mathcal{Q}\right)$ is a non-holonomic integrable prolongation of second order then, there exists a second-order non-holonomic integrable prolongation of $\overline{F}^{2}M$, $P^{1}\left(Q\right)$, such that
$$\overline{\mathcal{G}} P^{1}\left(Q\right) = \mathcal{P}^{1} \left( \mathcal{Q}\right),$$
if, and only if, $P$ is a morphism of Lie groupoids from the pair groupoid $M \times M$ to $\Pi^{1}\left(M , M\right)$.
\end{itemize}

\end{corollary}

This result could induce us to think that if $\overline{P}$ satisfies
$$\overline{\mathcal{G}} \overline{P} = \mathcal{P}^{1} \left( \mathcal{Q}\right),$$
then, $\overline{P}$ is a second-order non-holonomic prolongation but this is not true. We just have to impose the condition locally to prove it.\\

Finally, as for the simple media case, we can generalize the map $\overline{\mathcal{G}}$ to a map which takes non-holonomic $\overline{G}-$structures of second order on $M$ into reduced subgroupoids of $\tilde{J}^{1}\left(FM\right)$. Let $\overline{\omega}_{\overline{G}} \left(M\right)$ be a non-holonomic $\overline{G}-$structure of second order on $M$, then we consider the following set,
$$\overline{\mathcal{G}}\left(\overline{\omega}_{\overline{G}}\left(M\right)\right) = \{ L_{y}\cdot [L_{x}^{-1}]\ / \ L_{x}, L_{y} \in \overline{\omega}_{\overline{G}}\left(M\right)\},$$
where we are considering the equivalence class in $\tilde{J}^{1}\left(FM\right)$. We will denote $\overline{\mathcal{G}}\left(\overline{\omega}_{\overline{G}}\left(M\right)\right)$ by $\tilde{J}^{1}_{\overline{G}} \left(FM\right)$. $\tilde{J}^{1}_{\overline{G}} \left(FM\right)$ is a reduced subgroupoid of $\tilde{J}^{1} \left(FM\right)$. In fact, taking a local section of $\overline{\omega}_{\overline{G}}\left(M\right)$, 
$$\overline{P}_{U}: U \rightarrow \overline{\omega}_{\overline{G}} \left(U\right),$$
the map given by
$$
\begin{array}{rccl}
F_{U}  : & \tilde{J}^{1}\left(FU\right) & \rightarrow & \overline{F}^{2}U \times U\\
& \tilde{L}_{x,y}  &\mapsto &  \left( \tilde{L}_{x,y} \cdot [\overline{P}_{U}\left(x\right)], x \right)
\end{array}
$$
is a diffeomorphism which satisfies that $F_{U} \left(\tilde{J}^{1}_{\overline{G}} \left(FU\right)\right)= \overline{\omega}_{\overline{G}} \left(U\right) \times U$.\\
Analogously to parallelisms, we can prove that every reduced subgroupoid can be inverted by $\overline{\mathcal{G}}$ in a non-holonomic $\overline{G}-$structure of second order on $M$.\\
Fix $z_{0} \in M$ and $\overline{Z}_{0} \in \overline{F}^{2}M$ with $\overline{\pi}^{2} \left(\overline{Z}_{0}\right) = z_{0}$. Then, we define
\begin{equation}\label{128}
\overline{G} := \{ \overline{Z}_{0}^{-1}\cdot  \overline{g}_{z_{0}} \cdot \overline{Z}_{0} \ / \ \overline{g}_{z_{0}} \in \tilde{J}^{1}_{\overline{G}}\left(z_{0}\right) \} = \overline{Z}_{0}^{-1} \cdot \tilde{J}^{1}_{\overline{G}} \left(z_{0}\right) \cdot \overline{Z}_{0} \cong  \tilde{J}^{1}_{\overline{G}} \left(z_{0}\right)
\end{equation}
where $ \tilde{J}^{1}_{\overline{G}}\left(z_{0}\right)$ is the isotropy group of $\tilde{J}^{1}_{\overline{G}}\left(FM\right)$ over $z_{0}$. Therefore, $\overline{G}$ is clearly a Lie subgroup of $\overline{G}^{2}\left(n\right)$.\\
Then, we can generate a second-order non-holonomic $\overline{G}-$structure over $M$ in the following way
$$ \overline{\omega}_{\overline{G}} \left(M\right) := \{ L_{z_{0}, x} \cdot \overline{Z}_{0} \cdot \overline{g} \ / \ \overline{g} \in \overline{G}, \ L_{z_{0},x} \in \tilde{J}^{1}_{\overline{G}}\left(FM\right)_{z_{0}}\} \triangleq \tilde{J}^{1}_{\overline{G}}\left(FM\right)_{z_{0}} \cdot \overline{Z}_{0} .$$
Notice that $\tilde{J}^{1}_{\overline{G}} \left(FM\right)_{z_{0}}$ and $\overline{\omega}_{\overline{G}}\left(M\right)$ are clearly isomorphic.\\

Next, let $\overline{\omega}_{\overline{G}} \left(M\right)$ be an integrable second-order non-holonomic (resp. integrable prolongation) $\overline{G}-$structure on $M$. Using Proposition \ref{31}, Proposition \ref{104} and the fact of that $\overline{\mathcal{G}}$ preserves integrable sections (resp. integrable prolongations) we have that $\tilde{J}^{1}_{\overline{G}} \left(FM\right)$ is integrable (resp. an integrable prolongation).\\

Conversely, let $\tilde{J}^{1}_{\overline{G}} \left( FM \right)$ be an integrable (resp. integrable prolongation) Lie subroupoid of $\tilde{J}^{1} \left( FM \right)$. Then, we may construct an integrable (resp. integrable prolongation) second-order non-holonomic $\overline{G}-$structure $\overline{\omega}_{\overline{G}} \left(M\right)$ such that
\begin{equation}\label{149}
\overline{\mathcal{G}}\left(\overline{\omega}_{\overline{G}}\left(M\right)\right) = \tilde{J}^{1}_{\overline{G}} \left( FM \right).
\end{equation}
To do this we just have to use Proposition \ref{70} (resp. Proposition \ref{147}) and define it locally.
However, as you can imagine, not all  second-order non-holonomic $\overline{G}-$structure wich satisfies Eq. (\ref{149}) is integrable.

Now we have what we need to prove the relation which we had mentioned. Let be $\mathcal{B}$ a Cosserat continuum and a crystal frame $\overline{Z}_{0} \in \overline{F}^{2} \mathcal{B}$ at $z_{0}$. Then, we have defined the second-order non-holonomic $\overline{G}_{0}-$structure of uniform references as follows
$$\overline{\omega}_{\overline{G}_{0}} \left( \mathcal{B}\right) = \overline{\Omega}_{z_{0}} \left( \mathcal{B}\right) \cdot \overline{Z}_{0},$$
i.e.,
$$\overline{\mathcal{G}}\left(\overline{\omega}_{\overline{G}_{0}}\left( \mathcal{B}\right) \right) = \overline{\Omega} \left( \mathcal{B}\right).$$
Therefore, there exists $\overline{g} \in \overline{G}^{2}\left(n\right)$ such that the second-order non-holonomic $\overline{G}_{0}-$structure $\overline{\omega}_{\overline{G}_{0}}\left( \mathcal{B}\right)\cdot \overline{g}$ is a second-order non-holonomic integrable prolongation if, and only if, $\overline{\Omega} \left( \mathcal{B}\right)$ is a second-order non-holonomic integrable prolongation. 
So, using Proposition \ref{119} and Proposition \ref{85}, we have the following result:

\begin{proposition}
A Cosserat continuum $\mathcal{B}$ is homogeneous (resp. locally homogeneous) if, and only if, there exists a crystal reference $\overline{Z}_{0}$ such that $\mathcal{B}$ is homogeneous (resp. locally homogeneous) over $\overline{Z}_{0}$.
\end{proposition}
Hence, we our notion of homogeneity of a Cosserat medium $\mathcal{B}$ (which does not depend on a reference crystal) is equivalent to the existence of a configuration $\Phi$ such that $\mathcal{B}$ is homogeneous over the reference crystal $j^{1}_{e_{1},Z_{0}} \Phi^{-1}$ (in terms of the non-holonomic $\mathcal{G}-$structures of second order).

\end{document}